\documentclass{article} 
\usepackage{_styles}
\usepackage{_ref}

\nomenclature[01]{$(Q,R)$}{a bound quiver}
\nomenclature[02]{$\mathcal{L}$}{a complete set of representatives of isomorphism classes of indecomposable representations of $(Q,R)$}
\nomenclature[03]{$\mathcal{I}$}{subset of interval representations in $\mathcal{L}$}
\nomenclature[04]{$\cat{rep}(Q,R)$}{the category of representations of $(Q,R)$}

\title{\mytitle}

\makeatletter
\renewcommand\AB@affilsepx{,\,  \protect\Affilfont}
\makeatother

\author[1]{Yasuaki Hiraoka}
\author[2]{Ken Nakashima}
\author[3]{Ippei Obayashi}
\author[4]{Chenguang Xu}
\affil[1,4]{Kyoto University}
\affil[2]{Shimane University}
\affil[3]{Okayama University}

\date{}

\begin{document}
\maketitle

\begin{abstract}

A fundamental challenge in multiparameter persistent homology is the absence of a complete and discrete invariant.
To address this issue, we propose an enhanced framework that realizes a holistic understanding of a fully commutative quiver's representation via synthesizing interpretations obtained from intervals.
Additionally, it provides a mechanism to tune the balance between approximation resolution and computational complexity.
This framework is evaluated on commutative ladders of both finite-type and infinite-type.
For the former, we discover an efficient method for the indecomposable decomposition leveraging solely one-parameter persistent homology. 
For the latter, we introduce a new invariant that reveals persistence in the second parameter 
by connecting two standard persistence diagrams using interval approximations.
We subsequently present several models for constructing commutative ladder filtrations,
offering fresh insights into random filtrations and 
demonstrating our toolkit's effectiveness in analyzing the topology of materials.
\end{abstract}

{
\centerline{}
\begin{quote}
\end{quote}
}

\tableofcontents

\section{Introduction}
Topological data analysis (TDA) is a rapidly emerging field in applied mathematics that features 
leveraging algebraic topology tools to solve problems in data science \cite{carlsson2009topology}.
Persistent homology, a core component of TDA, examines homology modules derived from a filtration of topological spaces constructed based on a given dataset.
This filtration provides a multi-scale perspective on the underlying structure of the original data.
In this context, the resulting algebraic object can be viewed as a representation of a quiver. 

The successes of one-parameter persistent homology, 
evident in its fruitful applications across diverse domains, such as 
cosmology \cite{xu2019finding}, medical imaging \cite{oyama2019hepatic, chung2009persistence}, and material science \cite{hiraoka2016hierarchical,hong2019medium},
can be attributed to the structure theorem for finitely generated modules over a principal ideal domain,
which lays a solid foundation for utilizing persistence diagrams as a compact descriptor to encode the topological information of a dataset extractable through persistent homology.

Complex topological structures embedded in real-world datasets often require a multiparameter filtration of topological spaces to comprehensively capture and analyze their properties.
This led to the development of multiparameter persistent homology, 
a highly anticipated field promising potential breakthrough in our understanding of complex systems by tracking the evolution of topological features across multiple parameters. 
However, the non-existence of a discrete and complete invariant \cite{carlsson2007theory} in this situation poses a significant challenge,
making it crucial to develop new methodologies to advance this field further and enhance its accessibility to a broader spectrum of researchers and data practitioners.

\begin{figure}[htbp]
\centering
\[
\begin{tikzcd}[row sep=\tikzcdRowSep,column sep=\tikzcdColumnSep]
\substack{1\\ \circ}\arrow[r,leftrightarrow,shift right=0.9ex]
&\substack{2\\ \circ}\arrow[r,leftrightarrow,shift right=0.9ex]
&\substack{3\\ \circ}\arrow[r,leftrightarrow,shift right=0.9ex]
&\substack{\makebox[2em]{\phantom{0.1em}} \\ \cdots}\arrow[r,leftrightarrow,shift right=0.9ex]
&\substack{n\\ \circ}\\
\substack{\circ \\ 1'}\arrow[u]\arrow[r,leftrightarrow,shift left=0.9ex]
&\substack{\circ \\ 2'}\arrow[u]\arrow[r,leftrightarrow,shift left=0.9ex]
&\substack{\circ \\ 3'}\arrow[u]\arrow[r,leftrightarrow,shift left=0.9ex]
&\substack{\cdots \\ \makebox[2em]{}}\arrow[u]\arrow[r,leftrightarrow,shift left=0.9ex]
&\substack{\circ \\ n'}\arrow[u]
\end{tikzcd}
\]
\caption{
A commutative ladder with length $n$. 
The symbol $\leftrightarrow$ means either $\leftarrow$ or $\rightarrow$. 
Orientations in the two rows are assumed to be identical, and each square commutes.
}

\label{fig:cl_illu_1}
\end{figure}
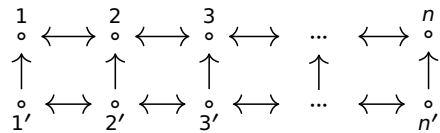

\textbf{Settings.} 
Our research provides a refined theoretical framework for understanding
representations of fully commutative quivers through interval subquivers.
The framework builds upon a generalized version of the boundary compression and interval approximation technique proposed in 
\cite{asashiba2019approximation}.
We validate our new tools on a specific family of fully commutative quivers known as commutative ladders \cite{escolar2016persistence},
characterized by a two-parameter configuration where the second parameter changes only once
(see Figure \ref{fig:cl_illu_1}).

Commutative ladders provide a feasible and valuable testbed for our framework.
When the ladder length of a commutative ladder is equal to or below four, 
it possesses a finite representation type, 
and we can maintain a complete discrete invariant, 
analogous to the scenario in one-parameter persistent homology.
For infinite-type commutative ladders, we introduce a novel invariant based on interval approximations. 
These two approaches enable the exploration of topological structures in datasets that can be fitted to a commutative ladder filtration.



\textbf{Related works.}
The study of multiparameter persistence has made considerable advancements in recent years.
Patel's work generalized persistence diagrams and demonstrated the feasibility of using Möbius inversion to compute them \cite{patel2018generalized}.
The concept of interval approximation, introduced by Asashiba et al. \cite{asashiba2019approximation}, 
serves as a building component for crafting several invariants in multiparameter persistence.
Notably, zigzag persistence and Möbius inversion have been explored in the computation of some important invariants.
For instance, 
Kim et al. proposed the generalized rank invariant in \cite{kim2021generalizedpd},
and Dey et al. demonstrated that the generalized rank invariant can be computed via boundary caps using zigzag persistence \cite{dey2021computing}.
Kim et al. showed how the bigraded Betti numbers can be calculated from the generalized rank invariants \cite{kim2021generalized}.
Botnan et al. proposed a visual representation of the rank invariant in multiparameter persistence modules \cite{botnan2021signed}. 



\textbf{Contributions.}
We present a novel framework for the study of multiparameter persistence modules.
Central to our approach is the introduction of ``tours'' and ``courses'', allowing us to track selected compositions of paths satisfying specific properties.
Building on this, we enrich the established concept of interval approximation as a linear combination of these courses. 
This refined notion of interval approximation offers greater flexibility to extract information from a given interval compared with existing methods.
As we apply the new framework to two-dimensional commutative grids,
a challenge arises: the exponential growth in the number of intervals makes the computation of
interval approximations impractical. 
To address this, we introduce the ``partial interval approximation'', 
an invariant designed to tune the balance between the number of examined intervals
and the resolution of the approximation reached.

We then study commutative ladders using the new framework.
Starting with finite-type cases,
we realize a more streamlined computation of the indecomposable decomposition using only one-parameter zigzag persistence, bypassing the need for 2D representation calculations.
This finding helps reveal several new types of courses and paves the way for exploring non-intervals, an aspect of multiparameter persistent homology previously under-investigated.
Turning our attention to the infinite-type scenario,
we propose a new invariant: the connected persistence diagram. 
It visualizes persistence in both directions by combining two standard persistence diagrams and then connects homology generators according to their vertical persistence measured by an interval approximation.

To complete the picture and facilitate applicability,
we introduce several models tailored for constructing commutative ladder filtrations of simplicial complexes.
These models encompass techniques to create filtrations from point cloud data
and random simplicial complexes.
We then exemplify the versatility of our toolkit through a series of computational demonstrations.
Specifically, we study the topological structures of random simplicial complexes and atomic arrangements using commutative ladder configurations.
Our computational outcomes highlight the effectiveness of the new framework 
and commutative ladders as a tool for studying complex data structures.
Notably, non-interval components exhibit a markedly lower
proportion in configurations derived from point cloud data.


\textbf{Outline.}
This paper is organized as follows.
In \textsection 2, we establish the relevant background and notations used throughout this paper.
In \textsection 3, we refine interval approximations, propose partial interval approximations, and demonstrate how our framework can be used to construct persistence diagrams of a slice in a 2D persistence module as proposed in RIVET \cite{lesnick2015interactive}.
In \textsection 4, we apply our framework to finite-type commutative ladders,
yielding an efficient method for computing any indecomposable decomposition,
and to infinite-type commutative ladders,
yielding a novel diagram to visualize the interval approximations. 
In \textsection 5, we introduce several models for building up commutative ladder filtrations.
In \textsection 6, we employ the new toolkit to analyze the topological properties of filtrations generated from the aforementioned models.
We conclude with an overview of our advancements in \textsection 7.

\textbf{Source code and related data.}
The source code will be available on this paper's homepage \cite{dedicatedWebpage}.

\section{Preliminaries}
This section reviews key concepts in quiver representations and persistent homology and fixes conventions.
Details and proofs can be found in \cite{assem2006elements}, \cite{asashiba2019approximation} and \cite{oudot2017persistence}. 
We adhere to a fixed base field $\Bbbk$ throughout this paper.



\subsection{Representations of Quivers with Relations}

\begin{defn}
\hfill
\begin{itemize}
\item 
A \emph{quiver} $Q=(Q_0,Q_1,\startmap,\targetmap)$ is a directed multigraph consisting of
a vertex set $Q_0$, 
an arrow set $Q_1$ 
and two maps $\startmap, \targetmap\colon Q_1\to Q_0$ 
assigning the start and target for each arrow in $Q_1$. 
This quiver is called \emph{finite} if both $Q_0$ and $Q_1$ are finite sets.
\item A vertex $u\in Q_0$ is referred to as a \emph{source} if it has no arrows pointing toward it, 
and a \emph{sink} if it has no arrows originating from it. 
\item  
A \emph{path} $p$ from vertex $u$ to vertex $v$ 
is a finite sequence of concatenable arrows in $Q_1$, 
written as $p=(v|\alpha_l\cdots\alpha_1|u)$,
where $\startmap(\alpha_1)=u$, $\targetmap(\alpha_i)=\startmap(\alpha_{i+1})$ for $i=1,\ldots,l-1$, and $\targetmap(\alpha_l)=v$.
Following the function composition convention, the arrows are ordered from right to left. 
If the starting vertex and the target vertex are clear from the context, 
they are omitted, and the path is written as  ${\alpha_l\cdots\alpha_1}$.
The maps $\startmap$ and $\targetmap$ can be extended to the set of paths,
by defining $\startmap(p)= \startmap(\alpha_1)$ and 
$\targetmap(p)= \targetmap(\alpha_l)$.
\item
A path's \emph{length} is the number of arrows it contains.
\item
Paths $p_1$ and $p_2$ are \emph{parallel},
denoted by $p_1\spara p_2$,
if they share the same starting and target vertices.
\item 
Every vertex $u\in Q_0$ is associated with a unique length-zero path, called the \emph{trivial path}, denoted as $\varepsilon_u=(u||u)$. 
Each vertex and its trivial path can be regarded as equivalent.
\item 
For a non-negative integer $l$,
define $Q_l$ as the set of paths in $Q$ of length $l$,
and $Q_{\geq n}\coloneqq\bigcup_{l\geq n}Q_l$.
\item 
We can construct a unital associative $\Bbbk$-algebra 
with the underlying vector space being the free $\Bbbk$-module generated by $Q_{\geq 0}$, 
and the product is given by the concatenation of paths. This associative algebra is 
called the \emph{path algebra} $\Bbbk Q$ of $Q$.
Also, from the concatenation of paths, a quiver can be naturally regarded as a category. 
\item A quiver $Q$ is said to be \emph{acyclic} if it does not contain any cycles, \textit{i.e.}, a non-trivial path that starts and ends at the same vertex.
\end{itemize}
\end{defn}

\begin{figure}[htbp]
\[
\begin{tikzcd}[row sep=\tikzcdRowSep,column sep=\tikzcdColumnSep]
\bullet \arrow[r] &\bullet &\bullet \arrow[l]\\
\end{tikzcd}
\]
\caption{
A type $\mathbb{A}_3$ quiver with orientation $\tau_3=(fb)$.
}
\label{fig:A3-quiver-fb}
\end{figure}
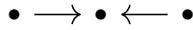

\begin{exa}
A quiver is of \emph{type $\mathbb{A}_n$} 
if its underlying graph is a linear graph with $n$ vertices. 
The orientation of arrows in a type $\mathbb{A}_n$ quiver is specified by a string $\tau_n$ consisting of $n-1$ letters $f$ and $b$, where $f$ stands for a forward arrow and $b$ a backward arrow (see Figure \ref{fig:A3-quiver-fb}).
We represent a type $\mathbb{A}_n$ quiver with its orientation $\tau_n$ as the pair $(A_n,\tau_n)$. 
When the orientation is implicit, we may simply write $A_n$.
An \emph{equi-oriented} type $\mathbb{A}_n$ quiver has all of its arrows pointing in the same direction, 
denoted as $\vv{A}_n$.
\end{exa}

\begin{defn}
Let $Q$ be a quiver.
A \emph{relation} $\rho$ in $Q$ with coefficients in $\Bbbk$ is a linear combination of parallel paths, written as
\[
\rho = \souwa{i=1}{m}\lambda_ip_i,
\]
where $\lambda_i\in \Bbbk$, $p_i\in Q_{\geq 1}$,
and $p_i \spara p_j$ for $i,j\in\Set{1,\ldots, m}$.
For  a set $R\coloneqq\Set{\rho_j}_{j\in J}$ of relations, the pair $(Q,R)$ is called a \emph{quiver with relations}. 
\end{defn}

Given a set of relations $\set{\rho_j}_{j\in J}$, one can generate a two-sided ideal $\mathcal{I} = \gnrt{\rho_j\mid j\in J}$ in the path algebra $\Bbbk Q$, and also define the quotient algebra $\Bbbk Q/\mathcal{I}$.

\begin{rmk}
The relation $\rho=0$ is regarded as the \emph{trivial relation}, 
and it holds that $\mathbb{k}Q\cong \mathbb{k}Q/\gnrt{0}$.
A relation of the form $\rho=\lambda p_1$ with $\lambda\neq 0$ is referred to as a 
\emph{zero relation}.
A relation of the form $p_1-p_2$ is called a 
\emph{commutativity relation}.
\end{rmk}




\begin{defn}
Let $Q$ be a quiver. 
Its \emph{full commutativity relations} is 
a set $R_Q^{\operatorname{fc}}$ defined as:
\[
R_Q^{\operatorname{fc}}=\Set{p_1-p_2 | p_1,p_2 \in Q_{\geq 1} \text{ and } p_1\spara p_2}.
\]
A quiver with these relations is called a \emph{fully commutative quiver}.
\end{defn}





\begin{rmk}
    Throughout this paper, we will focus exclusively on quivers with relations that are finite,  acyclic, and fully commutative unless otherwise specified.
\end{rmk}



\begin{defn}
Let $Q=(Q_0,Q_1,\startmap,\targetmap)$ be a quiver.
\begin{itemize}
    \item A \emph{subquiver} of $Q$ is a quiver $Q'=(Q_0^{'},Q_1^{'},\startmap',\targetmap')$
    such that $Q_0^{'}\subseteq Q_0$, $Q_1'\subseteq Q_1$, $\startmap'$ and $\targetmap'$ are the restrictions of $\startmap$ and $\targetmap$ to $Q_1'$ respectively.
    \item A subquiver $Q'$ is \emph{full} if every arrow $\alpha\in Q_1$ with 
    its start and target in $Q_0^{'}$ also belongs to $ Q_1^{'}$.
    \item A full subquiver $Q'$ of $Q$ is said to be \emph{convex} if, for any path in $Q$ that starts and ends in $Q'$, all intermediate vertices also belong to $Q'$.
    \item The \emph{convex hull} of a set $S \subseteq Q_0$, denoted by $\Conv(S)$, is the full subquiver of $Q$ whose vertices are all the vertices that lie on a path starting and ending with vertices in $S$.
    \item A subquiver is said to be \emph{connected} if its underlying graph is connected.
\end{itemize}
\end{defn}


The notions of convexity and convex hull can be naturally extended to quivers with relations.

\begin{defn}\label{defn:convex-hull-quiver-with-relations}
Let $(Q, R)$ be a quiver with relations and $Q'$ be a subquiver of $Q$. 
The set of induced relations $R'$ on $Q'$ consists of elements 
$\rho\in R$ such that all paths $p_i$ in $\rho=\sum_{i=1}^{m} \lambda_i p_i$
are in $Q'$.
We say $(Q',R')$ is \emph{convex} if for any non-zero path $p$ in $(Q,R)$ with its start and target in $Q'$, 
all vertices of the path are also in $Q'$. 
The convex hull of a set $S\subseteq Q_0$ is then the full subquiver of $Q$ whose 
vertices are all the vertices that lie on a non-zero path in $(Q, R)$ starting and ending with vertices in $S$.
\end{defn}



\begin{defn}
Let $G=(Q, R)$ be a quiver with relations.
An induced subquiver with relations $G'=(Q',R')$ is called an 
\emph{interval subquiver (with relations)}, or simply an \emph{interval} of $G$, if 
$Q'$ is a connected convex subquiver, and $R'$ does not contain any zero relations.
The set of all interval subquivers of $G$ is denoted as $\mathbb{I}_G$. It forms a partially ordered set by containment of the corresponding vertex sets.
Specifically, for two interval subquivers $I$ and $J$, $I\leq J$ if and only if $I_0\subseteq J_0$.
\end{defn}

\begin{rmk}
This paper's definition of an interval differs slightly from \cite[Definition 2.4]{asashiba2019approximation}. 
The cited paper defines an interval for quivers without relations, 
requiring only two conditions: convexity and connectedness.
However, our research engages with quivers with relations, 
and we accommodate quivers with zero relations for generality. 
As a result, 
a new condition ``does not contain any zero relations" is added.
\end{rmk}


 


\begin{exa}
Consider a quiver with relations $Q\coloneqq \clubsuit\xrightarrow{\alpha}\vardiamond\xrightarrow{\beta}\spadesuit$ 
and $R\coloneqq \Set{\beta\alpha}$.
Then, the set of all interval subquivers of $Q$ is
$\mathbb{I}_Q=\set{\clubsuit,\vardiamond,\spadesuit,
\clubsuit\xrightarrow{\alpha}\vardiamond,\vardiamond\xrightarrow{\beta}\spadesuit}$. 
Notice that $\clubsuit\xrightarrow{\alpha}\vardiamond\xrightarrow{\beta}\spadesuit$ does not qualify as an interval subquiver as it encompasses the zero relation $\beta\alpha$.
\end{exa}



Next, we introduce representations of quivers and quivers with relations.
The category of finite-dimensional $\Bbbk$-vector spaces is denoted as $\vect_{\Bbbk}$.
For brevity, $\vect$ is used since we are working over a fixed field.

\begin{defn}
Let $Q$ be a quiver and $G=(Q, R)$ be a quiver with relations.
\begin{itemize}
    \item A (finite-dimensional) \emph{representation} of $Q$ is a functor $M$ from $Q$ to $\vect$. 
    We denote the associated vector spaces as $M_u$ for $u\in Q_0$, and the morphisms as $M(\alpha)$ for $\alpha\in Q_1$.
    \item A \emph{representation} of $G$ is a representation $M$ of the underlying quiver $Q$
satisfying the additional condition that the evaluation of $M$ on each relation $\rho\in R$ vanishes. 
    \item The category of representations of $G$ is denoted as $\rep(G)$.
    \item The vector $\dimv M\coloneqq (\dim_{\Bbbk} M_u)_{u\in Q_0}$ is called the \emph{dimension vector} of $M$
    \item If the associated $\Bbbk$-algebra $\Bbbk Q/\gnrt{R}$ is a representation-finite algebra, $G$ is said to have a \emph{finite type}.
\end{itemize}
\end{defn}

\begin{defn}
    Let $G$ be a quiver with relations, $M$ be a representation in $\rep(G)$, and $G'$ be an induced subquiver with relations.
    A representation in $\rep(G')$ can be derived from $M$ by 
    restricting $M$ to the vertices and arrows of the induced subquiver. 
    This restricted representation is denoted as $\restr{M}{G'}$.
\end{defn}

\begin{defn}
Consider a quiver with relations $G$ and $M$ as a representation of $G$.
The \emph{support} of $M$, denoted as $\supp(M)$, is the full subquiver of $Q$ consisting of vertices $u$ for which $M_u\neq 0$. 
\end{defn}

\begin{defn}\label{defn:interval-and-interval-representation}
Given a quiver with relations $G$ and an interval subquiver $I$, 
the associated interval representation $V_I$ is defined as follows:
\[
\begin{minipage}{.35\linewidth}
    $(V_I)_u\coloneqq\begin{cases}
    \Bbbk & u\in I_0\\
    \mathbf{0} & \text{otherwise};
    \end{cases}$
\end{minipage}
\begin{minipage}{.45\linewidth}
    $V_I(\alpha)\coloneqq\begin{cases}
    \id_{\Bbbk} & \text{both $\startmap(\alpha)$ and $\targetmap(\alpha)$  are in $I_0$}\\
    \mathbf{0} & \text{otherwise}.
    \end{cases}$
\end{minipage}
\]
\end{defn}

\begin{exa}
Consider the following fully commutative quiver and interval subquiver:
\[
\begin{minipage}{.45\linewidth}
\centering
\[
G\coloneqq
\begin{tikzcd}[row sep=\tikzcdRowSep,column sep=\tikzcdColumnSep]
\overset{5}{\bullet} \arrow[r] 
&\overset{6}{\bullet}\arrow[r] 
&\overset{7}{\bullet}\arrow[r] 
&\overset{8}{\bullet}
\\
\underset{1}{\bullet} \arrow[r]\arrow[u]\arrow[ru,phantom,"\circlearrowleft"] 
&\underset{2}{\bullet}\arrow[r]\arrow[u]\arrow[ru,phantom,"\circlearrowleft"]
&\underset{3}{\bullet}\arrow[u]\arrow[r]\arrow[ru,phantom,"\circlearrowleft"] 
&\underset{4}{\bullet}\arrow[u]\\
\end{tikzcd},\]
\end{minipage}
\begin{minipage}{.45\linewidth}
\centering
$
I\coloneqq
\begin{tikzcd}[row sep=\tikzcdRowSep,column sep=\tikzcdColumnSep]
\overset{5}{\bullet} \arrow[r] 
&\overset{6}{\bullet}\arrow[r] 
&\overset{7}{\bullet}
\\
&\underset{2}{\bullet}\arrow[r]\arrow[u]\arrow[ru,phantom,"\circlearrowleft"]
&\underset{3}{\bullet}\arrow[u]\\
\end{tikzcd}$.
\end{minipage}
\]
The associated interval representation $V_I$ is 
$
\begin{tikzcd}[row sep=\tikzcdRowSepSmaller,column sep=\tikzcdColumnSepSmaller]
\Bbbk \arrow[r,"\id",pos=0.4] &\Bbbk\arrow[r,"\id",pos=0.4] & \Bbbk\arrow[r] &\mathbf{0}\\
\mathbf{0} \arrow[r]\arrow[u] &\Bbbk\arrow[r,"\id"',pos=0.4]\arrow[u,"\id"',pos=0.2]& \Bbbk\arrow[u,"\id"',pos=0.2]\arrow[r] &\mathbf{0}\arrow[u]\\
\end{tikzcd}.
$
Notice that a bijective correspondence exists between an interval representation and its dimension vector.
Therefore, we can represent $V_I$ by its dimension vector ${\dimv V_I=
\begin{pmatrix}
1&1&1&0\\
0&1&1&0
\end{pmatrix}}$. 
\end{exa}

\begin{defn}
Let $G$ be a quiver with relations.
A representation $M\in\rep(G)$ is said to be \emph{interval-decomposable} if it is isomorphic to a direct sum of interval representations of $G$.
\end{defn}

When we are working with a finite acyclic quiver $Q$ with a set of relations $R$, both the path algebra $\Bbbk Q$ and its quotient algebra $\Bbbk Q/\gnrt{R}$ have finite dimensions.
As a result, the representation category of $(Q, R)$ satisfies the unique decomposition theorem, also known as the Krull-Schmidt theorem.

\begin{thm}\label{thm:KS}
Let $\mathcal{L}$ be a complete set of representatives of the isomorphism classes of indecomposable representations of a 
quiver with relations $G=(Q, R)$.
For each representation $M\in \rep(G)$, there exists a unique function $d_M\colon \mathcal{L}\to \Z_{\geq 0}$
such that
\begin{equation}\label{eq:KS-thm}
M\cong \Directsum{L\in \mathcal{L}}{}L^{d_M(L)}.
\end{equation}
The function $d_M$ is referred to as the \emph{multiplicity function} of $M$, and the value $d_M(L)$ is called the \emph{multiplicity} of the indecomposable $L$ in $M$.
This isomorphism is referred to as the \emph{indecomposable decomposition} of $M$.
Moreover, $M$ is uniquely determined by $d_M$ up to isomorphism.
\end{thm}


\subsection{Commutative Grids and Commutative Ladders}


\begin{defn}
Let $Q=(Q_0,Q_1,\startmap,\targetmap)$ and $Q'=(Q'_0,Q'_1,\startmap',\targetmap')$ be two quivers. 
Their \emph{Cartesian product} $Q\times Q'$ is a quiver defined as follows:
\begin{itemize}
    \item The vertex set is the Cartesian product $Q_0\times Q'_0$.
    \item There exists an arrow from $(u,u')$ to $(v,v')$ if and only if either:
    \begin{itemize}
        \item $u=v$ and there exists an arrow $\alpha'$ from $u'$ to $v'$, denoted by $(u,\alpha')\in Q_0\times Q'_1$;
        \item $u'=v'$ and there exists an arrow $\alpha$ from $u$ to $v$, denoted by $(\alpha,u')\in Q_1\times Q'_0$.
    \end{itemize}
\end{itemize}
The \emph{tensor product} of $Q$ and $Q'$, denoted by $Q\tensor Q'$, is the quiver $Q\times Q'$ with the following relations:
\[
\Big((v,v')\mid(\alpha,v')(u,\alpha')\mid(u,u')\Big)
-
\Big((v,v')\mid(v,\alpha')(\alpha,u')\mid(u,u')\Big)
\]
for all $\alpha\colon u\to v \in Q_1$ and $\alpha'\colon u'\to v' \in Q'_1$ (see Figure \ref{fig:tensor-illus-1}).
\begin{figure}[htbp]
\[
\begin{tikzcd}[row sep=large,column sep=large]
\substack{(u,v')\\\bullet} \arrow[r,"{(\alpha,v')}"] &
\substack{(v,v')\\\bullet} \\
\substack{\bullet\\(u,u')} \arrow[r,"{(\alpha,u')}"']\arrow[u,"{(u,\alpha')}"] 
&\substack{\bullet\\(v,u')} \arrow[u,"{(v,\alpha')}"'] \\
\end{tikzcd}
\]
\caption{
Illustration of a relation in the tensor product construction.
}
\label{fig:tensor-illus-1}
\end{figure}
\end{defn}

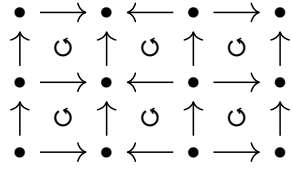
\begin{figure}[htbp]
\[
\begin{tikzcd}[row sep=\tikzcdRowSep,column sep=\tikzcdColumnSep]
\bullet \arrow[r] &\bullet & \bullet\arrow[l]\arrow[r] &\bullet\\
\bullet \arrow[r]\arrow[u] \arrow[ru,phantom,"\circlearrowleft"]&\bullet \arrow[u]\arrow[ru,phantom,"\circlearrowleft"]& \bullet\arrow[l]\arrow[r]\arrow[u] \arrow[ru,phantom,"\circlearrowleft"]&\bullet\arrow[u]\\
\bullet \arrow[r]\arrow[u]\arrow[ru,phantom,"\circlearrowleft"] &\bullet\arrow[u]\arrow[ru,phantom,"\circlearrowleft"]& \bullet\arrow[l]\arrow[u]\arrow[r]\arrow[ru,phantom,"\circlearrowleft"] &\bullet\arrow[u]\\
\end{tikzcd}
\]
\caption{
Commutative grid $G_{4,3}$ with orientation $(\tau_4,\tau_3)$ where $\tau_4=(fbf)$ and $\tau_3=(ff)$.
}
\label{fig:comg-exa1}
\end{figure}

\begin{defn}
    A \emph{two-dimensional fully commutative grid with orientation $(\tau_p,\tau_q)$} 
    is the tensor product of quivers $(A_p,\tau_p)$ and $(A_q,\tau_q)$, 
    written as $G_{p,q}\coloneqq A_p\tensor A_q$ (See Figure \ref{fig:comg-exa1}).
    We often refer to $G_{p,q}$ simply as a \emph{commutative grid} or a \emph{grid}.
    If both $A_p$ and $A_q$ are equi-oriented, 
    we have the \emph{equi-oriented commutative grid},
    represented as $\Gpq\coloneqq\vv{A}_p\tensor \vv{A}_q$.
\end{defn}

The vertices of $G_{p,q}$ can be depicted as a rectangular lattice with $p$ columns and $q$ rows,
with edges connecting vertically or horizontally adjacent vertices.
For equi-oriented grids, we will draw horizontal arrows pointing rightwards and vertical ones upwards without loss of generality.




Intervals in $\vv{G}_{p,q}$ exhibit staircase shapes,
and they can be parameterized as discussed above \cite[Proposition 21]{asashiba2018interval}.
To simplify the notation, we use $\mathbb{I}_{p,q}$ to represent $\mathbb{I}_{\Gpq}$.
\begin{equation}\label{eq:Ipq}
\mathbb{I}_{p,q}=
\Set{
\DUnion{i=s}{t}[b_i,d_i]_i
|
\begin{array}{l}
b_i,d_i,s,t\in \Z,
1\leq b_i\leq d_i\leq p,\ 
1\leq s\leq t\leq q,\ 
\\
\mbox{if $t>s$, then }
1\leq b_{i+1}\leq b_i\leq d_{i+1}\leq d_i \leq p
\ \forall\ i\in \Set{s,\ldots,t-1}
\end{array}
}.
\end{equation}



\begin{defn}
A \emph{commutative ladder} $\clad{\tau_n}$ is a commutative grid $G_{n,2}$ with orientation 
$(\tau_n,\tau_2)$.
We can use just $\tau_n$ to represent its orientation since the value of $\tau_2$ is inconsequential in this context.
An \emph{equi-oriented commutative ladder} of length $n$ is denoted as $\cln$.
\end{defn}

Commutative ladders are of significant interest 
in the transition from one-parameter to two-parameter persistence modules. 
Under specific conditions, 
these ladders exhibit representation-finiteness,
unlocking the possibility of computing all multiplicity functions.
On the other hand,
they remain representation-infinite in general situations, 
driving the need for new approaches.
A criterion for representation-finiteness of commutative ladders is given in \cite{escolar2016persistence}.

\begin{thm}\label{thm:commutative-ladder-rep-type}
For an arbitrary orientation $\tau_n$, the commutative ladder $\clad{\tau_n}$
is
\begin{enumerate}
    \item representation-finite if $n\leq 4$;
    \item representation-infinite if $n\geq 5$.
\end{enumerate}
\end{thm}



\subsection{Persistent Homology}


\begin{defn}
    A \emph{persistence module} refers to a representation of a quiver with relations $G=(Q, R)$.
    This corresponds to a module over the quotient algebra $\Bbbk Q/\gnrt{R}$.
    Several common families have established nomenclature as below.
    \begin{itemize}
    \item A \emph{one-parameter persistence module} is a representation of a  quiver
    $A_n$.
    \item A \emph{zigzag persistence module} designates a one-parameter persistence module, highlighting that the underlying quiver $A_n$ may not be equi-oriented.
    \item A \emph{two-parameter persistence module} is a representation of a commutative grid $G_{p,q}$.
    \item A \emph{persistent homology} is a persistence module obtained by taking a homology functor on a filtration of topological spaces.
    \end{itemize}
\end{defn}

In this paper, the terms ``persistence module'' and ``representation'' are used interchangeably.
Consider a quiver $A_n$ represented as 
    ${\underset{1}{\bullet}\leftrightarrow
        \underset{2}{\bullet}\leftrightarrow
        \cdots
        \leftrightarrow
        \underset{n}{\bullet}}$
    and a persistence module $M$ of $A_n$.
    According to Gabriel's theorem, $M$ is interval-decomposable \cite{gabriel1972unzerlegbare}.
    Then by Theorem \ref{thm:KS},
    $M$ is isomorphic to $\Directsum{I\in\mathbb{I}_{A_n}}V_I^{d_M(V_I)}$.
    
\begin{defn}
Given a one-parameter persistence module $M$ and its associated multiplicity function $d_M$,
the persistence diagram of $M$, denoted as $\PD(M)$, 
visualizes $d_M$ as a multiset of points 
in the two-dimensional integer lattice $\mathbb{Z}^2$.
Here, the multiplicity for $(x,y)$ with $x<y$ is $d_M(V_{[x,y-1]})$.
\end{defn}

\begin{rmk}
    For any point $(b,d)$ in the persistence diagram, we adopt the following conventions:
    \begin{itemize}
        \item The birth coordinate $b$ is inclusive, which means that the generator emerges at value $b$.
        \item The death coordinate $d$ is exclusive, indicating that value $d$ is the earliest point at which 
        the generator vanishes.
    \end{itemize}
\end{rmk}

\section{Refining Interval Approximations}\label{sec:techniques-ia}

In the study of multiparameter persistence modules, 
we aim to obtain the indecomposable decomposition of a given representation $M$, 
equivalent to the computation of the multiplicity function $d_M$.
However, direct computation is extremely challenging,
especially when $M$ is not interval-decomposable. 
To overcome this difficulty, we turn to the interval approximation method, first proposed in \cite{asashiba2019approximation}, which approximates the rank invariant of a representation $M$ via those of interval representations.

A crucial quantity in defining the interval approximation is compression.
It provides a lossy yet more manageable way to define invariants on $M$. 
This section presents a new mechanism for handling various types of compressions in a more general and flexible way.
We then stratify intervals within a general two-dimensional commutative grid, 
proposing the partial interval approximation to address the issue of an exponentially growing number of intervals.
Moreover, we show how the interactive visualization of a 2-D persistence module \cite{lesnick2015interactive} can be reformulated using the interval approximation.
Throughout this section, we consider a finite fully commutative acyclic quiver $G=(Q, R)$.




\subsection{Courses and Tours}
\begin{defn}[Course]
A \emph{course} on $G$ is a pair $(C,F)$
with $C\coloneqq(C_0,C_1,\startmap,\targetmap)$ being a connected quiver
and $F\colon C_0\to Q_0$ acting as a \emph{labeling map},
such that 
for any arrow $\alpha\in C_1$, there exists a path from $F(\startmap(\alpha))$ to $F(\targetmap(\alpha))$ in $G$.
The set of all courses on $G$ is denoted by $\course(G)$.
\end{defn}

\begin{exa}
Consider the fully commutative quiver $\vv{G}_{3,3}$ 
and a connected quiver 
$C\coloneqq\underset{1}{\bullet}\rightarrow\underset{2}{\bullet}\leftarrow\underset{3}{\bullet}$. 
We define two labeling maps $F$ and $F'$ as shown in Figure \ref{fig:compression-exa}.
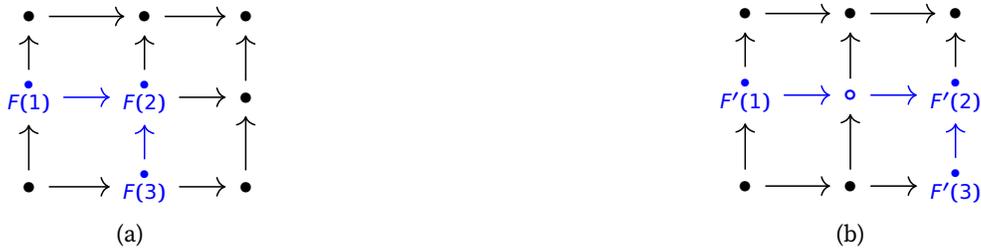
\begin{figure}[htbp]
\begin{subfigure}{0.38\textwidth}
\centering
\[
\begin{tikzcd}[row sep=\tikzcdRowSep,column sep=\tikzcdColumnSep]
\bullet \arrow[r] &\bullet \arrow[r]& \bullet \\
\textcolor{blue}{\substack{\bullet\\F(1)}} \arrow[r,blue]\arrow[u] &\textcolor{blue}{\substack{\bullet\\F(2)}} \arrow[r]\arrow[u]& \bullet\arrow[u] \\
\bullet \arrow[r]\arrow[u] &\textcolor{blue}{\substack{\bullet\\F(3)}} \arrow[r]\arrow[u,blue]& \bullet\arrow[u] \\
\end{tikzcd}
\]
\vspace{-24pt}
\caption{} \label{fig:compression-exa-a}
\end{subfigure}
\hspace*{\fill}
\begin{subfigure}{0.45\textwidth}
\centering
\[
\begin{tikzcd}[row sep=\tikzcdRowSep,column sep=\tikzcdColumnSep]
\bullet \arrow[r] &\bullet \arrow[r]& \bullet \\
\textcolor{blue}{\substack{\bullet\\F'(1)}} \arrow[r,blue]\arrow[u] &\textcolor{blue}{\circ} \arrow[r,blue]\arrow[u]& \textcolor{blue}{\substack{\bullet\\F'(2)}}\arrow[u] \\
\bullet \arrow[r]\arrow[u] &\bullet \arrow[r]\arrow[u]& \textcolor{blue}{\substack{\bullet\\F'(3)}}\arrow[u,blue] \\
\end{tikzcd}
\]
\vspace{-24pt}
\caption{} \label{fig:compression-exa-b}
\end{subfigure}
\vspace{12pt}
\caption{
(a) Labeling map {$F\colon C_0\to\left(\rvv{G}_{3,3}\right)_0$}.
(b) Another labeling map $F'$.
Notice that the vertex in the center is not in the image of $F'$.
}
\label{fig:compression-exa}
\end{figure}
\end{exa}

\begin{figure}[htbp]
\centering
\[
\begin{tikzcd}[row sep=\tikzcdRowSep,column sep=\tikzcdColumnSep]
\bullet \arrow[r] &\circ\arrow[r] & \circ\arrow[r] &\bullet\\
&\bullet\arrow[r]\arrow[u]& \circ \arrow[u]\arrow[r] &\circ\arrow[u]\\
\end{tikzcd}
\]
\caption{A visual representation of essential vertices. 
The interval $I$ comprises all vertices shown, where essential vertices are depicted as solid dots and non-essential vertices are shown as hollow dots.}
\label{fig:ss_1}
\end{figure}
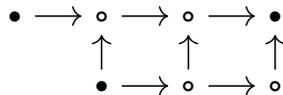

\begin{defn}[Essential Vertex]
Let $I$ be an interval subquiver of $G$. 
A vertex $v\in I_0$ is called
\emph{essential}
if $v$ is either a source or a sink in $I$ (see Figure \ref{fig:ss_1}).
The set of all essential vertices of $I$ is denoted by $E(I)$.
\end{defn}

\begin{rmk}
We note that the essential vertices defined here are called the ``source-sink-essential vertices'' in \cite[Definition 4.1]{asashiba2019approximation}, and they are the minimal information required to recover an interval subquiver. Specifically, an interval subquiver $I$ is the convex hull of $E(I)$.
\end{rmk}

The concept of essential vertices provides a
criterion to determine the containment relations between two intervals in $\mathbb{I}_G$. 
The following proposition justifies their name as being ``essential''. 

\begin{prop}\label{prop:essential-interval-containment}
Let $I,J$ be intervals of $G$.
If $E(I)$ is contained in $J_0$, then $I\leq J$.
\begin{proof}
This statement can be proved similarly as in \cite[Lemma 4.3]{asashiba2019approximation}.
\end{proof}

\end{prop}

Within an interval $I$, courses that visit all $E(I)$ are of particular importance.
\begin{defn}[Essential Course]
Let $(C, F)$ be a course on $G$. 
For an interval $I\in \mathbb{I}_G$, the course $(C,F)$ is said to be:
\begin{itemize}
    \item a \emph{course in $I$} if the image of $F$ is contained in $I_0$.
    \item an \emph{essential course in $I$}, or \emph{essential in $I$}, if it is a course in $I$ and all essential vertices of $I$ are contained in the image of $F$, \textit{i.e.}, $E(I)\subseteq F(C_0)\subseteq I_0$.
\end{itemize}
\end{defn}


\begin{prop}
Consider a course $(C,F)$ on $G$. 
If there exists an interval $I\in\mathbb{I}_G$
for which $(C,F)$ serves as an essential course, 
then $I$ is unique.
\begin{proof}
Suppose $(C, F)$ is an essential course in both $I$ and $J$. 
By definition, we have $E(I)\subseteq F(C_0)\subseteq I_0$ and $E(J)\subseteq F(C_0)\subseteq J_0$,
thus $E(I)\subseteq J_0$ and $E(J)\subseteq I_0$.
Then it follows from Proposition \ref{prop:essential-interval-containment} that $I = J$. 
\end{proof}
\end{prop}

\begin{rmk}
The converse statement is generally not true, 
as multiple essential courses can be defined within a fixed interval. 
Figure \ref{fig:tour_types} provides such an example.
\end{rmk}



\begin{defn}[Essential Assignment]
An \emph{essential assignment} is a map that assigns each interval $I\in \mathbb{I}_G$ an essential course in $I$.
This can be formally expressed as:
\[
\begin{aligned}
\xi\colon \mathbb{I}_G&\to \course(G)\\
I&\mapsto \mbox{$\xi(I)$ is an essential course in $I$}.
\end{aligned}
\]
\end{defn}


\begin{defn}[Tour]
A \emph{tour} on a course $(C,F)$ in $G$
is an additive functor
$\tour_{(C,F)}(-)$ 
that maps a representation $M\in\rep(G)$ 
to an object in $\rep(C)$ specified as below:
\[
\begin{aligned}
\tour_{(C,F)}(-): \rep(G) &\to \rep(C)\\
M&\mapsto \Bigl(M_{F(\startmap(\alpha))}\xrightarrow{M\bigl(F(\startmap(\alpha))\to F(\targetmap(\alpha))\bigr)}M_{F(\targetmap(\alpha))}\Bigr)_{\alpha\in C_1},
\end{aligned}
\]
where $M\bigl(F(\startmap(\alpha))\to F(\targetmap(\alpha))\bigr)$ 
represents the evaluation of $M$ on a path in $G$ from $F(\startmap(\alpha))$ to $F(t(\alpha))$,
which is well-defined by the full commutativity of $G$.
For a morphism $\varphi\colon M\to N$ in $\rep(G)$, a morphism 
$\tour_{(C,F)}(\varphi)$ in $\rep(C)$ is given by
$(M_{F(v)}\xrightarrow{\varphi(F(v))}N_{F(v)})$ for each $v\in C_0$.
\end{defn}

\begin{rmk}
Consider a representation $M\in\rep(G)$.
The choice of a quiver $C$ and a labeling map $F$ can significantly impact the analysis's feasibility.
Using $C=A_n$ will place the tour defined above in the category $\rep(A_n)$,
which usually makes the situation more tractable than working directly with a representation in $\rep(G)$. 
This process, however, can lead to information loss of $M$.
To compensate for it, we employ a set of tours to 
probe $M$,
with each tour offering partial information about $M$ from different perspectives. 
Collectively, they provide a more complete understanding of $M$,
where the amount of information loss varies and is based on the particular choice of courses and the method by which the tours are combined.
For instance, in the context of $\cln$, 
an essential assignment (referred to as a quiver morphism in the cited paper below)
that exhibits several appealing properties is demonstrated in \cite[Section 5.1]{asashiba2023approximation},
highlighting the value of this approach.
\end{rmk}



We introduce the Hasse quiver as an analog to the Hasse diagram in quivers, facilitating 
the characterization of the transitive reduction of paths.

\begin{defn}[Hasse Quiver]
Let $Q=(Q_0, Q_1)$ be an acyclic quiver and let $S\subseteq Q_0$ be a subset of its vertices.
The \emph{Hasse quiver} $\mathfrak{h}(S,Q)$ is a quiver derived from $S$ and $Q$ as follows:
\begin{itemize}
    \item The vertex set of $\mathfrak{h}(S,Q)$ is $S$.
    \item For each pair of distinct vertices $u,v\in S$, an arrow is drawn from $u$ to $v$ if there is a path from $u$ to $v$ in $Q$, and there is no third vertex $w\in S\setminus\Set{u,v}$ in any path from $u$ to $v$ in $Q$.
\end{itemize}
\end{defn}

\begin{exa}
We illustrate the concept of the Hasse quiver by capturing the idea of a ``compressed category'' within the context of $\Gpq$, as described in \cite{asashiba2019approximation}.
Given an interval subquiver $I\in \Ipq$, 
the selection of different subsets $S\subseteq I_0$ 
can lead to a variety of Hasse quivers. 
Consider the following three possibilities:
\begin{enumerate}
    \item $S_1\coloneqq I_0$;
    \item $S_2\coloneqq E(I)$;
    \item $S_3\coloneqq \operatorname{CC}(I)$,
    where $\operatorname{CC}$ is an operation that identifies ``corner-complete'' vertices.
    This subset of $I_0$ includes all vertices present in both a row and a column with essential vertices. 
    Note that all essential vertices are corner-complete.
\end{enumerate}
As an example, consider $\vv{G}_{4,3}$ and an interval $I\in \mathbb{I}_{4,3}$
in Figure \ref{fig:hasse-quiver-example-interval-ss-cc}.
The Hasse quivers on $E(I)$ and $\operatorname{CC}(I)$ are shown in (b) and (c), respectively. Notice that $\mathfrak{h}(S_1, I)$ equals $I$, emphasizing that the Hasse quiver is equal to the original quiver when all vertices are considered.
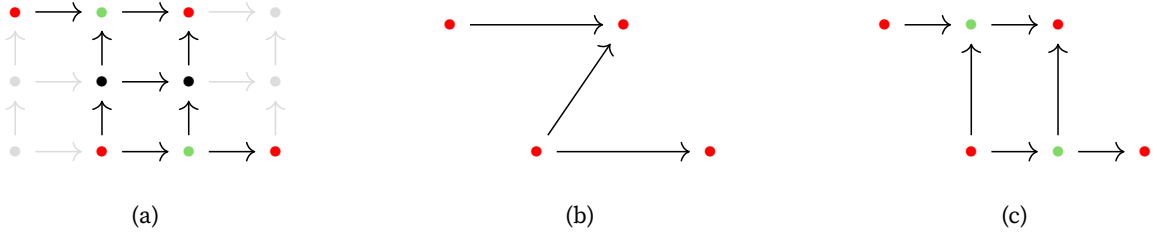
\begin{figure}[htbp]
\begin{subfigure}{0.3\textwidth}
\centering
\[
\begin{tikzcd}[row sep=\tikzcdRowSep,column sep=\tikzcdColumnSep]
\textcolor{red}{\bullet} \arrow[r] & \textcolor{pastel-green}{\bullet} \arrow[r] & \textcolor{red}{\bullet}  \arrow[r,color=dimgray] & \dimgray{\bullet}    \\
\dimgray{\bullet} \arrow[r,color=dimgray]\arrow[u,color=dimgray] & \bullet \arrow[u] \arrow[r] & \bullet \arrow[u] \arrow[r,color=dimgray] & \dimgray{\bullet} \arrow[u,color=dimgray]    \\
\dimgray{\bullet} \arrow[r,color=dimgray]\arrow[u,color=dimgray] & \textcolor{red}{\bullet} \arrow[u] \arrow[r] & \textcolor{pastel-green}{\bullet} \arrow[u] \arrow[r] & \textcolor{red}{\bullet} \arrow[u,color=dimgray] 
\end{tikzcd}
\]
\caption{} \label{fig:hasse-quiver-example-interval}
\end{subfigure}
\hspace*{\fill}
\begin{subfigure}{0.3\textwidth}
\centering
\[
\begin{tikzcd}[row sep=\tikzcdRowSep,column sep=\tikzcdColumnSep]
\textcolor{red}{\bullet} \arrow[rr] & & \textcolor{red}{\bullet}   &    \\
\arrow[u,phantom]& &  &    \\
 \arrow[u,phantom]& \textcolor{red}{\bullet} \arrow[uur] \arrow[rr]&  & \textcolor{red}{\bullet} 
\end{tikzcd}
\]
\caption{}\label{fig:hasse-quiver-example-ss}
\end{subfigure}
\hspace*{\fill}
\begin{subfigure}{0.3\textwidth}
\centering
\[
\begin{tikzcd}[row sep=\tikzcdRowSep,column sep=\tikzcdColumnSep]
\textcolor{red}{\bullet} \arrow[r] & \textcolor{pastel-green}{\bullet} \arrow[r] & \textcolor{red}{\bullet}  &   \\
\arrow[u,phantom]& &  &  \\
 \arrow[u,phantom]& \textcolor{red}{\bullet} \arrow[uu] \arrow[r] & \textcolor{pastel-green}{\bullet} \arrow[uu] \arrow[r] & \textcolor{red}{\bullet} 
\end{tikzcd}
\]
\caption{}\label{fig:hasse-quiver-example-cc}
\end{subfigure}
\caption{
(a)
A visual representation of the interval $I$. 
Vertices and arrows outside the interval $I$ are grayed out for clarity. 
Essential vertices of $I$ are marked red, while corner-complete vertices that are not essential are colored in green.
(b) The Hasse quiver of the essential vertices.
(c) The Hasse quiver of the corner-complete vertices.
}
\label{fig:hasse-quiver-example-interval-ss-cc}
\end{figure}
\end{exa}

We now show two examples demonstrating the adaptability and versatility of the concept of essential assignment.
These examples highlight how it serves as a unified framework for incorporating related definitions from existing literature.

\begin{exa}\label{exa:ess-comp-examples}
The three types of compressions (ss, cc, and tot) introduced in \cite{asashiba2019approximation} 
can be expressed as different essential assignments.
Consider $\Gpq$ and the following maps defined on its set of intervals:
\[
\begin{aligned}
&&\Ipq&\to \course(\Gpq)\\
\xi^{\operatorname{tot}}\colon&& I&\mapsto \bigl(\mathfrak{h}(I_0,I),\id\bigr),\\
\xi^{\operatorname{ss}}\colon&& I&\mapsto \bigl(\mathfrak{h}(E(I),I),\id\bigr),\\
\xi^{\operatorname{cc}}\colon&& I&\mapsto \bigl(\mathfrak{h}(\operatorname{CC}(I),I),\id\bigr).
\end{aligned}
\]
With these constructions, we can verify that
$\tour_{\xi^*(I)}(-)=\operatorname{Comp}^*_I(-)$ 
for $*=\operatorname{ss},\ \operatorname{cc},\ \operatorname{tot}$ as defined in the reference.
\end{exa}



\begin{exa}
We demonstrate how the concept of boundary cap from \cite[Definition 19]{dey2021computing} fits within our framework. 
Consider the interval $I$ depicted in Figure \ref{fig:boundary_cap_example_1}.
Its boundary cap, denoted by $\partial I$, 
can be regarded as a type $\mathbb{A}_n$ quiver.
It is constructed from all vertices in $E(I)$ 
and intermediate vertices on the boundary of $I$ 
to ensure connectivity.
Therefore, $\partial I$ can be expressed as an essential course in $I$.
As a result, we can reproduce it using an essential assignment that follows the same pattern.
\begin{figure}[htbp]
    \centering
    \includegraphics[width=0.332\textwidth]{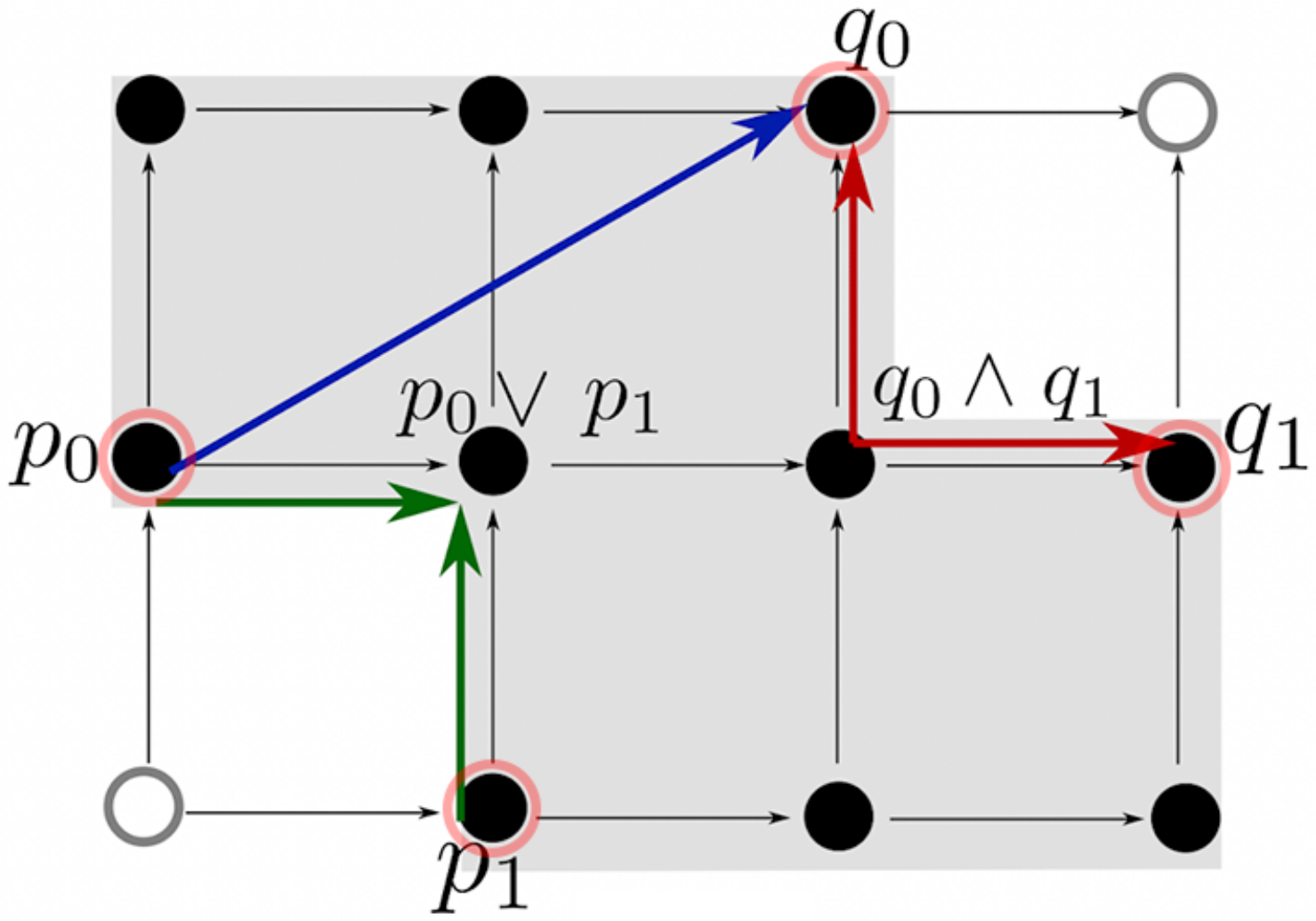}
    \caption{An example of a boundary cap $\partial I$ from \cite[Figure 2]{dey2021computing}. 
    The gray area represents the interval $I$, with essential vertices encircled.
    The sequence of arrows from $p_1$ to $q_1$ represents the boundary cap $\partial I$.
    }
    \label{fig:boundary_cap_example_1}
\end{figure}
\end{exa}


\subsection{\texorpdfstring{$\xi$}{ξ}-compressed Multiplicities and  Interval Approximations}
We are ready to introduce the concept of $\xi$-compressed multiplicity. This quantity captures the multiplicity of an interval representation $V_I$ within a representation $M$ relative to the designated tour $\xi(I)$.

\begin{defn}[$\xi$-compressed Multiplicity]
Let $M$ be a representation of $G$ and $\xi$ be an essential assignment.
For an interval $I\in\IG$,
we define the \emph{$\xi$-compressed multiplicity of $M$ on the interval $I$} as follows:
\[
c^\xi_M(I)\coloneqq d_{\tour_{\xi(I)}(M)}\bigl(
\tour_{\xi(I)}(V_I)
\bigr),
\]
where $d$ is the multiplicity function as in (\ref{eq:KS-thm}).
\end{defn}



\begin{rmk}
Consider an $I\in\IG$ with $\xi(I)=(C,F)$.
Using the definition of $V_I$ and noting that $C$ is a connected quiver,
we can easily verify that $\tour_{\xi(I)}(V_I)$ is an indecomposable representation. 
In particular, if
$C$ is a type $\mathbb{A}_n$ quiver,
then all indecomposable representations of it 
are interval representations.
The longest one among them can be represented as: 
\[
\begin{aligned}
\tour_{(C,F)}(V_I)&\cong
{\underset{1}{\Bbbk}\xleftrightarrow{\id}
\underset{2}{\Bbbk}\xleftrightarrow{\id}
\cdots 
\xleftrightarrow{\id}\underset{n}{\Bbbk}
}
\end{aligned}
\]
with $n=|C_0|$ 
,
which we denote as
$V_{\xi(I)}$
or
$V_{C}$.
\end{rmk}

\begin{exa}
    Consider the essential assignment $\xi^{\operatorname{tot}}$ in Example \ref{exa:ess-comp-examples}.
    Then $c^{\xi^{\operatorname{tot}}}_M(I)=\underline{d}^{\operatorname{tot}}_M(I)$
    in \cite[Remark 4.13]{asashiba2019approximation}, which is shown to be equal to the generalized rank invariant of \cite{kim2021generalizedpd}.
\end{exa}



\begin{prop}\label{prop:int-containment-multiplicity}
Let $\xi$ be an essential assignment. 
For any $I,J\in \mathbb{I}_G$ satisfying $I\leq J$, the following holds:
\[
\tour_{\xi(I)}(V_J)=\tour_{\xi(I)}(V_I).
\]
\begin{proof}
This equality is straightforward, as when we compute the value on the left-hand side, any vertices not included within $I$ can be ignored.
\end{proof}
\end{prop}



\begin{prop}\label{prop:compression-additive}
Let $\xi$ be an essential assignment and $I$ be an interval in $\mathbb{I}_G$.
For any $M,N\in\rep(G)$, we have:
\[
c_{M\oplus N}^{\xi}(I)=c_{M}^{\xi}(I)+c_{N}^{\xi}(I).
\]
\begin{proof}
Starting from the definition of the $\xi$-compressed multiplicity, we obtain
\[
c_{M\oplus N}^{\xi}(I) = d_{\tour_{\xi(I)}(M\oplus N)}\bigl(
\tour_{\xi(I)}(V_I)
\bigr).
\]
Using the additivity of the $\tour$ functor with respect to direct sums 
and the additivity of the multiplicity function $d$ 
with respect to direct sums of representations in the subscript, the expression expands to:
\[
d_{\tour_{\xi(I)}(M)\oplus\tour_{\xi(I)}(N)}\bigl(
\tour_{\xi(I)}(V_I)
\bigr)
=d_{\tour_{\xi(I)}(M)}\bigl(
\tour_{\xi(I)}(V_I)
\bigr)+
d_{\tour_{\xi(I)}(N)}\bigl(
\tour_{\xi(I)}(V_I)
\bigr),
\]
substituting back the definition of $\xi$-compressed multiplicity completes the proof.
\end{proof}
\end{prop}


\begin{prop}\label{prop:int-containment-multiplicity-value}
Let $\xi$ be an essential assignment and $I, J\in\IG$.
The $\xi$-compressed multiplicity function $c^\xi_{V_J}$ evaluates as:
\[
c^\xi_{V_J}(I)=
\begin{cases}
1 & I\leq J\\
0 & \text{otherwise}.
\end{cases}
\]
\begin{proof}
By definition,
$c^\xi_{V_J}(I)=d_{\tour_{\xi(I)}({V_J})}\bigl(
\tour_{\xi(I)}(V_I)
\bigr)$.

\noindent\textbf{Case 1:}
If $I\leq J$, 
then Proposition \ref{prop:int-containment-multiplicity} 
asserts that
$d_{\tour_{\xi(I)}(V_I)}\bigl(
\tour_{\xi(I)}(V_I)
\bigr)=1$.

\noindent\textbf{Case 2:}
Otherwise, 
there exists an essential vertex $u$ of $I$ but not in $J$,
as guaranteed by Proposition \ref{prop:essential-interval-containment}.
This implies that $(V_J)_u=\mathbf{0}$.
Given that $\xi(I)$ is an essential course in $I$, $u$ is visited by $\xi(I)$.
As a result, the associated vector space of $u$ in $\tour_{\xi(I)}(V_J)$ is zero,
but it is nonzero in $\tour_{\xi(I)}(V_I)$.
This leads to the zero multiplicity for $\tour_{\xi(I)}(V_I)$ in $\tour_{\xi(I)}(V_J)$.
\end{proof}
\end{prop}

The following lemma is intended to serve as a counterpart to \cite[Lemma 4.21]{asashiba2019approximation}.
\begin{lemma}\label{lemma:compression-key-lemma}
Let $\xi$ be an essential assignment and $I\in\IG$.
If $M$ is an interval-decomposable representation in $\rep(G)$, 
then the following equation holds:
\[
c_M^{\xi}(I) = \souwa{I\leq J\in \mathbb{I}_G}d_M(V_J).
\]
\begin{proof}
By replacing $M$ in the subscript of $c^{\xi}_M(I)$ with its indecomposable decomposition
$\Directsum{J\in\IG}V_J^{d_M(V_J)}$ 
and using the additivity property from Proposition \ref{prop:compression-additive}, 
we can express the left-hand side as $\souwa{J\in\IG}d_M(V_J)\cdot c_{V_J}^{\xi}(I)$. 
Then, by applying Proposition \ref{prop:int-containment-multiplicity-value},
we obtain the desired result.
\end{proof}
\end{lemma}

\begin{defn}[Cover and Join]\label{defn:cover-and-join}
Consider $\IG$ as a partially ordered set.
The \emph{cover} of an interval $I\in \IG$, 
denoted $\Cov I$, is the set of intervals $J\in \IG$
satisfying $I <J$ and there is no interval $L$ such that $I < L <J$.
The \emph{join} of a subset $S \subseteq \IG$, denoted $\bigvee S$, is the supremum of $S$ provided that it exists.
\end{defn}

When $G=\Gpq$, for each $I\in\IG$, the join operation is well-defined for every subset $S\subseteq \Cov I$, as shown in the discussion above \cite[Example 3.7]{asashiba2019approximation}.
Here, $\bigvee S$ equals the minimum interval containing the union of intervals in $S$. 
For simplicity, we use the convention 
$\bigvee \es =I$ for $\es \subseteq \Cov I$.

\begin{thm}\label{thm:multiplicity-from-c}
Consider an equi-oriented commutative grid $\Gpq$.
Let $\xi$ be an essential assignment and $I\in\Ipq$.
If $M$ is an interval-decomposable representation of $\Gpq$, then
\[
    d_M(V_I)=\souwa{S\subseteq \Cov I}(-1)^{\# S}\cdot c_M^{\xi}(\bigvee S).
\]
\begin{proof}
    This follows directly from the Möbius inversion theorem, as detailed in \cite[Section 5]{asashiba2019approximation}.
\end{proof}
\end{thm}

The definition of interval approximation below is motivated by the Möbius inversion of the formula above, 
where we define it for more general quivers and remove the interval-decomposable condition.

\begin{defn}[Interval Approximation]\label{defn:interval-approximation-general}
Let $M$ be a representation of $G$ and $\xi$ be an essential assignment on $\IG$.
The \emph{interval approximation of $M$ by $\IG$ via $\xi$-compressed multiplicity functions}
is an integer-valued function $\delta_M^{\xi}$
that satisfies
\begin{equation}\label{eq:intv-approx}
c_M^{\xi}(I) = \souwa{I\leq J\in \IG}\delta^\xi_M(J)
\end{equation}
for any $I\in \IG$. When the choice of $\xi$ is clear from the context, we refer to $\delta_M^{\xi}$ as the \emph{interval approximation}.
\end{defn}

\begin{rmk}\label{rmk:interval-approximation-gpq}
A function $\delta^{\xi}_M$ can always be constructed as follows.
First, we define $\delta^{\xi}_M(J)\coloneqq c_M^{\xi}(J)$
for each maximal $J$ in $\IG$.
Then, we iteratively trace down along the cover relations and set
$\delta^{\xi}_M(J)\coloneqq c_M^{\xi}(J)-\souwa{J<I\in\IG}\delta^{\xi}_M(I)$. 
However, if the join operation  is well-defined for any subset $S\subseteq \Cov I$ for all $I\in\IG$,
we can apply the Möbius transform to (\ref{eq:intv-approx}) and obtain
    \[
    \delta^{\xi}_M(I)=\souwa{S\subseteq \Cov I}(-1)^{\# S}\cdot c_M^{\xi}({\bigvee S}).
    \]
\end{rmk}

To conclude this subsection, we prove a theorem establishing that 
the interval approximation accurately recovers the rank function, 
thereby justifying its name as an approximation.

\begin{lemma}
Let $\xi$ be an essential assignment and $p$
be a non-zero path in $G$ starting from vertex $s$ and ending at vertex $t$.
Let $B$ denote the convex hull of vertices $\Set{s,t}$.
If $M$ is a representation of $G$, then the following equation holds:
\[
\souwa{I\in\IG}\delta_M^{\xi}(I)\cdot \rank\bigl(V_I(p)\bigr)
=c^{\xi}_M({B}).
\]
\begin{proof}
For any $I\in\IG$, the rank of the morphism $V_I(p)$ satisfies
\[
\rank\bigl(V_I(p)\bigr)=
\begin{cases}
    1 & V_I(p)=\id \,\Leftrightarrow\, \text{each vertex of $p$ is in $I$}\\
    0 & \text{otherwise}.
\end{cases}
\]
Since the convex hull $B$ is the unique minimum interval subquiver
that contains each vertex of $p$, we can reformulate the left-hand side as
\[
\souwa{I\in\IG}\delta_M^{\xi}(I)\cdot \rank\bigl(V_I(p)\bigr)
=\souwa{\substack{I\in\IG\\ V_I(p)=\id}}\delta_M^{\xi}(I)\\
=\souwa{B\leq I \in \IG}\delta_M^{\xi}(I)\\
=c^{\xi}_M(B),
\]
where the last equality is obtained by Definition \ref{defn:interval-approximation-general}.
\end{proof}
\end{lemma}


\begin{lemma}\label{lemma:rank-invariant-lemma-2}
    Assume the same conditions as in the previous lemma.
    Then the following equality holds:
    \[
        c^{\xi}_M(B)=\rank\bigl(M(p)\bigr).
    \]
    \begin{proof}
        Consider the following restrictions:
        \begin{itemize}
            \item Restrict $M\in\rep(G)$ to $\restr{M}{B}\in\rep(B)$.
            \item Use the symbol $\widetilde{B}$ when we view $B$ as an interval subquiver of $\mathbb{I}_B$.
            \item Restrict $V_B\in \rep(G)$ to $\restr{V_B}{B}\in \rep(B)$, which equals $V_{\widetilde{B}}$
            \item Restrict $\xi\colon \IG\to\course(G)$ 
            to $\widetilde{\xi}\colon \mathbb{I}_B\to \course(B)$.
        \end{itemize}
        The $\xi$-compressed multiplicity on the left-hand side can be reformulated as follows:
        \begin{equation}\label{eq:compression-to-restricted-compression}
        c^{\xi}_M(B)
        =
        d_{\tour_{\xi(B)}(M)}\bigl(\tour_{\xi(B)}(V_{B})\bigr)
        =
        d_{\tour_{\widetilde{\xi}(B)}(\restr{M}{B})}\bigl(\tour_{\widetilde{\xi}(B)}(\restr{V_{B}}{B})\bigr)
        =
        c^{\widetilde{\xi}}_{\restr{M}{B}}(\widetilde{B}).
        \end{equation}
        We will show that the 
        right-hand side
        above is equal to $r\coloneqq\rank\bigl(M(p)\bigr)$.
        It is easy to verify that $V_{\widetilde{B}}$ is an indecomposable projective-injective representation.
        Therefore, the multiplicity of $V_{\widetilde{B}}$ in $\restr{M}{B}$ can be expressed as below.

        \[
        \begin{aligned}
        d_{\restr{M}{B}}(V_{\widetilde{B}}) & =\dim_{\Bbbk} \Hom(V_{\widetilde{B}},\restr{M}{B})
        - \dim_{\Bbbk}\Hom\left(V_{\widetilde{B}}/\operatorname{soc}(V_{\widetilde{B}}),\restr{M}{B}\right) 
        & \text{(\cite[Theorem 3.4]{asashiba2017decomposition} )} \\
        & =\dim_{\Bbbk}M(s)-\Bigl(\dim_{\Bbbk}M(s)-\rank\bigl(M(p)\bigr)\Bigr) 
        & \text{(the discussion below \cite[(4.1)]{asashiba2019approximation})} \\
        & =\rank\bigl(M(p)\bigr) \\
        & = r.
        \end{aligned}
        \]
        As a result, $\restr{M}{B}$ can be written as a direct sum ${\restr{M}{B}=V_{\widetilde{B}}^{\oplus r}\oplus N}$, with $N\in\rep(B)$ not having $V_{\widetilde{B}}$ as a direct summand. 
        By applying the additivity (Proposition \ref{prop:compression-additive}) to this decomposition, 
        we have
        \[
        c^{\widetilde{\xi}}_{\restr{M}{B}}(\widetilde{B})
        =
        c^{\widetilde{\xi}}_{V_{\widetilde{B}}^{\oplus r}\oplus N}(\widetilde{B})
        =
        r\cdot c^{\widetilde{\xi}}_{V_{\widetilde{B}}}(\widetilde{B})
        +c^{\widetilde{\xi}}_{N}(\widetilde{B})
        \geq r.
        \]

        Next we show the reverse inequality $c^{\widetilde{\xi}}_{\restr{M}{B}}(\widetilde{B})\leq r$.
        We break down the proof into the following steps.
        \begin{enumerate}
            \item First we observe that $s$ is the only source of $B$.
            If $s$ were not a source, there would be a vertex $t\neq u\in B_0$ 
            with an arrow from $u$ to $s$.
            By $B$'s defining property,
            there would be a path from $s$ to $t$ passing through $u$,
            leading to a cycle containing $s$, 
            contradicting $G$'s acyclicity.
            The definition of the convex hull immediately implies its uniqueness.
            Similarly, $t$ is the only sink of $B$.
        \item Consider the essential course $\widetilde{\xi}(B)\coloneqq (C,F)$, where $C$ is a connected quiver and $F\colon C_0\to B_0$ is the labeling map. Since this is an essential course in $B$, and $s,t$ are essential vertices of $B$, there exists vertices $c_s,c_t\in C_0$
        such that $F(c_s)=s$ and $F(c_t)=t$.
        As $C$ is connected, we can find a subquiver $W$ of $C$ of type $\mathbb{A}_n$, 
        which starts from $c_s$ and ends at $c_t$.
        This subquiver induces a type $\mathbb{A}_n$ course $\omega\coloneqq(W,\restr{F}{W})$ in $B$, where
        $\restr{F}{W}$ denote the restriction of $F$ to $W_0$. 
        Observe that vertices in $\Im(\restr{F}{W})$ together with arrows between each adjacent pair can be expressed as 
        \begin{equation}\label{eq:main-lemma-proof-walk-labels}
        \underset{s}{\bullet}\rightarrow
        \underset{2}{\bullet}\leftrightarrow
        \underset{3}{\bullet}\leftrightarrow
        \cdots
        \leftrightarrow
        \underset{n-1}{\bullet}\rightarrow
        \underset{t}{\bullet},
        \end{equation}
        where the directions of the first and last arrows are fixed, as shown in the first step above.
        Since $W$ is a subquiver of $C$, the following inequality holds:
        \[
        \begin{aligned}
        c^{\widetilde{\xi}}_{\restr{M}{B}}(\widetilde{B})
        & =
        d_{\tour_{\widetilde{\xi}(B)}(\restr{M}{B})}\bigl(\tour_{\widetilde{\xi}(B)}(V_{\widetilde{B}})\bigr)
        \\
        & \leq 
        d_{\tour_{\omega}(\restr{M}{B})}\bigl(\tour_{\omega}(V_{\widetilde{B}})\bigr)\\
        &\eqqcolon x.\\
        \end{aligned}
        \]
        Therefore, it suffices to show that $x\leq r$.
        By definition, $\tour_{\omega}(V_{\widetilde{B}})^{\oplus x}$ is a summand of $\tour_{\omega}(\restr{M}{B})$,
        thereby we have the following section and retraction, where $\pi\circ\sigma=\id$:
        \[
        \begin{tikzcd}
            \tour_{\omega}(\restr{M}{B})\arrow[d,bend left=20,"\pi"]\\
            \tour_{\omega}(V_{\widetilde{B}})^{\oplus x} \arrow[u,bend left=20,"\sigma"]
        \end{tikzcd}.
        \]
        We can explicitly represent the sections and retractions between vector spaces 
        using the labels in (\ref{eq:main-lemma-proof-walk-labels}) as the two commutative diagrams below,
        where the orientation of horizontal arrows is identical to those in (\ref{eq:main-lemma-proof-walk-labels}): 
        \[
        \begin{tikzcd}[column sep=large]
            M_s\arrow[r]& 
            M_{2}\arrow[r,leftrightarrow]&
            M_{3}\arrow[r,leftrightarrow]&
            \cdots \arrow[r,leftrightarrow]&
            M_{n-1}\arrow[r]&
            M_t\\
            \Bbbk^{x}\arrow[r,"\id"']\arrow[u,"\sigma_s"]& 
            \Bbbk^{x}\arrow[r,leftrightarrow,"\id"']\arrow[u,"\sigma_2"]&
            \Bbbk^{x}\arrow[r,leftrightarrow,"\id"']\arrow[u,"\sigma_3"]&
            \cdots \arrow[r,leftrightarrow,"\id"']&
            \Bbbk^{x}\arrow[r,"\id"']\arrow[u,"\sigma_{n-1}"]&
            \Bbbk^{x}\arrow[u,"\sigma_t"]
        \end{tikzcd},
        \]
        \[
        \begin{tikzcd}[column sep=large]
            M_s\arrow[r]\arrow[d,"\pi_s"]& 
            M_{2}\arrow[r,leftrightarrow]\arrow[d,"\pi_2"]&
            M_{3}\arrow[r,leftrightarrow]\arrow[d,"\pi_3"]&
            \cdots \arrow[r,leftrightarrow]&
            M_{n-1}\arrow[r]\arrow[d,"\pi_{n-1}"]&
            M_t\arrow[d,"\pi_t"]\\
            \Bbbk^{x}\arrow[r,"\id"']& 
            \Bbbk^{x}\arrow[r,leftrightarrow,"\id"']&
            \Bbbk^{x}\arrow[r,leftrightarrow,"\id"']&
            \cdots \arrow[r,leftrightarrow,"\id"']&
            \Bbbk^{x}\arrow[r,"\id"']&
            \Bbbk^{x}
        \end{tikzcd}.
        \]
        Since each morphism in the lower row is bijective, we can reorient all arrows to be forward-going.
        With this adjustment, we formulate the following commutative diagram by 
        selecting sections from $\sigma_s$ to $\sigma_{n-1}$ together with retraction $\pi_t$, 
        and identity maps are also indexed for clarity:
       \[
        \begin{tikzcd}[column sep=large]
            M_s\arrow[r]& 
            M_{2}\arrow[r,leftrightarrow]&
            M_{3}\arrow[r,leftrightarrow]&
            \cdots \arrow[r,leftrightarrow]&
            M_{n-1}\arrow[r]&
            M_t\arrow[d,"\pi_t"]\\
            \Bbbk^{x}\arrow[r,"\id_x({s\rightarrow 2})"' {shift={(0,-0.4ex)}}]\arrow[u,"\sigma_s"]& 
            \Bbbk^{x}\arrow[r,"\id_x({2\rightarrow 3})"' {shift={(0,-0.4ex)}}]\arrow[u,"\sigma_2"]&
            \Bbbk^{x}\arrow[r,"\id_x({3\rightarrow 4})"' {shift={(0,-0.4ex)}}]\arrow[u,"\sigma_3"]&
            \cdots \arrow[r,"\id_x({n-2\rightarrow n-1})"' {shift={(0,-0.4ex)}}]&
            \Bbbk^{x}\arrow[r,"\id_x({n-1\rightarrow t})"' {shift={(0,-0.4ex)}}]\arrow[u,"\sigma_{n-1}"]&
            \Bbbk^{x}
        \end{tikzcd}.
        \]
        Let ${\id_x(s\rightarrow t)}$ denote the composition of all identity maps in the lower row.
        If we can prove ${\pi_t\circ M(p)\circ \sigma_s=\id_x(s\rightarrow t)}$,
        then it would follow that $r=\rank\bigl(M(p)\bigr)\geq \rank\bigl(\id_x(s\rightarrow t)\bigr)=x$.
        \item Now we prove the aforementioned equation in the second step.
        Let $M_1\coloneqq M_s$, $p_1\coloneqq p$, $\sigma_1\coloneqq \sigma_s$ and 
        $M_n\coloneqq M_t$, $p_n\coloneqq \id_{M_t}$, $\sigma_n\coloneqq \sigma_t$
        for easier indexing.
        By the definition of the convex hull $B=\Conv\left(\Set{s,t}\right)$, each vertex in it lies on a path from $s$ to $t$.
        Hence the following commutative diagram exists for ${i=2,3,\ldots,n}$:
        \[
        \begin{tikzcd}[column sep=huge]
            M_{i-1}\arrow[rr,bend left=40,"M(p_{i-1})"]\arrow[r,leftrightarrow,"M(i-1\leftrightarrow i)"']&M_i\arrow[r,bend right=20,"M(p_i)"']&M_t.
        \end{tikzcd}
        \]
        We follow the procedures below to reduce the map $M(p_{i-1})\circ\sigma_{i-1}$ for $i=2,3,\ldots,n-1$.
        \begin{itemize}
            \item If the arrow from $M_{i-1}$ to $M_i$ goes forward, then by the commutative diagram
            \[
            \begin{tikzcd}[column sep=large]
                M_{i-1}\arrow[r,"M(i-1\rightarrow i)"]\arrow[rr,bend left=40,"M(p_{i-1})"]&M_i\arrow[r,bend right=20,"M(p_i)"']&M_t\\
                \Bbbk^x\arrow[u,"\sigma_{i-1}"]\arrow[r,"\id_x(i-1\rightarrow i)"',pos=0.42]&\Bbbk^x\arrow[u,"\sigma_i"']
            \end{tikzcd},
            \]
            we have
            \[
                M(p_{i-1})\circ\sigma_{i-1}
                =M(p_i)\circ M(i-1\rightarrow i)\circ \sigma_{i-1}
                =M(p_i)\circ \sigma_i\circ \id_x(i-1\rightarrow i).
            \]
            \item If the arrow from $M_{i-1}$ to $M_i$ goes backward, then by the commutative diagram
            \[
            \begin{tikzcd}[column sep=large]
                M_{i-1}\arrow[rr,bend left=40,"M(p_{i-1})"]&M_i\arrow[l,"M(i-1\leftarrow i)"']\arrow[r,bend right=20,"M(p_i)"']&M_t\\
                \Bbbk^x\arrow[u,"\sigma_{i-1}"]\arrow[r,"\id_x(i-1\rightarrow i)"',pos=0.42]&\Bbbk^x\arrow[u,"\sigma_i"']
            \end{tikzcd},
            \]
            we have
            \[
                M(p_{i-1})\circ\sigma_{i-1}
                =M(p_{i-1})\circ M(i-1\leftarrow i)\circ \sigma_i\circ \id_x(i-1\rightarrow i)
                =M(p_{i})\circ \sigma_i\circ \id_x(i-1\rightarrow i).
            \]
        \end{itemize}
        
        Applying this reduction iteratively, we obtain
        \[
        \begin{aligned}
            \pi_t\circ M(p)\circ \sigma_s &= \pi_t\circ M(p_1)\circ \sigma_1 \\
            &=\pi_t \circ M(p_2)\circ \sigma_2 \circ \id_x(1\rightarrow 2)\\
            &\vdots \\
            &=\pi_t \circ M(p_{n})\circ \sigma_{n} 
            \circ \id_x(n-1\rightarrow n) 
            \circ \cdots 
            \circ \id_x(2\rightarrow 3)
            \circ \id_x(1\rightarrow 2)
            \\
            &=\pi_t\circ \sigma_t \circ \id_x(s\rightarrow t)\\
            &= \id_x(s\rightarrow t).
        \end{aligned}
        \]
        This concludes the proof.
        \end{enumerate}
    \end{proof}
\end{lemma}


\begin{thm}\label{thm:rank-interval-approximation}
Let $\xi$ be an essential assignment on $G$, $M$ be a representation in $\rep(G)$,
and $p$ be a path in $G$.
Then the following equation holds:
\[
\rank\bigl(M(p)\bigr)=\souwa{I\in\IG}\delta_M^{\xi}(I)\cdot \rank\bigl(V_I(p)\bigr).
\]
\begin{proof}
If $p$ is a zero path, then both sides are equal to zero.
For non-zero paths, the result follows directly from the two preceding lemmas.
\end{proof}
\end{thm}

\subsection{Approximation Series in 2D Grids}


In this subsection, 
we focus on the computation of interval approximations in a two-dimensional commutative grid $\Gpq$. 
Even though the number of interval representations in $\Ipq$ is finite, 
it increases exponentially with the size of $\Gpq$.
To address this problem, we define a series of partial approximations, 
each one necessitating more computational demands but also 
providing incrementally details.
Additionally,
we demonstrate that the interactive visualization of a 2D
persistence module by the barcodes of
1D affine slices  in \cite{lesnick2015interactive}
can be reformulated using an approximation in the series developed.



\begin{thm}\label{thm:num-intervals}\cite[Theorem 31]{asashiba2018interval}
The number of interval representations of
$\Gpq$ 
is given by
\[
    \# \Ipq
    =\souwa{h=1}{q}
    \souwa{w=1}{p}
    \frac{(q-h+1)(p-w+1)}{h+w-1}
    \cdot
    \binom{h+w-1}{h-1}
    \cdot
    \binom{h+w-1}{w-1}.
\]
In particular, 
$\# \Ipq$ is $O(p^{2q})$ 
for $p\geq q$.
\end{thm}


We can stratify intervals by enumerating the number of essential vertices
since an interval is fully determined by its essential vertices.
As a result, the number of essential vertices serves as an indicator of the interval's complexity.
 
\begin{defn}[$k$-essential Interval]
Given a positive integer $k$, 
the set of \emph{$k$-essential intervals} $\mathbb{E}_k$ of $\Gpq$
is a subset of $\Ipq$ that includes all intervals with exactly $k$ essential vertices, expressed as:
\[
\mathbb{E}_k=\Set{I\in\Ipq | \#E(I) = k}.
\]
Furthermore, we use $\mathbb{E}_{\leq k}$ to represent
the disjoint union of all $i$-essential intervals for $i\leq k$: 
\[
\mathbb{E}_{\leq k}\coloneqq \Union{i\leq k}\mathbb{E}_i.
\]
\end{defn}

\begin{exa}\label{exa:Ek-cardinality}
The description and cardinality of $\mathbb{E}_k$ for $k=1,2,3$ are detailed as follows:
\vspace{\baselineskip}
\renewcommand{\arraystretch}{1.5}
\noindent
\begin{tabularx}{\textwidth}{
    @{\hspace{1.5em}}
    >{\raggedright\arraybackslash} 
    X
    @{\quad\hspace{1.5em}}
    >{\raggedright\arraybackslash} 
    X
    @{}
  }
  $\mathbb{E}_1=\Set{I\in\Ipq | \mbox{$I$ is a vertex}}$,
  & 
  $\#\mathbb{E}_1=p\cdot q$,
  \\
  ${\mathbb{E}_2=\Set{I\in\Ipq | \mbox{$I$ is a line segment or rectangle}}}$,
  & 
  $\#\mathbb{E}_2=\binom{p}{2} \cdot q + \binom{q}{2}\cdot p+\binom{p}{2}\cdot\binom{q}{2}$,
  \\
  $\mathbb{E}_3=\Set{I\in\Ipq | 
    \#E(I)=3
    }$, &
    $\#\mathbb{E}_3=\frac{p q}{18}\cdot(p^2-1)\cdot(q^2-1)$.
    \\
\end{tabularx}

Note that vertices are not considered as line segments or rectangles in $\mathbb{E}_2$.
The three terms in the expression for $\#\mathbb{E}_2$
represent the number of horizontal line segments, 
vertical line segments, and rectangles, respectively.
For $\#\mathbb{E}_3$, we can count the ways of cutting a smaller rectangle from a larger non-degenerate rectangle of size subject to the condition that either the top-right or the bottom-left vertices of both rectangles coincide.
This enumeration is formulated by the sum
${\souwa{h=2}{q}\souwa{w=2}{p}2\cdot\underbrace{(p-w+1)\cdot(q-h+1)}_{\text{$\#$larger rectangles of size $w\times h$}}\cdot\underbrace{(w-1)\cdot(h-1)}_{\text{$\#$smaller rectangles}}}$, which then reduces to the expression above.
\end{exa}
Recall the definitions of cover and join from Definition \ref{defn:cover-and-join}.
For a subset $\sIpq$ of $\Ipq$, we define $\Cov_{\sIpq}J$ as the cover of $J$ 
induced by the partial order on $\sIpq$ for each $J\in\sIpq$.
Likewise, given an interval $J\in\sIpq$ and a subset $S\subseteq \Cov_{\sIpq} J$,
we use $\bigvee\limits_{\sIpq} S$ to denote the join of $S$ under the induced partial order, assuming it exists.
We can now generalize the concept of interval approximation from Definition \ref{defn:interval-approximation-general}
by broadening the criteria for intervals to be considered.


\begin{defn}[Partial Interval Approximation]\label{defn:partial-intv-approx}
Consider a representation $M\in\rep(\Gpq)$, a set $\sIpq\subseteq \Ipq$, and an essential assignment $\xi$ defined on $\sIpq$.
A \emph{partial interval approximation} of $M$ 
by $\sIpq$ via $\xi$-compressed multiplicities is defined as an integer-valued function 
$\prescript{}{\sIpq}\delta^{\xi}_M$
defined on intervals of $\sIpq$ that satisfies the following equation for any $I\in\sIpq$:
\[
c^{\xi}_M(I)=\souwa{I\leq \elesIpq\in \sIpq}\prescript{}{\sIpq}\delta^{\xi}_M(\elesIpq).
\]
\end{defn}


\begin{rmk}
In parallel to Remark \ref{rmk:interval-approximation-gpq}, 
$\prescript{}{\sIpq}\delta^{\xi}_M(\elesIpq)$ can be built up 
by firstly setting $\prescript{}{\sIpq}\delta^{\xi}_M(\elesIpq)\coloneqq c_M^{\xi}(\elesIpq)$
for each $\elesIpq$ that is a maximal element of $\sIpq$, then tracing down along the cover relations in iterative steps by setting ${\prescript{}{\sIpq}\delta^{\xi}_M(\elesIpq)\coloneqq c_M^{\xi}(\elesIpq)-\souwa{\elesIpq<I\in\sIpq}\prescript{}{\sIpq}\delta^{\xi}_M(I)}$. 
%
%
Similarly, if 
$\Cov_{\sIpq} \elesIpq$ and $\bigvee\limits_{\sIpq} S$ are well-defined for  any $\elesIpq\in\sIpq$ and $S \subseteq \Cov_{\sIpq} \elesIpq$,
we can use the Möbius inversion to express the partial interval approximation as:
\begin{equation*}
    \prescript{}{\sIpq}\delta^{\xi}_M(\elesIpq)=
    \souwa{S\subseteq\Cov_{\sIpq}\elesIpq} 
    (-1)^{\# S}\cdot c_M^{\xi}({\bigvee\limits_{\sIpq}S}).
\end{equation*}
\end{rmk}


\begin{rmk}\label{rmk:rectangle-approximation}
We refer to the partial interval approximation 
$\prescript{}{\mathbb{E}_{\leq 2}}\delta^{\xi}_M(\elesIpq)$
as the \emph{rectangle approximation}.
\end{rmk}

\begin{defn}[Rank Invariant]
Consider a representation $M$ in $\rep(\Gpq)$ and a subset of intervals ${\sIpq\subseteq \Ipq}$.
A partial interval approximation $\prescript{}{\sIpq}\delta^{\xi}_M$ is said to be \textit{rank invariant} if 
the rank of $M$ over any path $w$ in $\Gpq$ 
can be expressed as the following weighted sum:
\[
\rank\bigl(M(w)\bigr)=
\souwa{\elesIpq\in\sIpq}
\prescript{}{\sIpq}\delta_M^{\xi}(\elesIpq)\cdot \rank\bigl(V_\elesIpq(w)\bigr).
\]
\end{defn}

\begin{exa}
Any interval approximation $\delta^{\xi}_M$ is rank invariant, as shown by Theorem \ref{thm:rank-interval-approximation}.
\end{exa}


The compressed multiplicity function $c^{\xi}_M(\elesIpq)$ 
offers more flexibility compared to the rank function, stemming from the freedom to select $\elesIpq$ and $\xi$.
This flexibility permits the selection of a subset $\sIpq$ of intervals $\Ipq$ for approximation, 
leading to a generalization of the rank invariant property.



\begin{defn}[$k$-rank Invariant]
Let $k$ be a fixed non-negative integer. 
A partial interval approximation 
$\prescript{}{\sIpq}\delta^{\xi}_{M}$ 
is said to be
\emph{$k$-rank invariant} if
 the following equation holds for any $I\in \mathbb{E}_{\leq k+1}$:
\[
c^{\xi}_M(I)=\souwa{J\in \sIpq}\prescript{}{\sIpq}\delta^{\xi}_M(\elesIpq)\cdot c^{\xi}_{V_\elesIpq}(I).
\]
\end{defn}

We now establish the equivalence 
between rank invariance and $1$-rank invariance.
Consider a path ${s\xrightarrow{\alpha_{l}\cdots\alpha_1}t}$ in $\Gpq$ and the rectangle interval $I$ bounded by vertices $s$ and $t$.
By invoking Lemma \ref{lemma:rank-invariant-lemma-2}, we obtain:
\[
c^{\xi}_M(I)=
d_{\tour_{\xi(I)}(M)}\bigl(\tour_{\xi(I)}(V_I)\bigr)
=
\rank\bigl(M({s\xrightarrow{\alpha_{l}\cdots\alpha_1}t})\bigr).
\]
Meanwhile, Proposition \ref{prop:int-containment-multiplicity-value} indicates that $c^{\xi}_{V_\elesIpq}(I)$ is equal to $\rank\bigl(V_J({s\xrightarrow{\alpha_{l}\cdots\alpha_1}t})\bigr)$. 
Substituting these two equations back into the defining equation of being $1$-rank invariant completes the statement.

\begin{rmk}
We observe that $0$-rank invariance is equivalent to preserving the dimension vector of the original representation. 
As discussed above, $1$-rank invariance ensures that the rank of paths is maintained.
Furthermore,
$2$-rank invariance preserves information 
about more complex shapes, such as the L-shaped regions depicted in $\mathbb{E}_3$ of Example \ref{exa:Ek-cardinality}.
\end{rmk}

\begin{thm}
Let $M$ be a representation of $\Gpq$.
If $\mathbb{E}_{\leq k+1} \subseteq \sIpq$, 
then the partial interval approximation 
$\prescript{}{\sIpq}\delta^{\xi}_M$ is $k$-rank invariant.
In particular,
the interval approximation $\delta^{\xi}_M=\prescript{}{\mathbb{E}_{\leq\infty}}\delta^{\xi}_M$
is $k$-rank invariant with respect to 
all non-negative integers $k$.
\begin{proof}
Consider the defining equation
$
c^{\xi}_M(I)=\souwa{I\leq J\in \sIpq}\prescript{}{\sIpq}\delta^{\xi}_M(\elesIpq)
$.
Given that 
$
{c^{\xi}_{V_\elesIpq}(I)=
\begin{cases}
1 & I\leq \elesIpq\\
0 & \text{otherwise}
\end{cases}
}
$
by Proposition \ref{prop:int-containment-multiplicity-value},
we can incorporate a multiplier $c^{\xi}_{V_\elesIpq}(I)$ for each summand, 
and then change the summation range to $\elesIpq\in \sIpq$.
This adjusted summation aligns with the definition of being $k$-rank invariant, concluding the proof.
\end{proof}
\end{thm}

\begin{cor}
If \mbox{ }$\mathbb{E}_{\leq 2}\subseteq \sIpq$, then
$\prescript{}{\sIpq}\delta^{\xi}_M$ is rank invariant.
As a result, the rectangle approximation $\prescript{}{\mathbb{E}_{\leq 2}}\delta^{\xi}_M$
as defined in Remark \ref{rmk:rectangle-approximation} is also rank invariant. 
\end{cor}

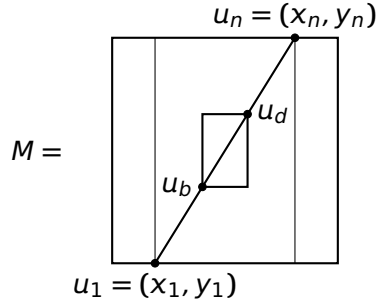
\begin{figure}[htbp]
    \centering
    \begin{tikzpicture}[scale=1]
        \draw (-1,1.5) node {$M=$};
        \draw[thick] (0,0) rectangle (3,3);
        \draw[thick] (0.57,0)--(2.43,3);
        \draw (0.57,0) node[anchor=north] {$u_1=(x_1,y_1)$};
        \filldraw (0.57,0) circle (.3ex);
        \draw (2.43,3) node[anchor=south] {$u_n=(x_n,y_n)$};
        \filldraw (2.43,3) circle (.3ex);
        \draw[thick] (1.2,1.01613) rectangle (1.8,1.98387);
        \draw (1.2,1.01613) node[anchor=east] {$u_b$};
        \filldraw (1.2,1.01613) circle (.3ex);
        \draw (1.8,1.98387) node[anchor=west] {$u_d$};
        \filldraw (1.8,1.98387) circle (.3ex);
        \draw (0.57,0)--(0.57,3);
        \draw (2.43,0)--(2.43,3);
        \draw (4,1.5) node {\phantom{$=M$}};
    \end{tikzpicture}
    \caption{Illustration of a 1D affine slice in a 2D grid.}
    \label{fig:slice-pd-example}
\end{figure}

Here we show that rectangle approximations of a 2D persistence module are equivalent to
1D affine slices as described in RIVET \cite{lesnick2015interactive}.
We use $\prescript{}{s}\rectangle^{\high{t}}$ in $\mathbb{E}_{\leq 2}$ to denote the rectangle interval bounded 
by the two vertices $s$ and $t$ if there exists a path from $s$ to $t$.
Consider the slice along the line connecting $u_1$ and $u_n$ in Figure \ref{fig:slice-pd-example},
and let $A_n$ denote the type $\mathbb{A}_n$ quiver defined by the slice. 
We denote the interval of $A_n$ connecting $s$ and $t$ as $[s,t]$.
The multiplicity of $V_{[u_b,u_d]}$ in the compressed representation $\tour_{\xi^{\operatorname{ss}}(\prescript{}{u_1}\rectangle^{\high{u_n}})}(M)\in \rep(A_n)$ 
can be calculated as:
\[
\begin{aligned}
    &d_{\tour_{\xi^{\operatorname{ss}}(\prescript{}{u_1}\rectangle^{\high{u_n}})}(M)}\bigl(V_{[u_b,u_d]}\bigr)\\
    =&
    \rank\bigl(M(u_b\to u_d)\bigr)
    -\rank\bigl(M(u_{b-1}\to u_d)\bigr)
    -\rank\bigl(M(u_{b}\to u_{d+1})\bigr)
    +\rank\bigl(M(u_{b-1}\to u_{d+1})\bigr)
    \\
    =&
    c_M^{\xi^{\operatorname{ss}}}\bigl({\prescript{}{u_b}\rectangle^{\high{u_d}}}\bigr)
    -c_M^{\xi^{\operatorname{ss}}}\bigl({\prescript{}{u_{b-1}}\rectangle^{\high{u_d}}}\bigr)
    -c_M^{\xi^{\operatorname{ss}}}\bigl({\prescript{}{u_b}\rectangle^{\high{u_{d+1}}}}\bigr)
    +c_M^{\xi^{\operatorname{ss}}}({\prescript{}{u_{b-1}}\rectangle^{\high{u_{d+1}}}}\bigr).
\end{aligned}
\]
The multiplicity of $V_{[u_b,u_d]}$ in the compressed representation $\tour_{\xi^{\operatorname{ss}}(\prescript{}{u_1}\rectangle^{\high{u_n}})}(M)$ 
is equal to the multiplicity of the interval $[u_b,u_d]$ in the persistence diagram of the slice,
thus both the rectangle approximation and RIVET's persistence diagram are determined by 
$\set{c^{\xi^{ss}}_M(I)}_{I\in\mathbb{E}_{\leq 2}}$.


\begin{figure}[htbp]
    \centering
\[
\begin{tikzcd}[row sep=normal,column sep=large]
O(pq)&O(p^2q^2)&O(p^3q^3)&\cdots&O(p^{2q})\\
\mathbb{E}_{\leq 1}\arrow[r,phantom, "\subseteq"] \arrow[d,phantom,"\rotatebox{-90}{$\rightsquigarrow$}"]
& \mathbb{E}_{\leq 2} \arrow[r,phantom, "\subseteq"]  \arrow[d,phantom,"\rotatebox{-90}{$\rightsquigarrow$}"]
& \mathbb{E}_{\leq 3} \arrow[r,phantom, "\subseteq"]  \arrow[d,phantom,"\rotatebox{-90}{$\rightsquigarrow$}"]
& \cdots \arrow[r,phantom, "\subseteq"]  
& \mathbb{E}_{\leq \infty} \arrow[l,phantom, "\subseteq"] \arrow[d,phantom,"\rotatebox{-90}{$\rightsquigarrow$}"]\\
\prescript{}{\mathbb{E}_{\leq 1}}\delta^{\xi}_M \arrow[r,phantom] \arrow[d,phantom,"\rotatebox{90}{$=$}"]
& 
\prescript{}{\mathbb{E}_{\leq 2}}\delta^{\xi}_M \arrow[r,phantom]  \arrow[d,phantom,"\rotatebox{90}{$=$}"]
&
\prescript{}{\mathbb{E}_{\leq 3}}\delta^{\xi}_M 
\arrow[r,phantom]
& \cdots \arrow[r,phantom]  
& \prescript{}{\mathbb{E}_{\leq \infty}}\delta^{\xi}_M \arrow[l,phantom] \arrow[d,phantom,"\rotatebox{90}{$=$}"]
\\
\dimv 
& 
\substack{
\mbox{RIVET's PD}/\\
\mbox{rectangle approximation}
}
& 
&
\cdots
&
\delta^{\xi}_M\\
\mbox{coarser}
\arrow[rrrr]
&&&&
\mbox{finer}
\end{tikzcd}
\]
\caption{Partial interval approximation series constructed from an increasing \mbox{sequence} of $k$-essential intervals.}
\label{fig:approximation-series}
\end{figure}
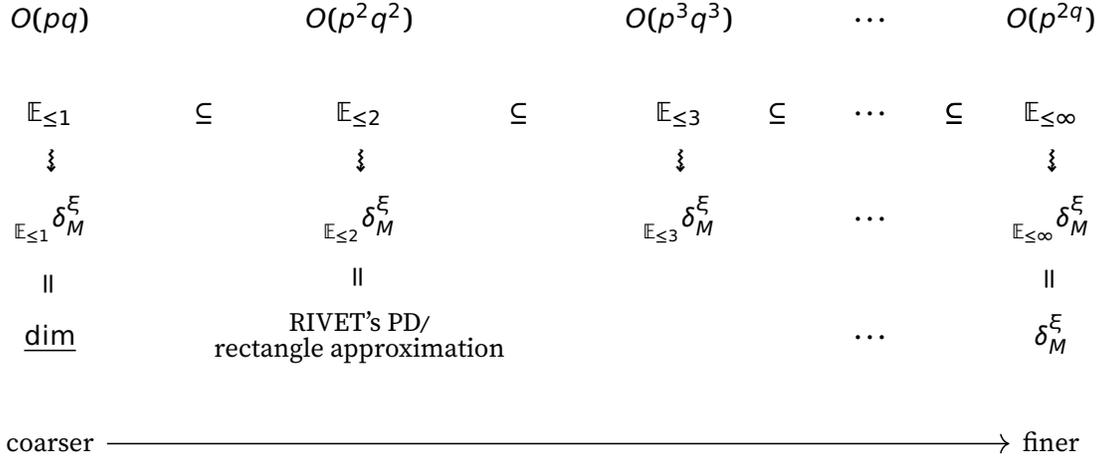

Based on these observations,
we construct an approximation series on $\Gpq$ 
employing the stratification provided by $k$-essential intervals, depicted in Figure \ref{fig:approximation-series}.
The top row displays
the order of the interval count for each corresponding set below.
The third row lists partial interval approximations, 
the fourth row shows names for specific invariants, 
and the bottom row illustrates the change in resolution.

It is worth noting that 
for any fixed $\Gpq$,
this approximation series stabilizes 
because the number of $k$-essential intervals in $\Ipq$ is finite.
The series starts with the dimension vector
and ends with the interval approximation.
RIVET's persistence diagram/rectangle approximation 
is located just to the right of the dimension vector.
As the series progresses, 
the corresponding partial interval approximation
provides more information about the morphisms of the representation $M$ since more intervals are included, 
at the cost of higher computational demands.

\section{Topological Invariants for Commutative Ladders}\label{sec:commutative-ladder}

In this section, we validate the effectiveness of our theoretical framework by addressing specific challenges associated with commutative ladders.
We first devise an algorithm that efficiently computes the indecomposable decomposition of any representation of an equi-oriented finite-type commutative ladder, 
then extend persistence diagrams to infinite-type cases. 
While our focus here is primarily on the equi-oriented cases, both approaches can seamlessly extend to an arbitrary orientation $\tau_n$.
\subsection{Finite-Type Commutative Ladders: Indecomposable Decomposition}\label{subsec:application-cl4}
Consider a finite-type commutative ladder $\cln$.
Let $\mathcal{L}$ be a complete set of representatives of the isomorphism classes of indecomposable $\cln$-modules.
Our objective is to compute the persistent homology
of a $\cln$-filtration $X$, represented as $M\coloneqq H_k(X)$.
This amounts to determine the multiplicity function $d_M(L)$ for every $L\in \mathcal{L}$. 
While existing theoretical frameworks 
such as \cite[Theorem 3.4]{asashiba2017decomposition} and \cite[Section 4]{escolar2016persistence} 
can handle this,
they both require the explicit representation of $M$ 
and the determination of a common basis for the hom-sets between multiple vector spaces, 
which can be computationally intensive.

In contrast, our algorithm circumvents the direct use of $M$ and instead utilizes the filtration $X$.
This approach transforms the original computation into numerous computations of zigzag persistent homology,
a tool where mature and fast algorithms are already available \cite{dey2022fast}.
As a result, we can significantly reduce computational demands and make it more tractable for practical implementations.


Let $G=(Q, R)$ be a fully commutative quiver and $\mathcal{L}$ be a finite subset of a complete set of representatives of the isomorphism classes of indecomposables in $\rep(G)$.
We say a representation ${M\in\rep(G)}$ is $\emph{$\mathcal{L}$-decomposable}$ if every indecomposable direct summand of $M$ is isomorphic to an element in $\mathcal{L}$.
For any $\mathcal{L}$-decomposable representation ${M\in\rep(G)}$ with $M=H_k(X)$, 
Theorem \ref{thm:KS}
guarantees the following decomposition:
\[
M\cong \Directsum{L\in \mathcal{L}}{}L^{d_M(L)}.
\]
Consider a family of functions $\set{f_j}_{j\in J}$ on $\rep(G)$
that are compatible with direct sum operations and isomorphism
(referred to as the \emph{compatibility condition} in the subsequent discussion).
Applying these functions to the decomposition of $M$ yields
the following set of equations:
\[
f_j(M)=\souwa{L\in\mathcal{L}} f_j(L) \cdot d_M(L),\quad \forall\,j\in J.
\]
These equations can be summarized into a matrix expression (assuming the column vector convention):
\begin{equation*}
    \Big(f_j(M)\Big)_{j\in J}
    =
    \Big(f_j(L)\Big)_{\substack{j\in J \\L\in\mathcal{L}}}
    \cdot 
    \Big(d_M(L)\Big)_{L\in \mathcal{L}}.
\end{equation*}
The coefficient matrix has $|J|$ rows and $|\mathcal{L}|$ columns, 
where we can assume $|J|=|\mathcal{L}|$
\footnote{To solve this linear system, 
a coefficient matrix of rank $|\mathcal{L}|$ is required. 
If $|J|<|\mathcal{L}|$, this system is underdetermined, and we need to add more functions.
If $|J|>|\mathcal{L}|$, we can remove surplus functions from $\set{f_j}_{j\in J}$ to equate $|J|$ and $|\mathcal{L}|$.}
.
If its rank is equal to $|\mathcal{L}|$, then there exists a left inverse to it, from which the multiplicity functions can be solved as
\[
\Big(d_M(L)\Big)_{L\in\mathcal{L}}=
\Big(f_j(L)\Big)^{-1}_{\substack{L\in\mathcal{L}\\ j\in J}}\cdot \Big(f_j(M)\Big)_{j\in J}.
\]
Notice that the inverse of the coefficient matrix is independent of the knowledge about $M$. 
For each new filtration $X$ and the associated $H_k(X)$,
we only need to recompute the vector ${\Big(f_j(M)\Big)_{j\in J}=\Big(f_j\circ H_k (X)\Big)_{j\in J}}$. 

\begin{exa}
This example examines the finite-type commutative ladder $\clt$. 
Let $\mathcal{L}$ denote a complete set of representatives of isomorphism classes of $\rep\bigl(\clt\bigr)$, 
where $|\mathcal{L}|=29$ by \cite{escolar2016persistence}.
Denoting the vertex set of $\clt$ as $V\coloneqq\Set{11,21,31,12,22,32}$,
 each representation $M$ can be formulated as:
\[
\begin{tikzcd}
M_{12}\arrow[r]
&M_{22}\arrow[r]
&M_{32}\\
M_{11}\arrow[r]\arrow[u]
&M_{21}\arrow[r]\arrow[u]
&M_{31}\arrow[u]
\end{tikzcd}.
\]
Consider the following functions defined on 
$\rep\bigl(\clt\bigr)$:
\begin{itemize}
    \item $f_{x}(M)\coloneqq d_{M_x}(\Bbbk)=\dim M_{x}$ for $x\in V$;
    \item $f_{x,y}(M)\coloneqq d_{M_x\to M_y}(\Bbbk\xrightarrow{\id}\Bbbk)=\rank \big(M_{x}\to M_{y}\big)$ for $x\neq y\in V$ if there exists a path from $x$ to $y$;
    \item $f_{x,y,z}(M)\coloneqq d_{M_x\to M_y\leftarrow M_z}(\Bbbk\xrightarrow{\id}\Bbbk\xleftarrow{\id}\Bbbk)$ 
    for distinct $x,y,z\in V$ if there exists a course ${(C=\underset{a}{\bullet}\rightarrow\underset{b}{\bullet}\leftarrow\underset{c}{\bullet},F)}$ on $\clt$ such that $F(a)=x,\ F(b)=y,\ F(c)=z$;
    \item $g_{x,y,z}(M)\coloneqq d_{M_x\leftarrow M_y\to M_z}(\Bbbk\xleftarrow{\id}\Bbbk\xrightarrow{\id}\Bbbk)$ 
    for distinct $x,y,z\in V$ if there exists a course ${(C=\underset{a}{\bullet}\leftarrow\underset{b}{\bullet}\rightarrow\underset{c}{\bullet},F)}$ on $\clt$ such that $F(a)=x,\ F(b)=y,\ F(c)=z$;
    \item $f_{x,y,z,w}(M)\coloneqq d_{M_x\to M_y\leftarrow M_z\to M_w}(\Bbbk\xrightarrow{\id}\Bbbk\xleftarrow{\id}\Bbbk\xrightarrow{\id}\Bbbk)$ 
    for distinct $x,y,z,w\in V$ if there exists a course ${(C=\underset{a}{\bullet}\rightarrow\underset{b}{\bullet}\leftarrow\underset{c}{\bullet}\rightarrow\underset{d}{\bullet},F)}$ on $\clt$ such that $F(a)=x,\ F(b)=y,\ F(c)=z,\ F(d)=w$.
\end{itemize}
All these functions meet the compatibility condition.
The linear space spanned by the family of functions of form $f_{x}$ and $f_{x,y}$ on $\mathcal{L}$ has a dimension of 16.
This dimension increases to 26 by including all functions defined by type $\mathbb{A}_3$ courses (\textit{i.e.}, functions of the form $f_{x,y,z}$ and $g_{x,y,z}$), and further to $|\mathcal{L}|$=29 upon adding three linearly independent functions of the form $f_{x,y,z,w}$\footnote{One possible choice is to set $(x,y,z,w)$ as in
$\Set{(12,22,21,31),(11,31,21,22),(21,22,12,32)}$.}.
This yields a coefficient matrix that can be used to compute the indecomposable decomposition of any representation of $\clt$.
\end{exa}



\begin{defn}[Zigzag Course with an Alternating Orientation]
A course $(C,F)$ is called a \emph{zigzag course with an alternating orientation} 
if $C$ is a type $\mathbb{A}_n$ quiver with orientation $\tau=(fbfbfb\cdots)$.
This definition accommodates the general one-parameter cases with any orientation
since consecutive arrows pointing in the same direction can be composed together, and if the first arrow points backward, it can be repeated twice.
For brevity, we use the term ``alternating zigzag course'' 
to refer to such a course.
\end{defn}

\begin{exa}\label{exa:zigzag-course-repeat-arrow}
    Consider a course $(\underset{x}{\bullet}\leftarrow\underset{y}{\bullet}\rightarrow\underset{z}{\bullet}, F)$
    with orientation type $(bf)$. 
    To transform it into an alternating zigzag course starting with a forward arrow, 
    we build a type $\mathbb{A}_4$ course $(\underset{a}{\bullet}\rightarrow\underset{b}{\bullet}\leftarrow\underset{c}{\bullet}\rightarrow\underset{d}{\bullet}, \widetilde{F})$,
    with the labeling map below:
    \[
    \widetilde{F}(a)\coloneqq F(y),\,
    \widetilde{F}(b)\coloneqq F(x),\, 
    \widetilde{F}(c)\coloneqq F(y),\,
    \widetilde{F}(d)\coloneqq F(z).
    \]
    This modified course effectively prepends the map $\underset{y}{\bullet}\rightarrow\underset{x}{\bullet}$ to the original course, and it satisfies the definition of an alternating zigzag course. 
    It can be easily verified that all computations using it in the $\xi$-compressed multiplicity 
    yield the same results as those obtained from the original course based on the definition.
\end{exa}

\begin{rmk}\label{rmk:zigzag-course-from-paths}
Recall that $\startmap$ and $\targetmap$
assign the start and target of a path.
For $n\geq 2$, a sequence of paths $(p_1,\ldots, p_{n-1})$ determines an alternating zigzag course $(C,F)$ of type $\mathbb{A}_n$ 
if it meets one of the following conditions:
\begin{itemize}
    \item If $n=2$, then the path $p_1$ always determines an alternating zigzag course, which is given by 
    ${C=\underset{1}{\bullet}\rightarrow\underset{2}{\bullet}}$ and $F(1)=\startmap(p_1)$, $F(2)=\targetmap(p_1)$.
    \item If $n> 2$, then it determines an alternating zigzag course if and only if the pattern below is satisfied:
    \[
    \underset{1}{\bullet}\xrightarrow{p_1}\underset{2}{\bullet}\xleftarrow{p_2}
    \underset{3}{\bullet}\xrightarrow{p_3}
    \cdots
    \xleftrightarrow{p_{n-1}}\underset{n}{\bullet},
    \]
    where $p_{n-1}$ points forwards if $n$ is even, and backwards if $n$ is odd. This condition can be expressed as:
    \begin{itemize}
        \item $\targetmap(p_i)=\targetmap(p_{i+1})$ for odd $i$ with $1\leq i\leq n-2$;
        \item $\startmap(p_i)=\startmap(p_{i+1})$ for even $i$ with $2\leq i\leq n-2$.
    \end{itemize}
    The alternating zigzag course $(C, F)$ can be easily read from the graphic pattern above.
\end{itemize}

\end{rmk}

Given an alternating zigzag course $(C, F)$, recall that $V_C$ denotes the longest interval representation in $\rep(C)$. 
We associate this course with the following function on $\rep(G)$:
\begin{equation}\label{eq:course-function}
f_{(C,F)}(M)\coloneqq d_{\tour_{(C,F)}(M)}(V_C).
\end{equation}
Notice that it satisfies the compatibility condition.
In particular, although the representation $M$ appears in the subscript, representation $\tour_{(C,F)}(M)=\tour_{(C,F)}\bigl(H_k(X)\bigr)$ is always of type $\mathbb{A}_n$. 
Therefore, we do not need to compute the representation $M$ in advance, significantly speeding up the computation.





Our current objective is to identify a sufficient number of alternating zigzag courses 
that can induce linearly independent functions.
To achieve this, we need to determine the labeling map 
\[
F\colon (A_n)_0\to G_0
\]
via a sequence of (non-trivial) paths $p_1,\ldots,p_{n-1}$ that form an alternating zigzag course $(A_n,F)$ 
for $n\leq N$, where $N$ is a preset limit.
Since $G$ is a finite acyclic quiver,
the count of
alternating zigzag courses (hence labeling maps) is finite
for any fixed positive integer $n$.
This allows for an exhaustive search for all courses.
Appendix \ref{appendix:enum-zigzag-courses} 
provides an illustrative enumeration algorithm
and
Appendix \ref{appendix:bfs-zigzag-courses}
shows a more efficient breadth-first search algorithm.

After gathering all alternating zigzag courses up to length $N$, we associate each course with a function defined by (\ref{eq:course-function}). Notice that multiple courses can result in the same function, and the generated functions can be linearly dependent (for example, the first path in an alternating zigzag course can be repeated twice, similar to the steps described in Example \ref{exa:zigzag-course-repeat-arrow}).
An algorithm of this process is attached in Appendix \ref{appendix:linearly-indep-picking}.

\begin{figure}[htbp]
\centering
\begin{tabular}{ 
>{\centering\arraybackslash} m{3cm} 
>{\centering\arraybackslash} m{2cm}
>{\centering\arraybackslash} m{3cm}
>{\centering\arraybackslash} m{4cm}}
Alternating zigzag course 
& Type & Labeling map & Graphic diagram 
\\\hline\hline
&&& \\[-1em]
$\underset{a}{\bullet}\textcolor{orange-like}{\rightarrow}\underset{b}{\bullet}\textcolor{purple-like}{\leftarrow}\underset{c}{\bullet}$& $\mathbb{A}_3$ &  
$\substack{a\mapsto (1,2),\ b\mapsto (4,2)\\c\mapsto (2,1)\hspace*{3.7em}}$&\includegraphics[width=80pt]{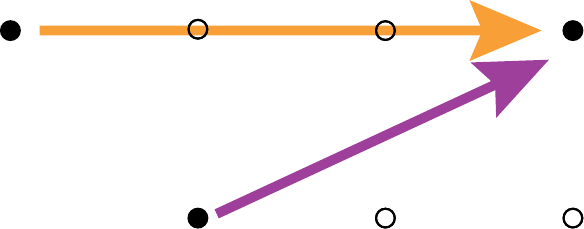}  \\\hline
&&& \\[-1em]
$\underset{a}{\bullet}\textcolor{orange-like}{\rightarrow}
\underset{b}{\circ}\textcolor{purple-like}{\leftarrow}
\underset{c}{\bullet}\textcolor{blue-like}{\rightarrow}
\underset{d}{\bullet}
$& $\mathbb{A}_4$ &  
$\substack{a\mapsto (1,2),\ b\mapsto (2,2)\\c\mapsto (2,1),\ d\mapsto (4,2)}$&\includegraphics[width=80pt]{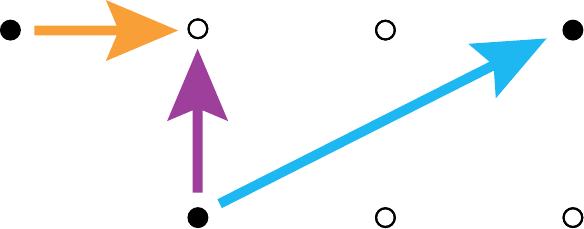}  \\\hline
&&& \\[-1em]
$\underset{a}{\bullet}\textcolor{orange-like}{\rightarrow}
\underset{b}{\circ}\textcolor{purple-like}{\leftarrow}
\underset{c}{\bullet}\textcolor{blue-like}{\rightarrow}
\underset{d}{\bullet}
$& $\mathbb{A}_4$ &  
$\substack{a\mapsto (2,1),\ b\mapsto (3,2)\\c\mapsto (1,2),\ d\mapsto (4,2)}$&\includegraphics[width=80pt]{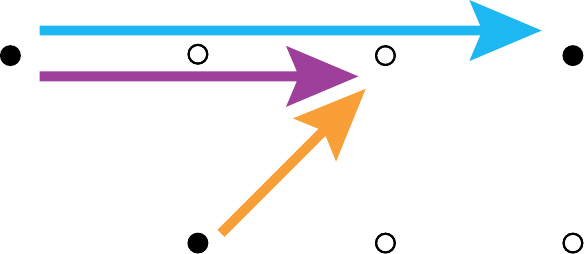}  \\\hline
&&& \\[-1em]
$\underset{a}{\bullet}\textcolor{orange-like}{\rightarrow}
\underset{b}{\bullet}\textcolor{purple-like}{\leftarrow}
\underset{c}{\bullet}\textcolor{blue-like}{\rightarrow}
\underset{d}{\circ}\textcolor{green-like}{\leftarrow}
\underset{e}{\circ}
$& $\mathbb{A}_5$ &  
$\substack{a\mapsto (2,1),\ b\mapsto (4,2)\\c\mapsto (1,2),\ d\mapsto (3,2)\\e\mapsto (3,1)\hspace*{3.7em}}$&\includegraphics[width=80pt]{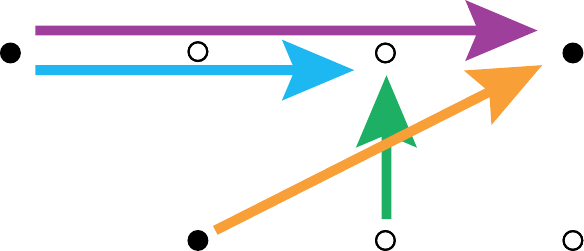}  \\\hline
&&& \\[-1em]
$\underset{a}{\bullet}\textcolor{orange-like}{\rightarrow}
\underset{b}{\circ}\textcolor{purple-like}{\leftarrow}
\underset{c}{\circ}\textcolor{blue-like}{\rightarrow}
\underset{d}{\bullet}\textcolor{green-like}{\leftarrow}
\underset{e}{\bullet}\textcolor{red-like}{\rightarrow}
\underset{f}{\bullet}
$& $\mathbb{A}_6$ &  
$\substack{a\mapsto (2,1),\ b\mapsto (4,1)\\c\mapsto (3,1),\ d\mapsto (3,2)\\ e\mapsto (1,2),\ f\mapsto (4,2)}$&\includegraphics[width=80pt]{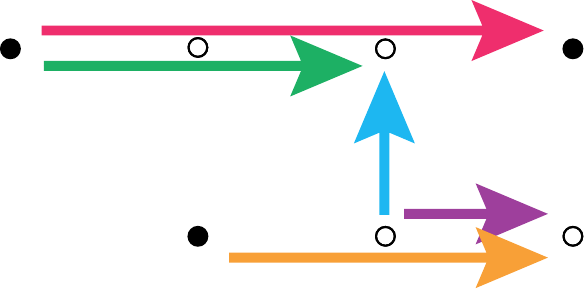} 
\end{tabular}
\caption{Five alternating zigzag courses used for solving the indecomposable decomposition of representations in $\clf$. 
They are essential courses in the interval depicted in Figure \ref{fig:ss_1}, with essential vertices shown as solid dots.}
\label{fig:tour_types}
\end{figure}

We execute the algorithms on $\cln$ for $n=2,3,4$ 
and successfully find sufficient linearly independent functions to solve their indecomposable decompositions.
Specifically, for $\clf$, choosing $N=6$ serves as an adequate preset limit, 
and a detailed list of the obtained 76 alternating zigzag courses can be found in \cite{dedicatedWebpage}.
As an illustration, we show five such courses in Figure \ref{fig:tour_types}.
The vertices in each course are labeled with letters, while those in 
$\clf$ are labeled with their Cartesian coordinates.
Several alternating zigzag courses in this figure
do not belong to the three types of compressions defined in \cite{asashiba2019approximation},
 demonstrating the capability of our framework 
to extract deeper insights compared to the existing methods.


\subsection{Infinite-Type Commutative Ladders: Connected Persistence Diagram}\label{subsec:cln}




Our discussion in the previous subsection relies on 
the $\mathcal{L}$-decomposability condition with a predetermined finite set of isomorphism classes $\mathcal{L}$,
where the finiteness requirement cannot be generalized to representations of a general $\cln$.
To broaden our framework's applicability, 
we introduce a novel invariant crafted specifically for general commutative ladders
by leveraging their unique two-row structure. 
This invariant captures the topological information in the two rows extractable via one-parameter persistence in the horizontal direction and measures the vertical persistence of 
generators along these two rows, using a selected essential assignment.

By setting the parameters in (\ref{eq:Ipq}) to $p=n,q=2$ and categorizing intervals into subsets based on their support, 
we can obtain an indexing for intervals in $\cln$ as below:
\begin{equation}\label{eq:cln-intervals}
\begin{aligned}
\mathbb{I}\coloneqq\mathbb{I}_{n,2} & = 
\Set{[b_1,d_1]_1 | 1\leq b_1\leq d_1\leq n}\\
& \sqcup
\Set{[b_2,d_2]_2 | 1\leq b_2\leq d_2\leq n}
\\ & \sqcup 
\Set{[b_1,d_1]_1\sqcup[b_2,d_2]_2 | 1\leq b_2\leq b_1\leq d_2\leq d_1\leq n}\\
&=\Id\sqcup\Iu\sqcup\Iud,
\end{aligned}
\end{equation}
where $\mathbb{I}_1$ contains all intervals entirely supported in the lower row,
$\mathbb{I}_2$ contains all intervals entirely supported in the upper row,
and
$\mathbb{I}_{2/1}$ contains intervals with support that bridges the two rows.
An element of form $[b_i,d_i]_i$ in $\Id$ or $\Iu$ has type $\mathbb{A}_{d_i-b_i+1}$,
and a typical element in $\Iud$ is shown below:
\[
\begin{tikzcd}[row sep=\tikzcdRowSepSmaller,column sep=\tikzcdColumnSepSmaller]
&\overset{b_2}{\bullet}\arrow[r]&\cdots\arrow[r]&\bullet\arrow[r]&\cdots\arrow[r]&\overset{d_2}{\bullet} \\
&&                               &\underset{b_1}{\bullet}\arrow[r]\arrow[u]&\cdots\arrow[r]&\bullet\arrow[u]\arrow[r]&\cdots\arrow[r]&\underset{d_1}{\bullet}
\end{tikzcd},
\]
where vertices $b_1$ and $b_2$ can share a column, as can $d_1$ and $d_2$.

%

\begin{defn}
Consider an essential assignment $\xi$ on $\mathbb{I}$
and a persistence module $M$ of $\cln$ as:
\begin{equation}\label{eq:cln-filtration}
    \begin{aligned}
    M\coloneqq
    \begin{tikzcd}[row sep=\tikzcdRowSepLarger,column sep=\tikzcdColumnSepLarger]
    M_{(1,2)}\arrow[r]&M_{(2,2)}\arrow[r]&\cdots\arrow[r]&M_{(n,2)}\\
    M_{(1,1)}\arrow[r]\arrow[u]&M_{(2,1)}\arrow[r]\arrow[u]&\cdots\arrow[r]&M_{(n,1)}\arrow[u]
    \end{tikzcd}
    \end{aligned}.
\end{equation}
%
We define an auxiliary function $\widetilde{\delta}^\xi_M \colon \mathbb{I}\to \Z$
constructed from the interval approximation $\delta^\xi_M$, as per the steps delineated below.
\begin{itemize}
    \item For intervals in $\mathbb{I}_1$, with fixed indices $b_1$ and $d_1$, we define
    \[
    \widetilde{\delta}^\xi_M\bigl([b_1,d_1]_1 \bigr) = 
    \delta_M^{\xi}\bigl([b_1,d_1]_1\bigr) +
    \souwa{\substack{[b_2,d_2]_2 \sqcup [b_1,d_1]_1\in \mathbb{I}_{\nicefrac{2}{1}}}}
    \delta_M^{\xi}\bigl([b_2,d_2]_2 \sqcup [b_1,d_1]_1\bigr).
    \]
    \item For intervals in $\mathbb{I}_2$, given $b_2$ and $d_2$, we define
    \[
    \widetilde{\delta}^\xi_M\bigl([b_2,d_2]_2 \bigr) = 
    \delta_M^{\xi}\bigl([b_2,d_2]_2\bigr) +
    \souwa{\substack{ [b_2,d_2]_2 \sqcup [b_1,d_1]_1\in \mathbb{I}_{\nicefrac{2}{1}}}}
    \delta_M^{\xi}\bigl([b_2,d_2]_2 \sqcup [b_1,d_1]_1\bigr).
    \]
    \item For intervals in $\mathbb{I}_{\nicefrac{2}{1}}$, we define
    \[
    \widetilde{\delta}^\xi_M\bigl([b_2,d_2]_2 \sqcup [b_1,d_1]_1\bigr) \coloneqq  \delta_M^{\xi}\bigl([b_2,d_2]_2 \sqcup [b_1,d_1]_1\bigr).
    \]
\end{itemize}
\end{defn}

Notice that the original interval approximation can be retrieved from the associated $\widetilde{\delta}^\xi_M$.
The following proposition  
explains the summations involved in the definition. It shows that $\widetilde{\delta}^\xi_M$ encodes
all the information available in the corresponding multiplicity functions when we only look at the upper or lower row.

\begin{prop}\label{prop:cPD-upper-lower-PD-equivalence}
    Let $\xi$ be an essential assignment and $M$ be a persistence module of $\cln$.
    Referring to (\ref{eq:cln-filtration}), 
    we denote the lower and upper rows of $M$ as
        \[
        \begin{array}{lcr}
            {
            \underline{M} \coloneqq 
                \begin{tikzcd}[column sep=small, ampersand replacement=\&]
                    M_{(1,1)}\arrow[r] \& M_{(2,1)}\arrow[r] \& \cdots\arrow[r] \& M_{(n,1)}
                \end{tikzcd}
            }
            & 
            \text{and} 
            &
            {
            \overline{M} \coloneqq 
                \begin{tikzcd}[column sep=small, ampersand replacement=\&]
                    M_{(1,2)}\arrow[r] \& M_{(2,2)}\arrow[r] \& \cdots\arrow[r] \& M_{(n,2)}
                \end{tikzcd},
            }
        \end{array}
        \]
    respectively. 
    Both $\underline{M}$ and $\overline{M}$ are one-parameter persistence modules 
    over quiver $\vv{A}_n$.
    Define $\mathcal{V}$ as the function mapping intervals in $\Id\sqcup\Iu$
     to representations of $\vv{A}_n$:
    \[
    \begin{aligned}
        \mathcal{V}\colon \Id\sqcup\Iu &\to \rep(\vv{A}_n)\\
        I&\mapsto V_I.
    \end{aligned}
    \]
    Then, the following equations hold:
    \[
    d_{\underline{M}}\circ\mathcal{V} =\widetilde{\delta}^\xi_M\big|_{\Id}
    \qquad \text{and} \qquad
    d_{\overline{M}}\circ\mathcal{V} =\widetilde{\delta}^\xi_M\big|_{\Iu}.
    \]
    In other words, for any interval $I\in\Id\sqcup\Iu$, its multiplicity in $\widetilde{\delta}^\xi_M$ aligns with its multiplicity in the corresponding one-parameter persistence module, irrespective of the choice of $\xi$.
    \begin{proof}
    We here only show the statement for $d_{\overline{M}}\circ\mathcal{V}$.
     By padding the persistence module $M$ with zero vector spaces on both sides 
     (which ensures that ${M_{(x,1)}=\mathbf{0}}$ and ${M_{(x,2)}=\mathbf{0}}$ for $x<1$ or $x>n$),
     rank calculations involving indices outside of the range $x\in [1,n]$ become well-defined.
     Consider an interval $[b,d]_2\in\Iu$, the multiplicity of its associated interval module is
     \[
     \begin{aligned}
     d_{\overline{M}}\circ\mathcal{V}\bigl([b,d]_2\bigr)
     &=d_{\overline{M}}(V_{[b,d]_2})\\
     &=\rank\Bigl(M\bigl((b,2)\to (d,2)\bigr)\Bigr)\\
     &-\rank\Bigl(M\bigl((b-1,2)\to (d,2)\bigr)\Bigr)-\rank\Bigl(M\bigl((b,2)\to (d+1,2)\bigr)\Bigr)\\
     &+\rank\Bigl(M\bigl((b-1,2)\to (d+1,2)\bigr)\Bigr).
     \end{aligned}
     \]
     By Theorem \ref{thm:rank-interval-approximation}, this equation rearranges to:
     \[
     \begin{aligned}
      &\souwa{I\in\mathbb{I}}\delta_M^{\xi}(I)\cdot \rank\Bigl(V_I\bigl((b,2)\to (d,2)\bigr)\Bigr)\\
         -&\souwa{I\in\mathbb{I}}\delta_M^{\xi}(I)\cdot \rank\Bigl(V_I\bigl((b-1,2)\to (d,2)\bigr)\Bigr)
         -\souwa{I\in\mathbb{I}}\delta_M^{\xi}(I)\cdot \rank\Bigl(V_I\bigl((b,2)\to (d+1,2)\bigr)\Bigr)\\
         +&\souwa{I\in\mathbb{I}}\delta_M^{\xi}(I)\cdot \rank\Bigl(V_I\bigl((b-1,2)\to (d+1,2)\bigr)\Bigr).\\
     \end{aligned}
     \]
     We define an alternating sum function $\Omega_{[b,d]_2}$ on the set of intervals $\mathbb{I}$:
     \[
     \begin{aligned}
     \Omega_{[b,d]_2}(I)\coloneqq &
     \rank\Bigl(V_I\bigl((b,2)\to (d,2)\Bigr)\\
     -&\rank\Bigl(V_I\bigl((b-1,2)\to (d,2)\bigr)\Bigr)
     -\rank\Bigl(V_I\bigl((b,2)\to (d+1,2)\bigr)\Bigr)\\
     +&\rank\Bigl(V_I\bigl((b-1,2)\to (d+1,2)\bigr)\Bigr).
     \end{aligned}
     \]
     It is easy to verify the following relations:
     \[
     \Omega_{[b,d]_2}(I)
     =
     \begin{cases}
     1 & \text{if $I=[b,d]_2$}\\
     1 & \text{if $I$ is of form $[b,d]_2\sqcup [b_1,d_1]$}\\
     0 & \text{otherwise}.
     \end{cases}
     \]
    The original equation can be further reduced to:
     \[
     \begin{aligned}
        &d_{\overline{M}}\circ\mathcal{V}\bigl([b,d]_2\bigr)\\
        = & \souwa{I \in\mathbb{I}}\delta_M^{\xi}(I) \cdot
         \Omega_{[b,d]_2}(I)\\
         =& \souwa{\substack{I \in\mathbb{I}\\\Omega_{[b,d]_2}(I)=1 }}
         \delta_M^{\xi}(I) \cdot
         \Omega_{[b,d]_2}(I)\\\\
         =& \delta_M^{\xi}\bigl({[b,d]_2}\bigr)+
            \souwa{I=[b,d]_2 \sqcup [b_1,d_1]_1\in \mathbb{I}_{\nicefrac{2}{1}}}
            \delta_M^{\xi}(I) \\
            =& \widetilde{\delta}^\xi_M\bigl([b,d]_2 \bigr).
     \end{aligned}
     \]
    \end{proof} 
\end{prop}


To visualize $\widetilde{\delta}^\xi_M$,
we first adopt an approach that retains a clear relationship to standard persistence diagrams 
by plotting within two complementary isosceles right triangles that together form a square.
We introduce the following points and line segments in $\mathbb{Z}^2$ for upcoming discussions:
\begin{itemize}
    \item ${\mathbb{T}_1\coloneqq\Set{(d,b)\in\mathbb{Z}_+\times\mathbb{Z}_+| d>b}}$, the lower triangular region in the first quadrant;
    \item ${\mathbb{T}_2\coloneqq\Set{(b,d)\in\mathbb{Z}_+\times\mathbb{Z}_+| b<d}}$, the upper triangular region in the first quadrant;
    \item $\mathbb{T}_{\layerV}\coloneqq\Set{\big((b_2,d_2),(d_1,b_1)\big)\in\mathbb{T}_2\times \mathbb{T}_1|-1\leq \frac{d_2-b_1}{b_2-d_1}\leq 0}$, 
    the set of line segments connecting points from $\mathbb{T}_2$ and $\mathbb{T}_1$ with slope\footnote{This condition on the slope of the line segment results from the index condition of intervals in $\Iud$ specified in (\ref{eq:cln-intervals}).} within the range $[-1,0]$.
\end{itemize}

\begin{defn}[Connected Persistence Diagram]
    Let $\xi$ be an essential assignment and $M$ be a persistence module of $\cln$.
    A \emph{connected persistence diagram} of $M$, denoted as $\cPD^{\xi}(M)$, 
    visualizes the interval approximation $\delta^{\xi}_M$ via $\widetilde{\delta}^\xi_M$ in the two-dimensional integer lattice $\mathbb{Z}^2$ as a multiset comprising both points and line segments,
    where the multiplicity of an element can be negative.
    Represented as a function, elements with a non-trivial multiplicity are given by:
    \begin{equation*}
    \begin{aligned}
    \cPD^{\xi}(M)\colon \mathbb{T}_1\sqcup\mathbb{T}_2\sqcup\mathbb{T}_{\layerV} &\to \mathbb{Z}\\
    \mathbb{T}_1\ni (d,b)&\mapsto  \widetilde{\delta}^\xi_M\big|_{\Id}\big([b,d-1]_1\big)\\
    \mathbb{T}_2 \ni (b,d)&\mapsto  \widetilde{\delta}^\xi_M\big|_{\Iu}\big([b,d-1]_2\big)\\
    \mathbb{T}_{\layerV}\ni \big((b_2,d_2),(d_1,b_1)\big)&\mapsto  \widetilde{\delta}^\xi_M\big([b_2,d_2-1]_2\sqcup[b_1,d_1-1]_1\big)
    \end{aligned}.
    \end{equation*}
\end{defn}

\begin{figure}[htbp]
    \centering
    \includegraphics[width=0.52\textwidth]{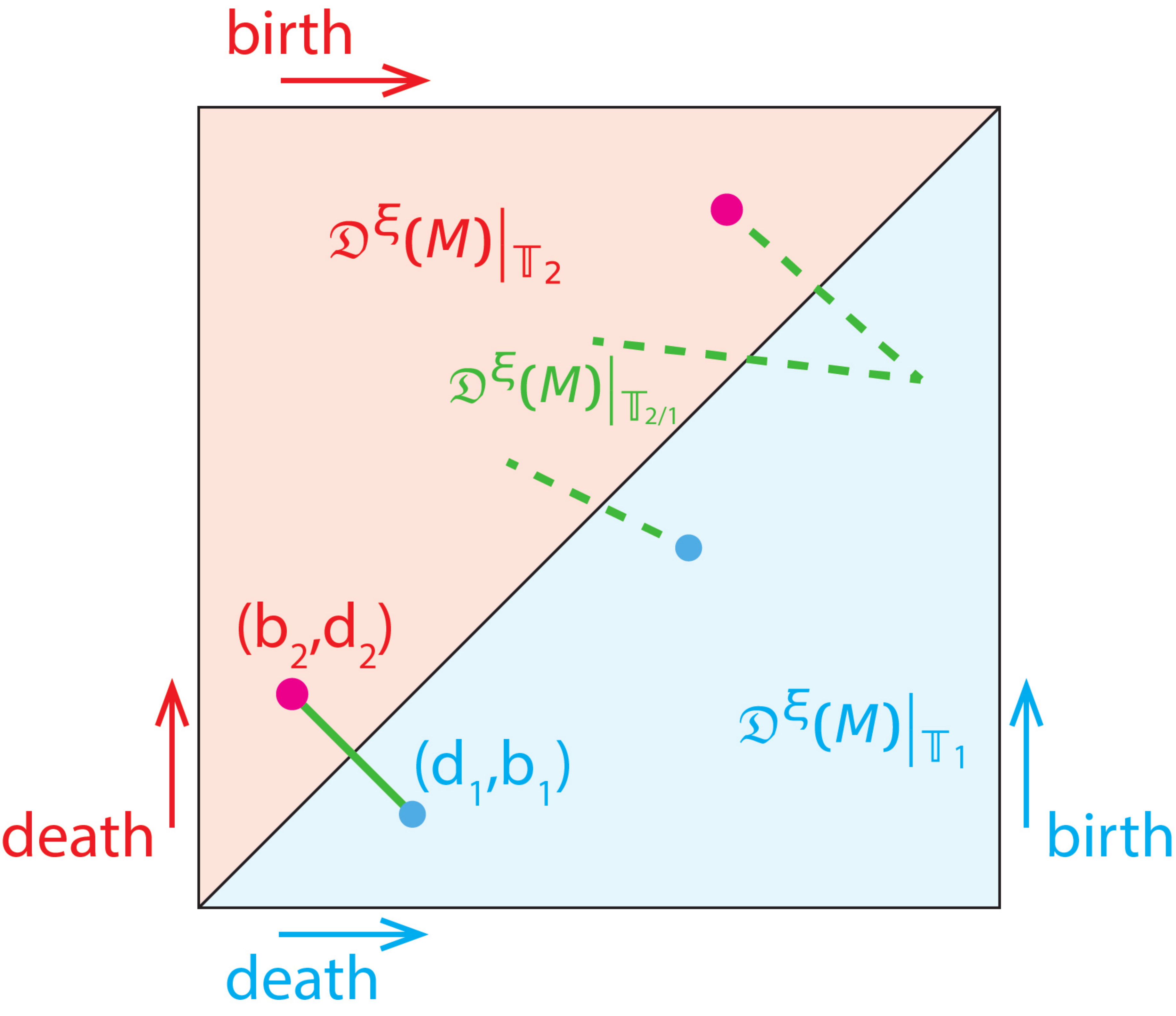}
    \caption{Schematic illustration of a connected persistence diagram. 
    A dashed line indicates negative multiplicity.
    }
    \label{fig:cPD-illustration-0}
\end{figure}


\begin{rmk}
In a connected persistence diagram, 
multiplicities on $\mathbb{T}_1$ and $\mathbb{T}_2$ measure the horizontal persistence,
and multiplicities on $\mathbb{T}_{\layerV}$ can measure the vertical persistence, which is crucial for revealing information that is not accessible via only one-parameter persistent homology on the two rows. 
Notice that the two endpoints of a line segment in $\mathbb{T}_{\layerV}$ often coincide with points from both $\mathbb{T}_1$ and $\mathbb{T}_2$. 
However, as illustrated in Figure \ref{fig:cPD-illustration-0}, this overlap is not requisite.
\end{rmk}

\begin{rmk}
We choose the source-sink essential assignment $\xi^{\operatorname{ss}}$ for our subsequent computations since its simplicity makes it an ideal starting point. 
We note that other essential assignments could offer more refined or specific insights. 
A visualization of $\xi^{\operatorname{ss}}$ on representatives of intervals 
classified by the number of essential vertices is provided in Figure \ref{fig:ss-essential-assignment-four-types}.
\begin{figure}[htbp]
    \begin{minipage}{0.23\linewidth}
    \centering
    \[
    \begin{tikzcd}[row sep=\tikzcdRowSep,column sep=\tikzcdColumnSep]
    \circ&\bullet\\
    \bullet\arrow[ur,color=LightBlack]&\circ
    \end{tikzcd}
    \]
    \end{minipage}
    \begin{minipage}{0.23\linewidth}
    \centering
    \[
    \begin{tikzcd}[row sep=\tikzcdRowSep,column sep=\tikzcdColumnSep]
    \bullet\arrow[r,-,color=LightBlack]&\circ\arrow[r,color=LightBlack]&\bullet\\
    &\bullet\arrow[ur]&\circ
    \end{tikzcd}
    \]
    \end{minipage}
    \begin{minipage}{0.23\linewidth}
    \centering
    \[
    \begin{tikzcd}[row sep=\tikzcdRowSep,column sep=\tikzcdColumnSep]
    \circ&\bullet&\\
    \bullet\arrow[ur]\arrow[r,-,color=LightBlack]&\circ\arrow[r,color=LightBlack]&\bullet
    \end{tikzcd}
    \]
    \end{minipage}
    \begin{minipage}{0.23\linewidth}
    \centering
    \[
    \begin{tikzcd}[row sep=\tikzcdRowSep,column sep=\tikzcdColumnSep]
    \bullet\arrow[r,-,color=LightBlack]&\circ\arrow[r,color=LightBlack]&\bullet&\\
    &\bullet\arrow[ur,color=LightBlack]\arrow[r,-,color=LightBlack]&\circ\arrow[r,color=LightBlack]&\bullet
    \end{tikzcd}
    \]
    \end{minipage}
    \caption{Illustration of intervals $I\in\Iud$ and the courses defined by $\xi^{\operatorname{ss}}(I)$. 
    Solid dots represent essential vertices, with the number of essential vertices ranging from 2 to 4 in each example.
    }
    \label{fig:ss-essential-assignment-four-types}
\end{figure}
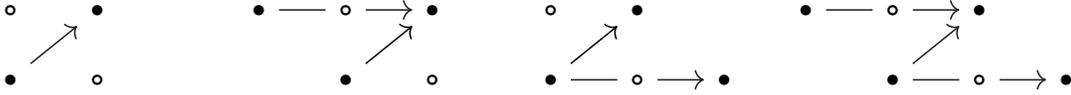
\end{rmk}

\begin{exa}\label{example:cPD-SC-example}
To highlight the effectiveness of connected persistence diagrams, 
we present two configurations with topological structures 
that are distinguishable by connected persistence diagrams
but not by standard persistence diagrams at any homology dimension. 
Figure \ref{fig:cPD-example-ab} depicts each configuration as a filtration of Čech complexes in the horizontal direction.
For clarity, we omit the associated disks and display all simplicial complexes at four critical radii. 
We also directly use the radius as our parameter instead of indexing them with natural numbers.
The symbol $x^-$ denotes a value infinitesimally less than $x$ for $x\in\mathbb{R}$.

Let the two representations be $M_a\coloneqq H_1(X_a)$ and $M_b\coloneqq H_1(X_b)$,
and the connected persistence diagrams be 
$D_a\coloneqq\cPD^{\xi^{\operatorname{ss}}}(M_a)$
and $D_b\coloneqq\cPD^{\xi^{\operatorname{ss}}}(M_b)$.
In both cases, the only interval in $\Iud$ that might yield a non-zero value under the map of $D_a$ or $D_b$
is $I_0\coloneqq[\sqrt{3},2^-]_2\sqcup[\sqrt{3},2^-]_1$.
It is easy to verify that $\Cov I_0=\Set{[(\sqrt{3})^-,2^-]_2\sqcup[\sqrt{3},2^-]_1,[\sqrt{3},2^-]_2\sqcup[\sqrt{3},2]_1}$.
Applying the definition of connected persistence diagram yields
\begin{equation*}
\begin{aligned}
D_c(I_0)=&\delta_{M_c}^{\xi^{\operatorname{ss}}}(I_0)\\
=&\souwa{S\subseteq \Cov I_0}(-1)^{\#S}c^{\xi^{\operatorname{ss}}}_{M_c}(\bigvee S)\\
=&c^{\xi^{\operatorname{ss}}}_{M_c}(I_0)
-c^{\xi^{\operatorname{ss}}}_{M_c}\bigl([(\sqrt{3})^-,2^-]_2\sqcup[\sqrt{3},2^-]_1\bigr)\\
&-c^{\xi^{\operatorname{ss}}}_{M_c}\bigl([\sqrt{3},2^-]_2\sqcup[\sqrt{3},2]_1\bigr)
+c^{\xi^{\operatorname{ss}}}_{M_c}\bigl([(\sqrt{3})^-,2^-]_2\sqcup[\sqrt{3},2]_1\bigr).
\end{aligned}
\end{equation*}
The last three terms all vanish in both cases,
because either ${(M_c)_{((\sqrt{3})^-,2)}=\mathbf{0}}$ or ${(M_c)_{(2,1)}=\mathbf{0}}$.
Continuing with the computation, we find that
$D_c(I_0)=c^{\xi^{\operatorname{ss}}}_{M_c}(I_0)=d_{\tour_{\xi^{\operatorname{ss}}(I_0)}(M_c)}\big(\tour_{\xi^{\operatorname{ss}}(I_0)}(V_{I_0})\big)=d_{(M_c)_{(\sqrt{3},1)}\rightarrow (M_c)_{(2^-,2)}}\big(\Bbbk\xrightarrow{\id}\Bbbk\big)$.
This yields a value of $1$ when $c=a$ and $0$ when $c=b$.



\begin{figure}[htbp]
	\centering
        \vspace*{14pt}
	\begin{subfigure}[c]{0.55\linewidth}
            \begin{tikzpicture}[vertex/.style={circle, fill=node-color, minimum size=1.5mm, inner sep=0pt},thick,scale=0.65,every node/.style={scale=0.9}]
                \def\shiftX{4};
                \def\shiftY{-3};
                \def\lineWidth{0.37};
                \coordinate (A) at (0,0);
                \coordinate (B) at (1.732,0);
                \coordinate (C) at (0.866,1.5);
                \coordinate (D) at (0.866,0.5);
                \coordinate (E) at (-0.866,1.5);
                \coordinate (F) at (0, 1);

                \node at (0, 2.2) {$r = 0$};
                \node at (\shiftX, 2.2) {$r = 1$};
                \node at (2*\shiftX, 2.2) {$r = \sqrt{3}$};
                \node at (3*\shiftX, 2.2) {$r = 2$};

                \foreach \x in {0, 1, 2} {
                    \foreach \y in {0, 1} {
                          \pgfmathsetmacro{\arrowstart}{2+ \x * \shiftX}
                          \pgfmathsetmacro{\arrowend}{2.6+ \x  * \shiftX}
                          \pgfmathsetmacro{\arrowheight}{0.8 + \y * \shiftY}
                          \draw [line width=0.4mm, ->, arrow-color] (\arrowstart,\arrowheight) -- (\arrowend,\arrowheight);
                    }
                }

                \foreach \x in {0,1,2,3}{
                    \pgfmathsetmacro{\arrowlocationX}{0.45+\x * \shiftX}
                    \pgfmathsetmacro{\arrowstart}{1.9+\shiftY}
                    \pgfmathsetmacro{\arrowend}{2.5+\shiftY}
                    \draw [line width=0.4mm, ->, arrow-color] (\arrowlocationX,\arrowstart) -- (\arrowlocationX,\arrowend);
                }

                \node[vertex] at (A) (A12) {};
                \node[vertex] at (B) (B12) {};
                \node[vertex] at (C) (C12) {};
                \node[vertex] at (E) (E12) {};
                \node[vertex] at (F) (F12) {};

                
                \node[vertex] at ([shift={(\shiftX,0)}] A) (A22) {};
                \node[vertex] at ([shift={(\shiftX,0)}] B) (B22) {};
                \node[vertex] at ([shift={(\shiftX,0)}] C) (C22) {};
                \node[vertex] at ([shift={(\shiftX,0)}] E) (E22) {};
                \node[vertex] at ([shift={(\shiftX,0)}] F) (F22) {};
                \draw[draw=edge-color,line width=\lineWidth mm] (E22) -- (F22);
                \draw[draw=edge-color,line width=\lineWidth mm] (C22) -- (F22);
                \draw[draw=edge-color,line width=\lineWidth mm] (A22) -- (F22);

                \fill[fill=area-color] ([shift={(2*\shiftX,0)}] E) -- ([shift={(2*\shiftX,0)}] C) -- ([shift={(2*\shiftX,0)}] A) -- cycle;
                \node[vertex] at ([shift={(2*\shiftX,0)}] A) (A32) {};
                \node[vertex] at ([shift={(2*\shiftX,0)}] B) (B32) {};
                \node[vertex] at ([shift={(2*\shiftX,0)}] C) (C32) {};
                \node[vertex] at ([shift={(2*\shiftX,0)}] E) (E32) {};
                \node[vertex] at ([shift={(2*\shiftX,0)}] F) (F32) {};
                \draw[draw=edge-color,line width=\lineWidth mm] (F32) -- (E32) -- (C32) -- (F32) -- (A32) -- (C32) -- (B32) -- (A32) --(E32);

                \fill[fill=area-color] ([shift={(3*\shiftX,0)}] E) -- ([shift={(3*\shiftX,0)}] C) -- ([shift={(3*\shiftX,0)}] A) -- cycle;
                \fill[fill=area-color] ([shift={(3*\shiftX,0)}] A) -- ([shift={(3*\shiftX,0)}] B) -- ([shift={(3*\shiftX,0)}] C) -- cycle;
                \node[vertex] at ([shift={(3*\shiftX,0)}] A) (A42) {};
                \node[vertex] at ([shift={(3*\shiftX,0)}] B) (B42) {};
                \node[vertex] at ([shift={(3*\shiftX,0)}] C) (C42) {};
                \node[vertex] at ([shift={(3*\shiftX,0)}] E) (E42) {};
                \node[vertex] at ([shift={(3*\shiftX,0)}] F) (F42) {};
                \draw[draw=edge-color,line width=\lineWidth mm] (F42) -- (E42) -- (C42) -- (F42) -- (A42) -- (C42) -- (B42) -- (A42) --(E42);

                \node[vertex] at ([shift={(0,\shiftY)}] A) (A11) {};
                \node[vertex] at ([shift={(0,\shiftY)}] B) (B11) {};
                \node[vertex] at ([shift={(0,\shiftY)}] C) (C11) {};

                \node[vertex] at ([shift={(\shiftX,\shiftY)}] A) (A21) {};
                \node[vertex] at ([shift={(\shiftX,\shiftY)}] B) (B21) {};
                \node[vertex] at ([shift={(\shiftX,\shiftY)}] C) (C21) {};

                \node[vertex] at ([shift={(2*\shiftX,\shiftY)}] A) (A31) {};
                \node[vertex] at ([shift={(2*\shiftX,\shiftY)}] B) (B31) {};
                \node[vertex] at ([shift={(2*\shiftX,\shiftY)}] C) (C31) {};
                \draw[draw=edge-color,line width=\lineWidth mm] (A31) -- (B31) -- (C31) -- (A31);

                \fill[fill=area-color] ([shift={(3*\shiftX,\shiftY)}] A) -- ([shift={(3*\shiftX,\shiftY)}] B) -- ([shift={(3*\shiftX,\shiftY)}] C) -- cycle;
                \node[vertex] at ([shift={(3*\shiftX,\shiftY)}] A) (A41) {};
                \node[vertex] at ([shift={(3*\shiftX,\shiftY)}] B) (B41) {};
                \node[vertex] at ([shift={(3*\shiftX,\shiftY)}] C) (C41) {};
                \draw[draw=edge-color,line width=\lineWidth mm] (A41) -- (B41) -- (C41) -- (A41);
                
            \end{tikzpicture}
            \vspace*{1pt}
		\caption*{(a-1) Filtration $X_a$. }
	\end{subfigure}
	\hspace*{16pt}
        \hfill
	\begin{subfigure}[c]{0.35\linewidth}
 \begin{tikzpicture}
	\draw[->] (-0.3,0) -- (2.5,0) node[right] {$r$};
	\draw[->] (0,-0.3) -- (0,2.5) node[above] {$r$};
	
	\foreach \x [count=\i] in {1.5,1.732,2}{
            \pgfmathsetmacro{\coord}{4*(\x-1.5)}
		\draw (\coord,-0.2) -- (\coord,0.2);
            \ifnum\i=2
                  \node[below] at (\coord, -0.2) {$\sqrt{3}$};
                \else
                  \node[below] at (\coord, -0.2) {\x};
                \fi
	}
	
	\foreach \y [count=\i] in {1.5,1.732,2}{
            \pgfmathsetmacro{\coord}{4*(\y-1.5)}
		\draw (-0.2,\coord) -- (0.2,\coord);
                \ifnum\i=2
                  \node[left] at (-0.2,\coord) {$\sqrt{3}$};
                \else
                  \node[left] at (-0.2,\coord) {\y};
                \fi
	}
	
	\draw[gray, dashed] (0,0) -- (2,2);

        \pgfmathsetmacro{\sqrtthree}{4*(1.732-1.5)}
        
        \draw[black] (\sqrtthree,2) -- (2,\sqrtthree);
	
 \filldraw (\sqrtthree,2) circle (2pt);
	\node[above right] at (\sqrtthree,2) {$(\text{birth}=\sqrt{3},\text{death}=2)$};
	
	\filldraw (2,\sqrtthree) circle (2pt);
	\node[below right] at (2,\sqrtthree) {$(\text{death}=2,\text{birth}=\sqrt{3})$};
\end{tikzpicture}
            \vspace*{0.5pt}
		\caption*{(a-2) Connected persistence diagram $D_a$.}
	\end{subfigure}
	\vskip\baselineskip
	\vspace*{14pt}
	\begin{subfigure}[c]{0.55\linewidth}
            \begin{tikzpicture}[vertex/.style={circle, fill=node-color, minimum size=1.5mm, inner sep=0pt},thick,scale=0.65,every node/.style={scale=0.9}]
                \def\shiftX{4};
                \def\shiftY{-3};
                \def\lineWidth{0.37};
                \coordinate (A) at (0,0);
                \coordinate (B) at (1.732,0);
                \coordinate (C) at (0.866,1.5);
                \coordinate (D) at (0.866,0.5);
                \coordinate (E) at (-0.866,1.5);
                \coordinate (F) at (0, 1);

                \node at (0, 2.2) {$r = 0$};
                \node at (\shiftX, 2.2) {$r = 1$};
                \node at (2*\shiftX, 2.2) {$r = \sqrt{3}$};
                \node at (3*\shiftX, 2.2) {$r = 2$};

                \foreach \x in {0, 1, 2} {
                    \foreach \y in {0, 1} {
                          \pgfmathsetmacro{\arrowstart}{2+ \x * \shiftX}
                          \pgfmathsetmacro{\arrowend}{2.6+ \x  * \shiftX}
                          \pgfmathsetmacro{\arrowheight}{0.8 + \y * \shiftY}
                          \draw [line width=0.4mm, ->, arrow-color] (\arrowstart,\arrowheight) -- (\arrowend,\arrowheight);
                    }
                }

                \foreach \x in {0,1,2,3}{
                    \pgfmathsetmacro{\arrowlocationX}{0.45+\x * \shiftX}
                    \pgfmathsetmacro{\arrowstart}{1.9+\shiftY}
                    \pgfmathsetmacro{\arrowend}{2.5+\shiftY}
                    \draw [line width=0.4mm, ->, arrow-color] (\arrowlocationX,\arrowstart) -- (\arrowlocationX,\arrowend);
                }
                    
                \node[vertex] at (A) (A12) {};
                \node[vertex] at (B) (B12) {};
                \node[vertex] at (C) (C12) {};
                \node[vertex] at (D) (D12) {};
                \node[vertex] at (E) (E12) {};

                \node[vertex] at ([shift={(\shiftX,0)}] A) (A22) {};
                \node[vertex] at ([shift={(\shiftX,0)}] B) (B22) {};
                \node[vertex] at ([shift={(\shiftX,0)}] C) (C22) {};
                \node[vertex] at ([shift={(\shiftX,0)}] D) (D22) {};
                \node[vertex] at ([shift={(\shiftX,0)}] E) (E22) {};
                \draw[draw=edge-color,line width=\lineWidth mm] (C22) -- (D22);
                \draw[draw=edge-color,line width=\lineWidth mm] (A22) -- (D22);
                \draw[draw=edge-color,line width=\lineWidth mm] (B22) -- (D22);

                \fill[fill=area-color] ([shift={(2*\shiftX,0)}] A) -- ([shift={(2*\shiftX,0)}] B) -- ([shift={(2*\shiftX,0)}] C) -- cycle;
                \node[vertex] at ([shift={(2*\shiftX,0)}] A) (A32) {};
                \node[vertex] at ([shift={(2*\shiftX,0)}] B) (B32) {};
                \node[vertex] at ([shift={(2*\shiftX,0)}] C) (C32) {};
                \node[vertex] at ([shift={(2*\shiftX,0)}] D) (D32) {};
                \node[vertex] at ([shift={(2*\shiftX,0)}] E) (E32) {};
                \draw[draw=edge-color,line width=\lineWidth mm] (D32) -- (C32) -- (B32) -- (D32) -- (A32) -- (C32) -- (E32) -- (A32) --(B32);

                \fill[fill=area-color] ([shift={(3*\shiftX,0)}] E) -- ([shift={(3*\shiftX,0)}] C) -- ([shift={(3*\shiftX,0)}] A) -- cycle;
                \fill[fill=area-color] ([shift={(3*\shiftX,0)}] A) -- ([shift={(3*\shiftX,0)}] B) -- ([shift={(3*\shiftX,0)}] C) -- cycle;
                \node[vertex] at ([shift={(3*\shiftX,0)}] A) (A42) {};
                \node[vertex] at ([shift={(3*\shiftX,0)}] B) (B42) {};
                \node[vertex] at ([shift={(3*\shiftX,0)}] C) (C42) {};
                \node[vertex] at ([shift={(3*\shiftX,0)}] D) (D42) {};
                \node[vertex] at ([shift={(3*\shiftX,0)}] E) (E42) {};
                \draw[draw=edge-color,line width=\lineWidth mm] (D42) -- (C42) -- (B42) -- (D42) -- (A42) -- (C42) -- (E42) -- (A42) --(B42);

                \node[vertex] at ([shift={(0,\shiftY)}] A) (A11) {};
                \node[vertex] at ([shift={(0,\shiftY)}] B) (B11) {};
                \node[vertex] at ([shift={(0,\shiftY)}] C) (C11) {};

                \node[vertex] at ([shift={(\shiftX,\shiftY)}] A) (A21) {};
                \node[vertex] at ([shift={(\shiftX,\shiftY)}] B) (B21) {};
                \node[vertex] at ([shift={(\shiftX,\shiftY)}] C) (C21) {};

                \node[vertex] at ([shift={(2*\shiftX,\shiftY)}] A) (A31) {};
                \node[vertex] at ([shift={(2*\shiftX,\shiftY)}] B) (B31) {};
                \node[vertex] at ([shift={(2*\shiftX,\shiftY)}] C) (C31) {};
                \draw[draw=edge-color,line width=\lineWidth mm] (A31) -- (B31) -- (C31) -- (A31);

                \fill[fill=area-color] ([shift={(3*\shiftX,\shiftY)}] A) -- ([shift={(3*\shiftX,\shiftY)}] B) -- ([shift={(3*\shiftX,\shiftY)}] C) -- cycle;
                \node[vertex] at ([shift={(3*\shiftX,\shiftY)}] A) (A41) {};
                \node[vertex] at ([shift={(3*\shiftX,\shiftY)}] B) (B41) {};
                \node[vertex] at ([shift={(3*\shiftX,\shiftY)}] C) (C41) {};
                \draw[draw=edge-color,line width=\lineWidth mm] (A41) -- (B41) -- (C41) -- (A41);
                
            \end{tikzpicture}
            \vspace*{1pt}
		\caption*{(b-1) Filtration $X_b$. }
	\end{subfigure}
	\hspace*{16pt}
        \hfill
	\begin{subfigure}[c]{0.35\linewidth}
 \begin{tikzpicture}
	\draw[->] (-0.3,0) -- (2.5,0) node[right] {$r$};
	\draw[->] (0,-0.3) -- (0,2.5) node[above] {$r$};
	
	\foreach \x [count=\i] in {1.5,1.732,2}{
            \pgfmathsetmacro{\coord}{4*(\x-1.5)}
		\draw (\coord,-0.2) -- (\coord,0.2);
            \ifnum\i=2
                  \node[below] at (\coord, -0.2) {$\sqrt{3}$};
                \else
                  \node[below] at (\coord, -0.2) {\x};
                \fi
	}
	
	\foreach \y [count=\i] in {1.5,1.732,2}{
            \pgfmathsetmacro{\coord}{4*(\y-1.5)}
		\draw (-0.2,\coord) -- (0.2,\coord);
                \ifnum\i=2
                  \node[left] at (-0.2,\coord) {$\sqrt{3}$};
                \else
                  \node[left] at (-0.2,\coord) {\y};
                \fi
	}
	
	\draw[gray, dashed] (0,0) -- (2,2);

        \pgfmathsetmacro{\sqrtthree}{4*(1.732-1.5)}
	
	\filldraw (\sqrtthree,2) circle (2pt);
	\node[above right] at (\sqrtthree,2) {$(\text{birth}=\sqrt{3},\text{death}=2)$};
	
	\filldraw (2,\sqrtthree) circle (2pt);
	\node[below right] at (2,\sqrtthree) {$(\text{death}=2,\text{birth}=\sqrt{3})$};
\end{tikzpicture}
            \vspace*{0.5pt}
		\caption*{(b-2) Connected persistence diagram $D_b$.}
	\end{subfigure}
	\caption{
            \mbox{(a-1)} 
            A $\cln$-filtration at critical radii, where all arrows indicate inclusion maps.
            \mbox{(a-2)} 
            The upper and lower triangular regions represent persistence diagrams at dimension one for respective rows. 
            The connecting line is 
            consistent with the observation that 
            the two generators share a homologous cycle.
            \mbox{(b-1)} In this filtration,
            the generator in the lower row becomes a boundary when mapped to the upper row.
            \mbox{(b-2)} 
            The two generators corresponding to the cycles in the upper row and lower row are disconnected, 
            indicating that they are irrelevant cycles.
        }
	\label{fig:cPD-example-ab}
\end{figure}

%
\end{exa}

Another visualization method to consider is the \emph{layered presentation} of the connected persistence diagram. 
This approach overlays the two standard persistence diagrams, and line segments are also drawn in the upper triangular region.
Although superimposed generators become less discernible, 
this method offers a clearer insight into 
how specific generators persist in the vertical direction in certain contexts.
We further illustrate this method in Section \ref{sec:silica}.

\begin{rmk}
The \emph{extended diagram} presented in \cite{edelsbrunner2022computational} also utilizes two triangular regions. Although visually similar, this structure differs from the connected persistence diagram, and they are not directly related.
\end{rmk}

\section{Models for Commutative Ladder Filtrations of Simplicial Complexes}

In this part, we introduce models to construct commutative ladder filtrations
at the simplicial complex level, enhancing the applicability of the topological invariants discussed in the prior section.


\subsection{A General Model}\label{subsec:general-model}


A \emph{filtered simplicial complex} is a pair $(K,f)$,
where $K$ is a simplicial complex and 
${f\colon K\to \extR\coloneqq\Set{\mathbb{R},\pm\infty}}$ 
is a \emph{filter} satisfying:
\[
\mbox{$f(\tau)\leq f(\sigma)$ if $\tau$ is a face of $\sigma$}.
\]
Every number $r_i\in \extR$ can be associated with a simplicial complex through the preimage
$f^{-1}\big((-\infty,r_i]\big)$.
For a strictly increasing sequence $r_1<r_2<\cdots <r_n$ in 
$\extR$, this yields a sequence of simplicial complexes:
\[
f^{-1}\big((-\infty,r_1]\big)\subseteq f^{-1}\big((-\infty,r_2]\big)\subseteq\cdots\subseteq f^{-1}\big((-\infty,r_n]\big),
\]
which forms a filtration of sublevel sets. 

We now consider a simplicial complex $K$ equipped with two filters $f_1$ and $f_2$, 
subject to the condition $f_1(\sigma)\geq f_2(\sigma)$ for any $\sigma\in K$.
Given a sequence $r_1<r_2<\cdots<r_n$, 
the pair $(K,f_y)$ with $y\in \Set{1,2}$ yields a sequence as above by defining
$X_{(i,y)}\coloneqq f_y^{-1}\big((-\infty,r_i]\big)$.
They play the role of the lower and upper rows in a $\cln$-filtration, respectively.
The condition $f_1\geq f_2$ ensures the vertical inclusions and the commutativity can also be verified.
Finally, the $\cln$-filtration specified by the triplet $(K,f_1,f_2)$ can be formulated as
\begin{equation}\label{eq:cln-filt-general-model}
    \begin{aligned}
    X\coloneqq
    \begin{tikzcd}[row sep=\tikzcdRowSepLarger,column sep=\tikzcdColumnSepLarger]
    X_{(1,2)}\arrow[r]&X_{(2,2)}\arrow[r]&\cdots\arrow[r]&X_{(n,2)}\\
    X_{(1,1)}\arrow[r]\arrow[u]&X_{(2,1)}\arrow[r]\arrow[u]&\cdots\arrow[r]&X_{(n,1)}\arrow[u]
    \end{tikzcd}
    \end{aligned}.
\end{equation}



\subsection{Thinning Models for Point Cloud Data}\label{subsec:model-point-cloud}


%

Let $P$ be a point cloud and $P'$ be a subset of $P$. 
Let $\check{C}(P,r)$ denote the Čech complex constructed on the point cloud $P$ with ball radius $r$.
For a sequence $r_1<r_2<\cdots<r_n$ in $\extR$, we have the following filtration:
%
\begin{equation}\label{eq:cln-filt-point-cloud-model}
    \begin{tikzcd}[row sep=\tikzcdRowSepLarger,column sep=\tikzcdColumnSepLarger]
    \check{C}(P,r_1)\arrow[r]
    &\check{C}(P,r_2)\arrow[r]
    &\cdots\arrow[r]
    &\check{C}(P,r_n)\\
    \check{C}(P',r_1)\arrow[r]\arrow[u]
    &\check{C}(P',r_2)\arrow[r]\arrow[u]
    &\cdots\arrow[r]
    &\check{C}(P',r_n)\arrow[u]
    \end{tikzcd},
\end{equation}
where the definition of the Čech complex guarantees the commutativity.
To incorporate this construction into the general model introduced above, 
we set the simplicial complex $K$ to be $\check{C}(P,r_n)$,
and define the induced filter $f_2$ for $\sigma\in K$ as:
\[
f_2(\sigma)=\min \Set{r_i\in\extR|\sigma \in \check{C}(P,r_i)}.
\]
Given that each simplex $\sigma$ can be expressed in the form $\Set{v_0,v_1,\ldots,v_d}$ with each $v_i\in P$, the filter $f_1$ on $K$ can be defined as:
\[
f_1(\sigma)
=
\begin{cases}
    f_2(\sigma)& \sigma\subseteq P',\\
    +\infty & \sigma \centernot\subseteq P'.
\end{cases}
\]
The filtration (\ref{eq:cln-filt-point-cloud-model}) can then be retrieved 
from the triplet $(K,f_1,f_2)$ following the steps in Section \ref{subsec:general-model}.

\begin{rmk}
Thinning in topological data analysis, introduced in \cite{yasuaki_ippei_kazuto_2019}, 
can incorporate a variety of patterns,
such as removing points forming certain shapes or employing a density function.
An example can be found in Section \ref{subsec:FCC-HCP}.
Besides what is described here, there exist other methods for creating a filtration from a point cloud $P$, 
including non-constant point removal or removing higher-dimensional simplices. 
The multi-cover method described in \cite{osang2021topological} can also generate a $\cln$-filtration consistent with our general model.
\end{rmk}


\subsection{Two Models from Random Simplicial Complexes}\label{subsec:model-abs-complex}


\subsubsection{Clique Complex Model}\label{subsubsec:CC-model}


Our first random model employs the Erdős-Rényi random graph process
for generating commutative ladders.
Consider the complete graph with $m$ vertices,
denoted by $K_m=(V_m, E_m)$,
where $V_m$ is the set of vertices and 
$E_m$ is the set of edges.
We associate each edge $e\in E_m$ with 
two independent random variables $T_e$ and $\widetilde{T}_e$, 
both follow the standard uniform distribution $U(0,1)$. 
Given $t\in[0,1]$, we define the following two increasing stochastic processes of subgraphs of $K_m$:
\[
    K_m^1(t)\coloneqq\big(V_m,\set{e\in E_m|T_e\leq t}\big),\qquad
    K_m^2(t)\coloneqq\big(V_m,\set{e\in E_m|T_e\cdot\widetilde{T}_e\leq t}\big).
\]
The stochastic process $K_m^1$ is called the Erdős-Rényi random graph process \cite{erdHos1960evolution}.
Recall that the clique complex $\Delta(G)$ 
of a graph $G$
is the abstract simplicial complex that includes all cliques (\textit{i.e.} complete subgraphs) of $G$. 
The two stochastic processes defined above induce two increasing stochastic processes of simplicial complexes
$\Delta(K_m^{1})$ and $\Delta(K_m^{2})$.
Notice that the condition $\Delta(K_m^{1})(t)\subseteq \Delta(K_m^{2})(t)$ holds for any $t\in[0,1]$, as ensured by the multiplication in $K^2_m$'s definition.
By selecting
$n$ values $t_1<t_2<\cdots<t_n\in [0,1]$ 
and setting $X_{(i,j)}\coloneqq \Delta(K_m^{j})(t_i)$,
we obtain a random $\cln$-filtration.

To fit it into (\ref{eq:cln-filt-general-model}) in the general model,
we set $K$ to be $\Delta(K_m)$.
For $j=1,2$, each filter $f_j$ on $K$ is defined as below for $\sigma=\Set{v_0,\ldots,v_d}$:
\[
f_1(\sigma)\coloneqq \max\limits_{0\leq<p<q\leq d}t_{v_pv_q}
\qquad\mbox{and}\qquad
f_2(\sigma)\coloneqq \max\limits_{0\leq<p<q\leq d}t_{v_pv_q}\cdot \widetilde{t}_{v_pv_q}.
\]
In this formulation, the subscript $v_pv_q$ represents the edge incident to vertices $v_p$ and $v_q$.
The terms $t_e$ and $\widetilde{t}_{e}$ are realizations of $T_e$ and $\widetilde{T}_e$, respectively.

\subsubsection{\texorpdfstring{$d$}{d}-Linial-Meshulam Model}\label{subsubsec:dLM-model}
We introduce another model adapted from the $d$-Linial-Meshulam process \cite{hiraoka2017minimum,linial2006homological}.
This stochastic process can be seen as a generalization of the Erdős-Rényi random graph process. 
It begins with a skeleton of a simplicial complex 
and then progressively adds simplices that are one-dimensional higher than the initial skeleton.
We outline the components required for this model as follows.
\begin{itemize}
    \item A vertex set $[m]\coloneqq\Set{1,\ldots,m}$.
    \item The largest abstract simplicial complex over $[m]$, denoted as $\Delta\coloneqq \Delta_{m-1}$, which has dimension $m-1$.
    \item A chosen dimension $d$ where $1\leq d\leq m-1$.
    \item The $(d-1)$-skeleton of $\Delta$,
    denoted as $\Delta^{(d-1)}$, 
    comprising all simplices in $\Delta$ having a dimension no greater than $d-1$.
    \item The set of all $d$-simplices in $\Delta$, represented as $\Delta_d$.
    \item Each $\sigma\in \Delta_d$ is associated with two independent random variables 
    $T_\sigma$ and $\widetilde{T}_\sigma$ that follow the standard uniform distribution $U(0,1)$.
    \item Two increasing stochastic processes of simplicial complexes:
    \[
    \mathcal{K}^{(d)}_{1}(t)\coloneqq\Delta^{(d-1)}\sqcup \Set{\sigma\in\Delta_d | T_{\sigma}\leq t},\qquad
    \mathcal{K}^{(d)}_{2}(t)\coloneqq\Delta^{(d-1)}\sqcup \Set{\sigma\in\Delta_d | T_{\sigma}\cdot \widetilde{T}_{\sigma}\leq t}.
    \]
\end{itemize}
For both $j=1,2$, 
we have $\mathcal{K}^{(d)}_{j}(0)=\Delta^{(d-1)}$ 
at the start of each filtration almost surely
and $\mathcal{K}^{(d)}_{j}(1)=\Delta^{(d)}$ at the end.
Analogous to the previous model,
selecting $n$ values $t_1<t_2<\cdots<t_n$ in $[0,1]$ 
and setting $X_{(i,j)}\coloneqq \mathcal{K}_{j}^{(d)}(t_i)$
yields a random $\cln$-filtration.
A realization of this random $\cln$-filtration 
can be established using the triplet $(K,f_1,f_2)$ as outlined in the general model similarly.




\section{Experiments and Analysis}
In this section, we discuss computational results
about the newly introduced invariants.
We begin by analyzing the occurrence of non-intervals
in the three different models outlined in the previous section,
then show the findings by applying our tools
to material structures.






\subsection{Non-intervals in Three Different Models}
This subsection explores the under-investigated realm of non-interval representations within multiparameter persistence.
We use the toolkit in Section \ref{sec:commutative-ladder} 
to analyze the models in Section \ref{subsec:model-point-cloud} and \ref{subsec:model-abs-complex}. 
We will explore persistence modules of both finite-type and infinite-type commutative ladders. 

\vspace*{2pt}
\noindent\textbf{Finite-type Commutative Ladders}

Recall from Theorem \ref{thm:commutative-ladder-rep-type} that the representation type of $\cln$ is finite when $n\leq 4$.
We will focus on $\clf$ 
because any representation of $\operatorname{CL}(1)$ or $\operatorname{CL}(2)$ is interval-decomposable,
and the non-interval representations of $\clt$ can be embedded into $\rep\bigl(\clf\bigr)$.

\begin{figure}[htbp]
    \centering
    \begin{adjustbox}{width=1.22\textwidth,center}
    \includegraphics[width=\textwidth]{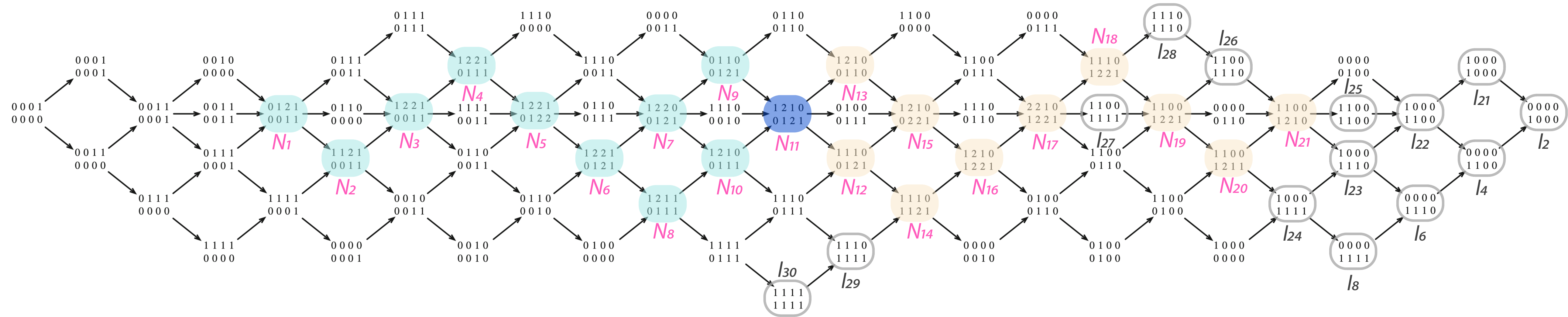}
    \end{adjustbox}
    \caption{The $76$ representatives in $\mathcal{L}$. 
    Non-intervals are numbered such that $N_{11-i}$ is isomorphic to the dual of $N_{11+i}$ for $i\in\set{0,\ldots,10}$.
    }
    \label{fig:AR-quiver-cl4}
\end{figure}

Let $\mathcal{L}$ be a complete set of representatives of the isomorphism classes of indecomposables in $\rep\bigl(\clf\bigr)$.
The Auslander-Reiten quiver
\footnote{
The Auslander-Reiten quiver of $\rep\bigl(\clf\bigr)$ consists of elements in $\mathcal{L}$ as vertices and irreducible morphisms among them as arrows. We use it to represent elements in $\mathcal{L}$ here.
See \cite[Chapter 4]{assem2006elements} for details.} 
reveals that there are 21 non-intervals and 55 intervals in $\mathcal{L}$ \cite[Fig. 17]{escolar2016persistence}. 
Non-intervals and some intervals are indexed as shown in Figure \ref{fig:AR-quiver-cl4} for easier reference in subsequent discussions, where each class is represented by its dimension vector.

\begin{rmk}
There are two non-interval representatives of isomorphism classes in $\clt$, represented by dimension vectors
$\begin{pmatrix}
1&2&1\\
0&1&1
\end{pmatrix}$
 and 
$\begin{pmatrix}
1&1&0\\
1&2&1
\end{pmatrix}$.
They can be embedded into either $N1$ and $N9$ or $N13$ and $N21$ in $\rep\bigl(\clf\bigr)$, as shown in Figure \ref{fig:AR-quiver-cl4}.
\end{rmk}




\noindent\textbf{Infinite-type Commutative Ladders} 

Unlike the scenario in the finite-type cases, 
our understanding of non-intervals in infinite-type cases is limited.
Despite this,
the connected persistence diagrams still provide some clues about them.
The property 
below establishes 
a negative multiplicity 
as an indicator of the presence 
of non-interval summands.

\begin{prop}\label{prop:negative-not-interval-decomposable}
Let $M$ be a representation of $\cln$.
    If there exists an essential assignment $\xi$ and an interval $I$ of $\cln$ such that $\delta_M^{\xi}(I)<0$,
    then $M$ has an indecomposable non-interval direct summand.
\begin{proof}
    The negativity of $\delta_M^{\xi}(I)$ 
    implies $d_M(V_I)\neq \souwa{S\subseteq \Cov I}(-1)^{\# S}c_M^{\xi}({\bigvee S})$ by Remark \ref{rmk:interval-approximation-gpq}.
    Therefore, $M$ is not an interval-decomposable representation by Theorem \ref{thm:multiplicity-from-c}.
\end{proof}
\end{prop}

\subsubsection{On Point Cloud Model and Clique Complex Model}
We demonstrate a comparative analysis of the results from the Point Cloud Model in Section \ref{subsec:model-point-cloud}
and the Clique Complex Model in Section \ref{subsubsec:CC-model}.
Throughout this discussion, 
we maintain the notations established in the individual model's specification.
Implementation details for building filtrations are listed below.
\renewcommand{\arraystretch}{1.5}
\noindent
\begin{tabularx}{\textwidth}{
    @{\hspace{1.5em}}
    >{\leavevmode\llap{\textbullet~}\raggedright}
    X
    @{\quad\hspace{1.5em}}
    >{\leavevmode\llap{\textbullet~}\raggedright\arraybackslash}
    X
    @{}
  }
  \multicolumn{1}{X}{\centering\bfseries Point Cloud Model} &
    \multicolumn{1}{X}{\centering\bfseries Clique Complex Model} \\
  \midrule
  Construct a set $P$ comprising uniformly randomly distributed points confined 
  in the unit cube of $\mathbb{R}^2$ or $\mathbb{R}^3$,
  where the number of points ranges from $5$ to $30$.
  Subsequently, create a point cloud $P'$ by randomly choosing a non-empty proper subset of $P$.
  & Construct a complete graph with $m$ vertices where $m\in[4,30]$. 
  Subsequently, generate samples of the random variables $T_e$ and $\widetilde{T}_e$ for $e\in E_m$.
  \\
  Set the homology dimension $d$ as $1$ if $P\subseteq \mathbb{R}^2$.  
  Otherwise, assign $d$ as either $1$ or $2$.
  Apply one-parameter persistent homology $H_d$ functor on $\check{C}(P)$ and $\check{C}(P')$ to get their critical filtration values.
  Collect these values and 
  denote their disjoint union with duplicates removed as $\mathscr{V}\coloneqq \Set{t_1,t_2,\ldots,t_{\abs{\mathscr{V}}}}$.
  Discard this configuration if $\abs{\mathscr{V}}<4$.
  & 
  Assign the homology dimension $d$ to be $1$ or $2$.
  Apply one-parameter persistent homology $H_d$ functor on $\Delta(K^1_m)$ and $\Delta(K^2_m)$ to get their critical filtration values.
  Collect these values and 
  denote their disjoint union with duplicates removed as $\mathscr{V}\coloneqq \Set{t_1,t_2,\ldots,t_{\abs{\mathscr{V}}}}$.
  Discard this configuration if $\abs{\mathscr{V}}<4$.
   \\
   For the finite-type cases, choose four radii randomly from $\mathscr{V}$.
   For the infinite-type cases, select an integer $n$ such that $4\leq n\leq \min(\abs{\mathscr{V}},50)$,
   and build a $\cln$-filtration by choosing $n$ radii randomly.
   &
   As detailed in the corresponding entry in the left column.
\end{tabularx}
\vspace*{\baselineskip}

\vspace*{\baselineskip}
We begin by examining the interval-decomposability of non-trivial representations.
For a representation of $\clf$, we compute its indecomposable decomposition 
to identify the presence of any non-interval components.
This allows us to determine the proportion of the number of representations that are not interval-decomposable compared to the total number of non-trivial representations.
For a general representation of $\cln$,
we turn to its connected persistence diagram and look for the presence of negative multiplicities, providing a lower-bound estimate of the aforementioned proportion.
Table \ref{table:point-cloud-non-interval-freq-cl4} presents the proportions obtained from various settings,
using either $\clf$ (for indecomposable decomposition) or cPD (for connected persistence diagram) as specified above.

\setlength{\heavyrulewidth}{1.5pt}
\begin{table}[htbp]
\centering
\begin{tabular}{ccccccc}
\toprule
 & \multicolumn{2}{c}{Point Cloud Model: $\mathbb{R}^2$} & \multicolumn{2}{c}{Point Cloud Model: $\mathbb{R}^3$} & \multicolumn{2}{c}{Clique Complex Model}\\
\midrule
 & $\clf$ & cPD & $\clf$ & cPD & $\clf$ & cPD\\
\cmidrule(r){2-3} \cmidrule(lr){4-5} \cmidrule(l){6-7}
$H_1$ & 0.034\% & 3.38\% & 0.097\% & 8.83\% & 1.727\% & 34.74\% \\
$H_2$ & N/A & N/A & 0.042\% & 1.25\% & 0.001\% & 0.037\% \\
\bottomrule
\end{tabular}
\vspace*{\baselineskip}
\caption{The proportion or its lower-bound of non-trivial representations that are not interval-decomposable for the Point Cloud Model embedded in $\mathbb{R}^2$, $\mathbb{R}^3$, and the Clique Complex Model. 
Tools used are indicated in the second row. 
\protect\footnotemark
}
\label{table:point-cloud-non-interval-freq-cl4}
\end{table}
\footnotetext{The number of total non-trivial representations in the row with $H_1$ are 132k, 132k, 66k, 66k, 211k, 221k respectively, and 49k, 50k, 141k, 142k for the row with $H_2$.}

\begin{figure}[htbp]
    \centering
    \includegraphics[width=0.78\textwidth]{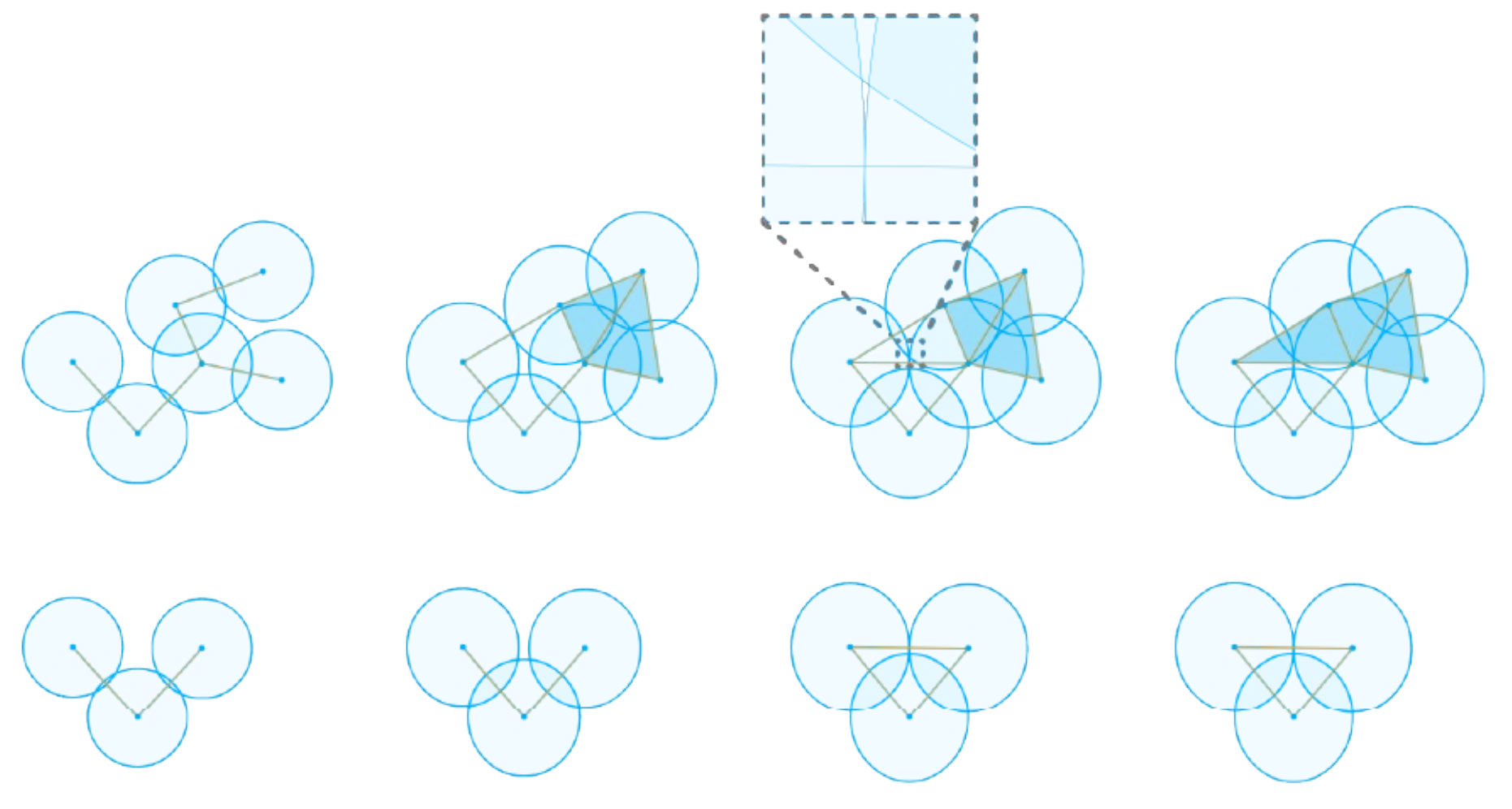}
    \caption{An example of a $\clf$-filtration.
    Its one-dimensional homology module corresponds to $N1$ in Figure \ref{fig:AR-quiver-cl4}.
    The new hole, encased in the square, is crucial for this non-interval 
    despite its extremely short lifetime.
    While the generator in the lower row also has a very short lifetime, 
    it slightly outlives the prior generator, leading to a non-interval component.
    This fragile structure is prone to perturbations, explaining the infrequent occurrence of 
    non-intervals in the Point Cloud Model.}
    \label{fig:non-intervals-cl4-configuration-example}
\end{figure}

Upon comparing the data from both models in homology dimension one using $\clf$,
it is evident that the Point Cloud Model exhibits a significantly lower proportion of
non-interval-decomposable representations compared to the Clique Complex Model.
This relatively infrequent occurrence of non-intervals in the Point Cloud Model
can be primarily attributed to the geometry of Čech complexes, as illustrated in Figure \ref{fig:non-intervals-cl4-configuration-example}.
Additionally, it is worth noting that the proportion of non-interval-decomposable representations in $\mathbb{R}^3$ is higher than that in $\mathbb{R}^2$, which is likely due to increased tolerance for perturbations in $\mathbb{R}^3$ owing to its additional degree of freedom.

\begin{figure}[htbp]
\begin{subfigure}[t]{0.49\linewidth}
\centering
\includegraphics[width=\textwidth]{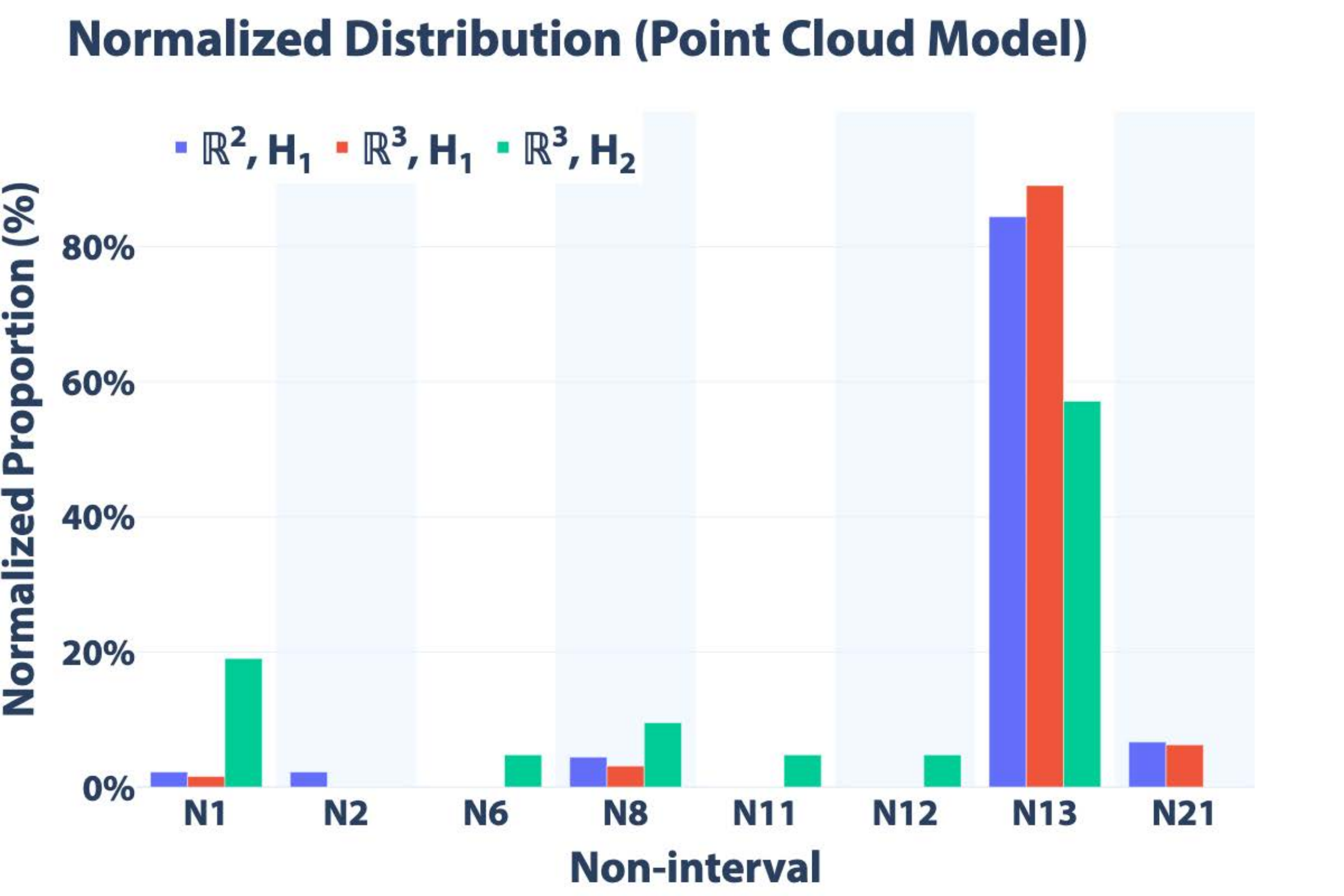}
\caption{} \label{fig:point-cloud-model-cl4-dist}
\end{subfigure}
\hspace*{\fill}
\begin{subfigure}[t]{0.49\linewidth}
\centering
\includegraphics[width=\textwidth]{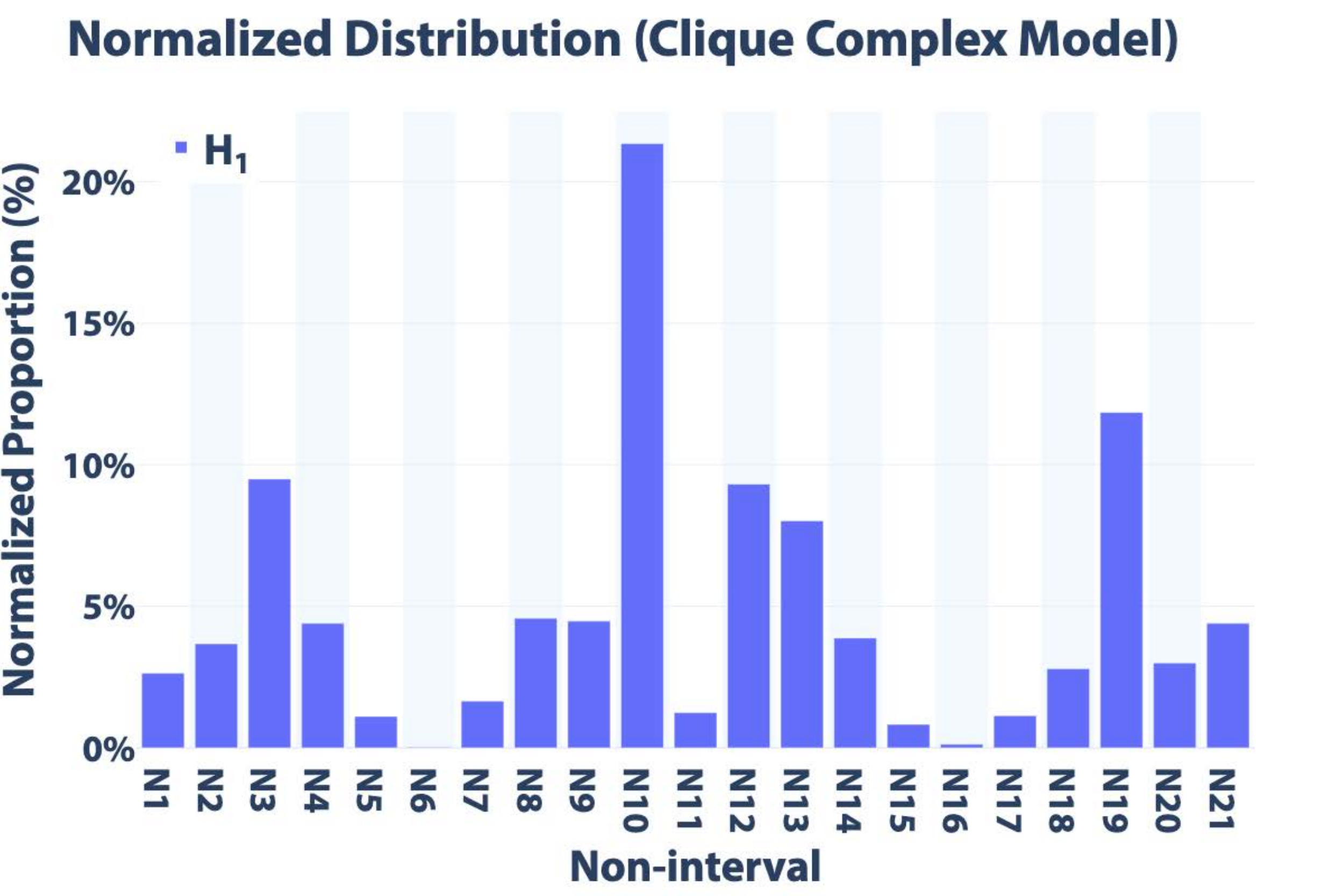}
\caption{} \label{fig:erdos-renyi-model-cl4-dist}
\end{subfigure}
\caption{
Normalized proportions are presented across various settings, 
only showing non-interval representatives observed in the computational results.
(a) The distribution of non-intervals in the Point Cloud Model.
(b)
The distribution of non-intervals in the Clique Complex Model with homology dimension one.
\protect\footnotemark
}
\label{fig:non-intervals-dist-pc-er}
\end{figure}

Besides the differences in interval-decomposability,
non-intervals also exhibit an uneven distribution across the set of all 21 representatives in different models.
Figure \ref{fig:non-intervals-dist-pc-er} illustrates the proportion of each representative's sum of multiplicities, normalized against the total multiplicity of all the non-intervals.

In the Point Cloud Model, 
only eight distinct representatives appear among all representatives.
Non-intervals predominantly concentrate on $N13$, which has one of the lowest 1-norms of the dimension vector among all representatives.
In contrast, in the Clique Complex Model, 
the diversity of non-intervals increases, and they are also less concentrated.
All 21 representatives are observed, 
with no apparent preference towards representatives possessing the lowest 1-norms of dimension vectors.

\footnotetext{
Among all the 141k non-trivial representations obtained using $H_2$, only one representation is not interval-decomposable, consisting of an $N14$ as its summand.}

\begin{figure}[htbp]
\begin{subfigure}[t]{0.49\linewidth}
\centering
\includegraphics[width=\textwidth]{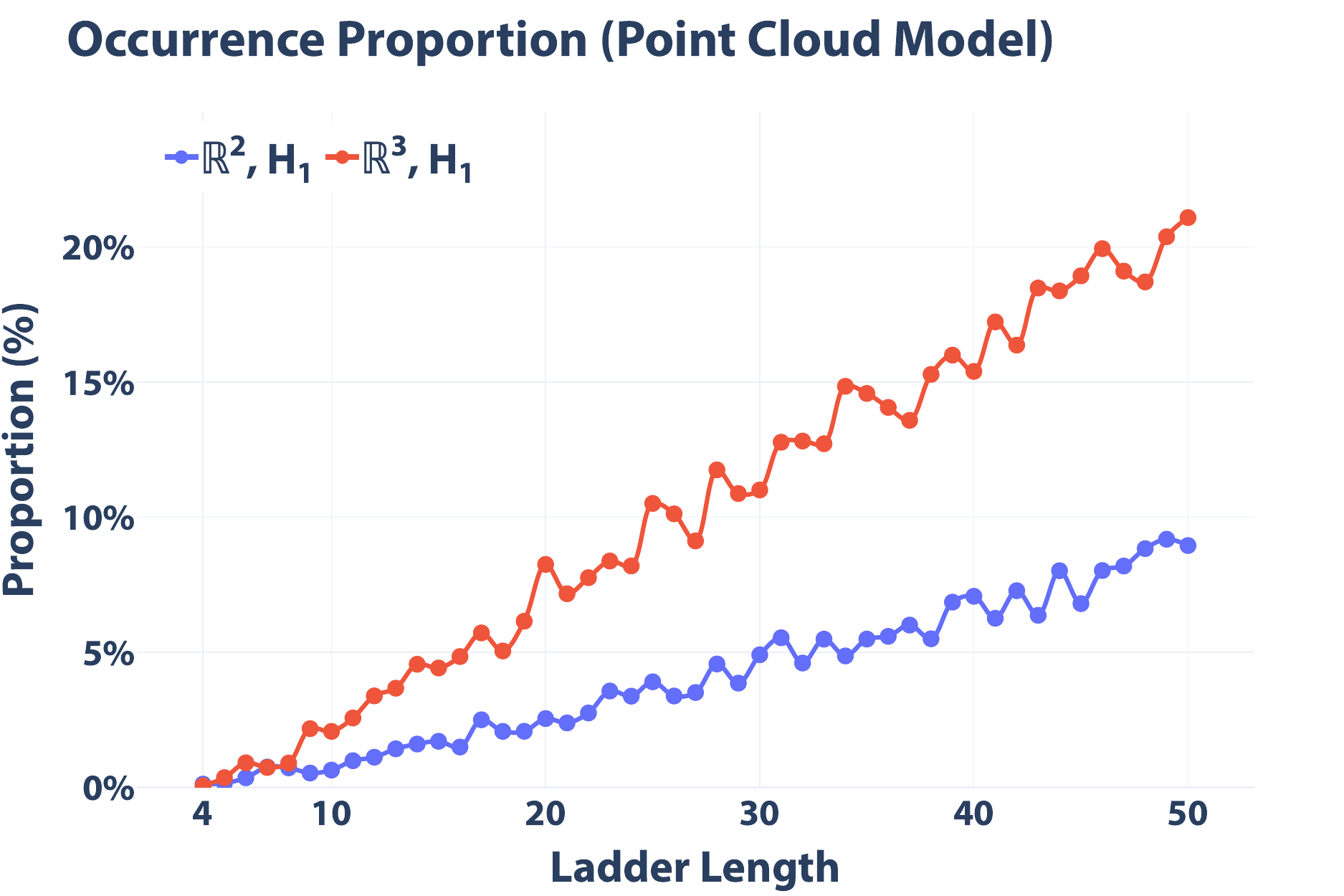}
    \caption{}
 \label{fig:point-cloud-model-cln-detection-rate}
\end{subfigure}
\hspace*{\fill}
\begin{subfigure}[t]{0.49\linewidth}
\centering
\includegraphics[width=\textwidth]{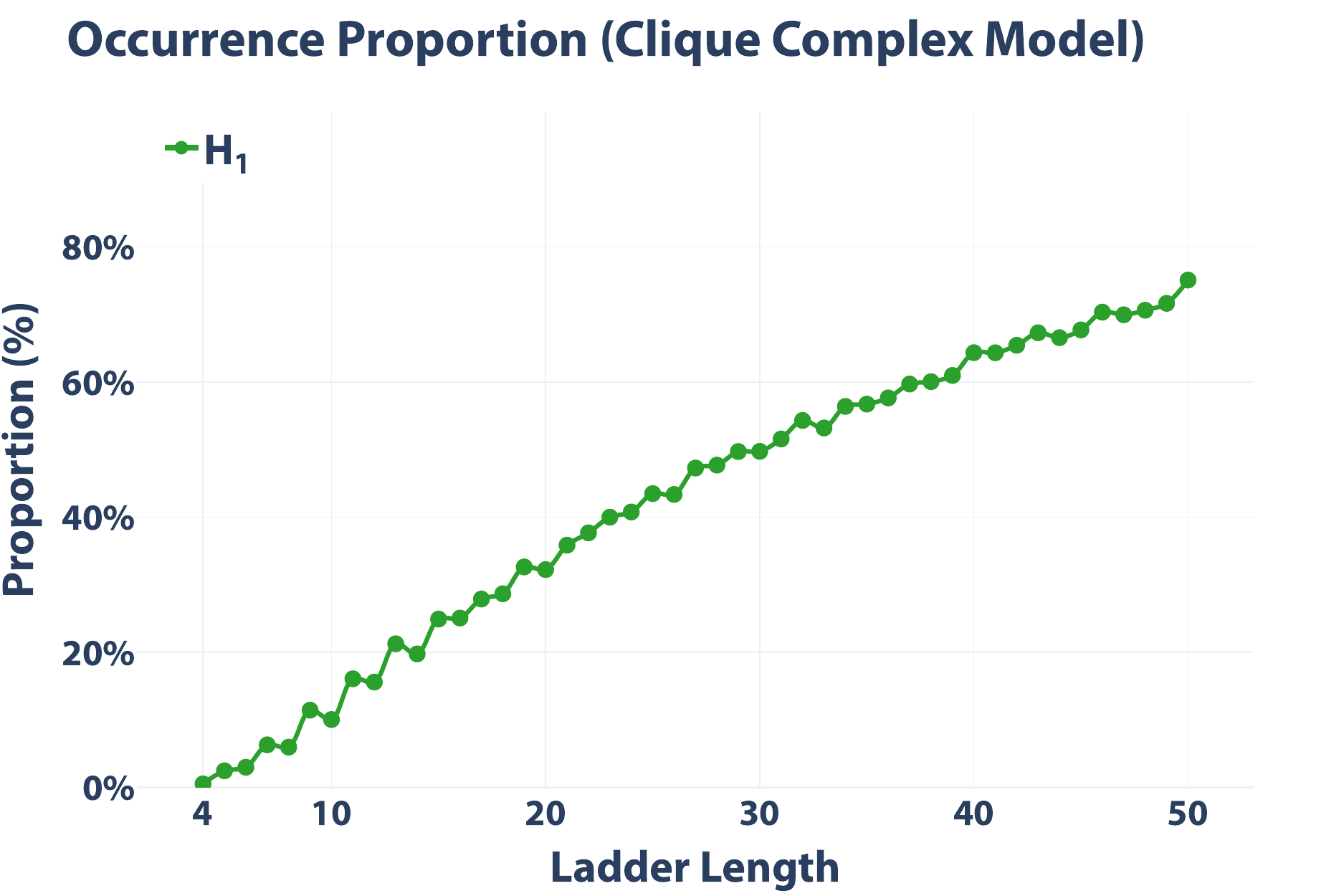}
\caption{}
\label{fig:clique-complex-model-cln-detection-rate}
\end{subfigure}
\caption{
Occurrence proportions of negative multiplicities in connected persistence diagrams of
homology dimension one.
(a) For the Point Cloud Model in $\mathbb{R}^2$ and $\mathbb{R}^3$.
(b) For the Clique Complex Model.
\protect\footnotemark
}
\label{fig:ladder-lenght-occurrence-proportion}
\end{figure}

Next, we look at the effect of ladder length $n$ on the detectability of a non-interval component.
Figure \ref{fig:ladder-lenght-occurrence-proportion} 
shows the ratios of the connected persistence diagrams with a negative multiplicity. 

\footnotetext{Results for homology dimension two in both the Point Cloud Model and the Clique Complex Model are omitted due to insufficient data points after grouping by ladder length.}

With shorter ladder lengths, many details in the original representation are lost, 
leading to a lower proportion of the occurrence of negative multiplicities. 
On the other hand, as we increase the length, more topological information is being scrutinized,
increasing the likelihood of finding a negative multiplicity.
We also see that the Clique Complex Model demonstrates a much higher proportion than the Point Cloud Model, aligning with the observation obtained using indecomposable decomposition.

\subsubsection{On \texorpdfstring{$d$}{d}-Linial-Meshulam Model}

\begin{figure}[htbp]
    \centering
    \includegraphics[width=0.78\textwidth]{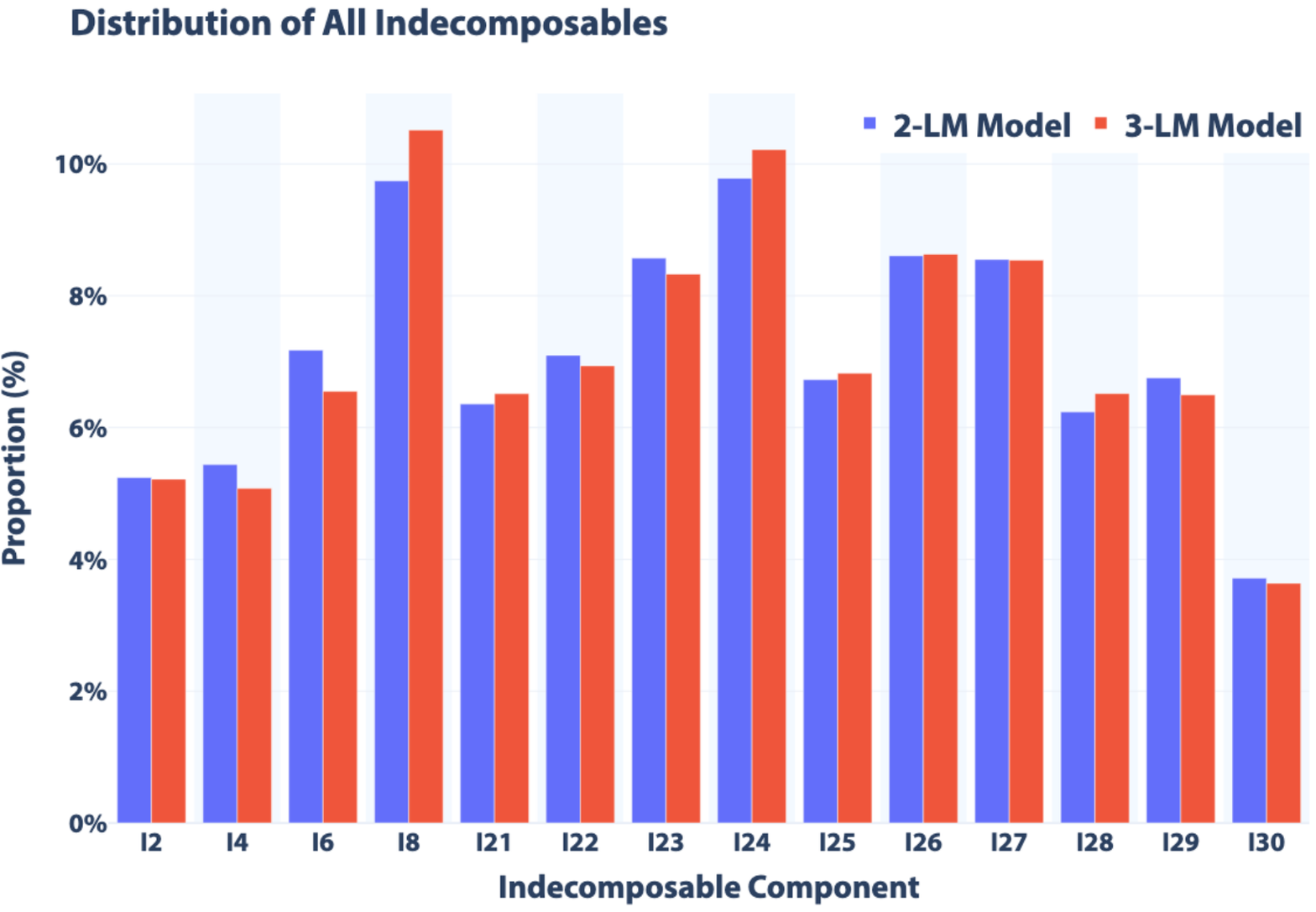}
    \caption{The distribution of non-vanishing indecomposable components in the computational results.\protect\footnotemark }
    \label{fig:dLM-cl4-distribution}
\end{figure}

This part analyzes the $d$-Linial-Meshulam Model outlined in Section \ref{subsubsec:dLM-model},
with the homology dimension fixed at $d-1$.
Our numerical computations focus on cases where $d=2$ and $d=3$, 
with the number of vertices $m$ varying between $4$ and $20$.
The methodology for creating a commutative ladder filtration 
is similar to the approaches used in the prior two models.
Figure \ref{fig:dLM-cl4-distribution} summarizes 
the proportions of each representative's sum of multiplicities within the indecomposable decomposition 
obtained using $\clf$.
\footnotetext{Computed on 15k 2-LM Model instances and 14k 3-LM Model instances}

Two immediate observations from the computational outcomes are:
\begin{enumerate}
    \item There exist only interval indecomposable components.
    \item All present intervals are anchored at the bottom-left vertex of $\clf$ (refer to Figure \ref{fig:AR-quiver-cl4} for details).
\end{enumerate}
This is not merely a coincidence but an inherent characteristic of the $d$-Linial-Meshulam Model.
In the construction process, 
we always start with a $(d-1)$ skeleton $\Delta^{(d-1)}$,
this implies that all the $(d-1)$-homologous cycles of $\Delta$
are present at critical value $t=0$.
Throughout the stochastic process, $d$-simplices from $\Delta_d$ are being incorporated,
which only serves to fill $(d-1)$-cycles without introducing new ones.
It is easy to see that homologous cycles neither merge nor split during this process.
Consequently, all representations created from this process are interval-decomposable and, in particular, pivoted at the bottom left vertex.
Therefore, the inherent structure of the $d$-Linial-Meshulam Model 
makes it a systematic approach to generate 
interval-decomposable representations pivoted at the left bottom vertex.



\subsection{On Exploring Material Structures}
This subsection uses connected persistence diagrams to uncover the topological properties inherent in amorphous and crystalline structures. 
Our analysis offers a unique perspective on the topological characteristics of material structures, 
highlighting features invisible to one-parameter persistent homology. 



\subsubsection{Silica Thinning Analysis}\label{sec:silica}
\begin{figure}[htbp]
\begin{subfigure}{0.32\textwidth}
\centering
\includegraphics[width=\linewidth]{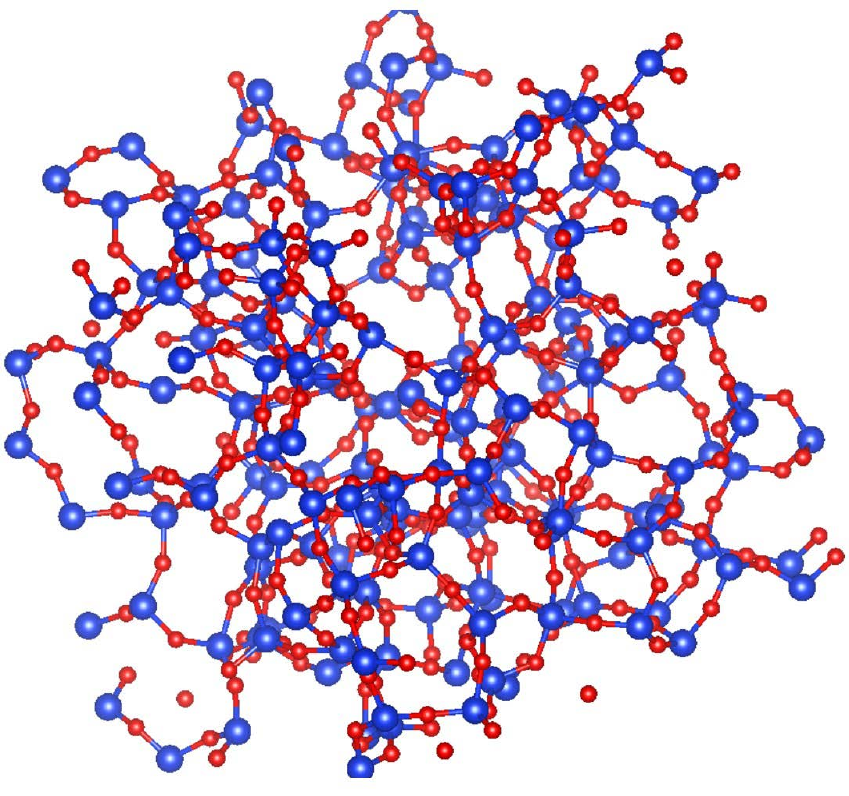}
\caption{} \label{fig:silica_input}
\end{subfigure}
\hspace*{\fill}
\begin{subfigure}{0.32\linewidth}
\centering
\includegraphics[width=\linewidth]{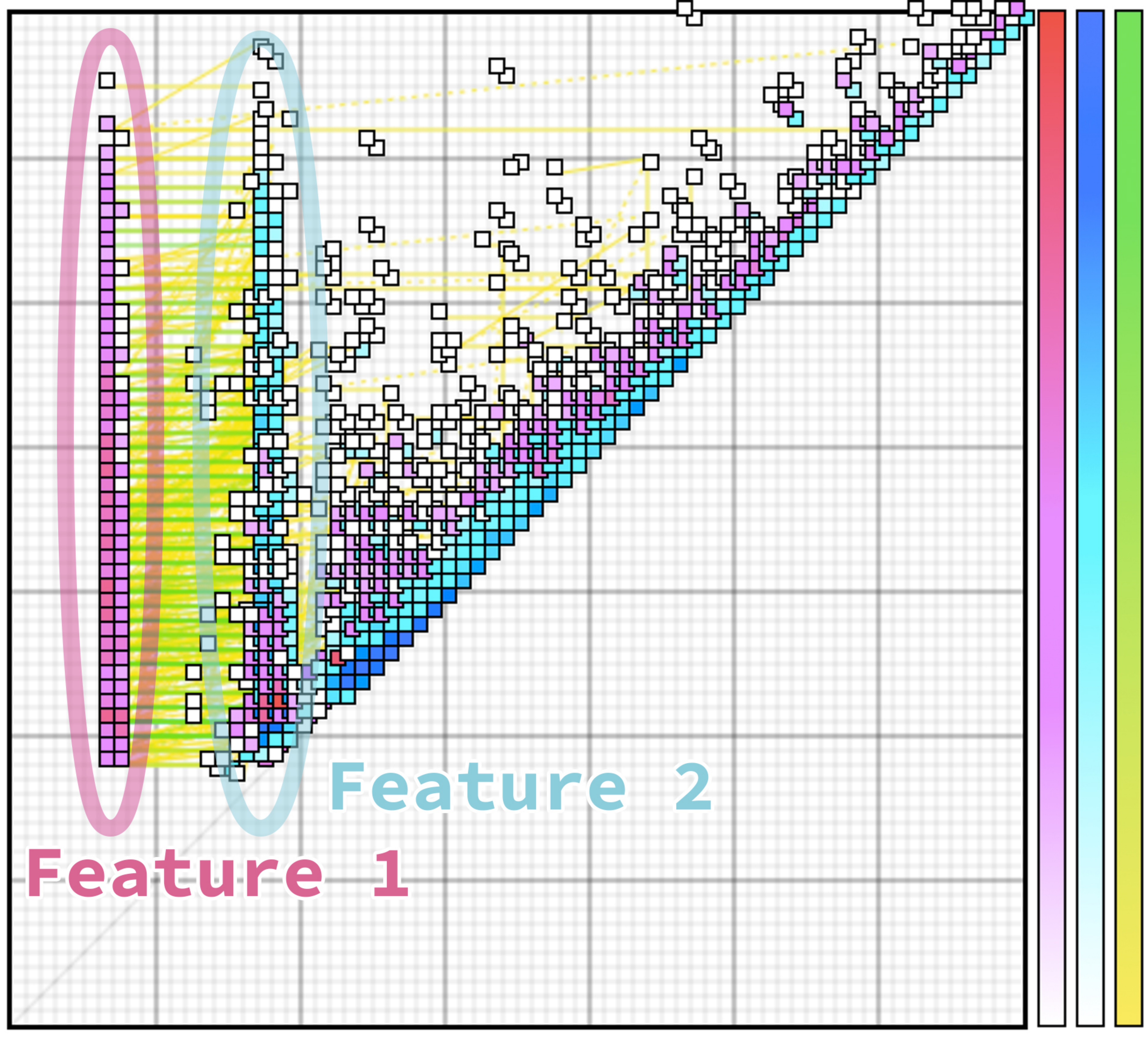}
\caption{} \label{fig:example-silica-case1}
\end{subfigure}
\hspace*{\fill}
\begin{subfigure}{0.32\linewidth}
\centering
\includegraphics[width=\linewidth]{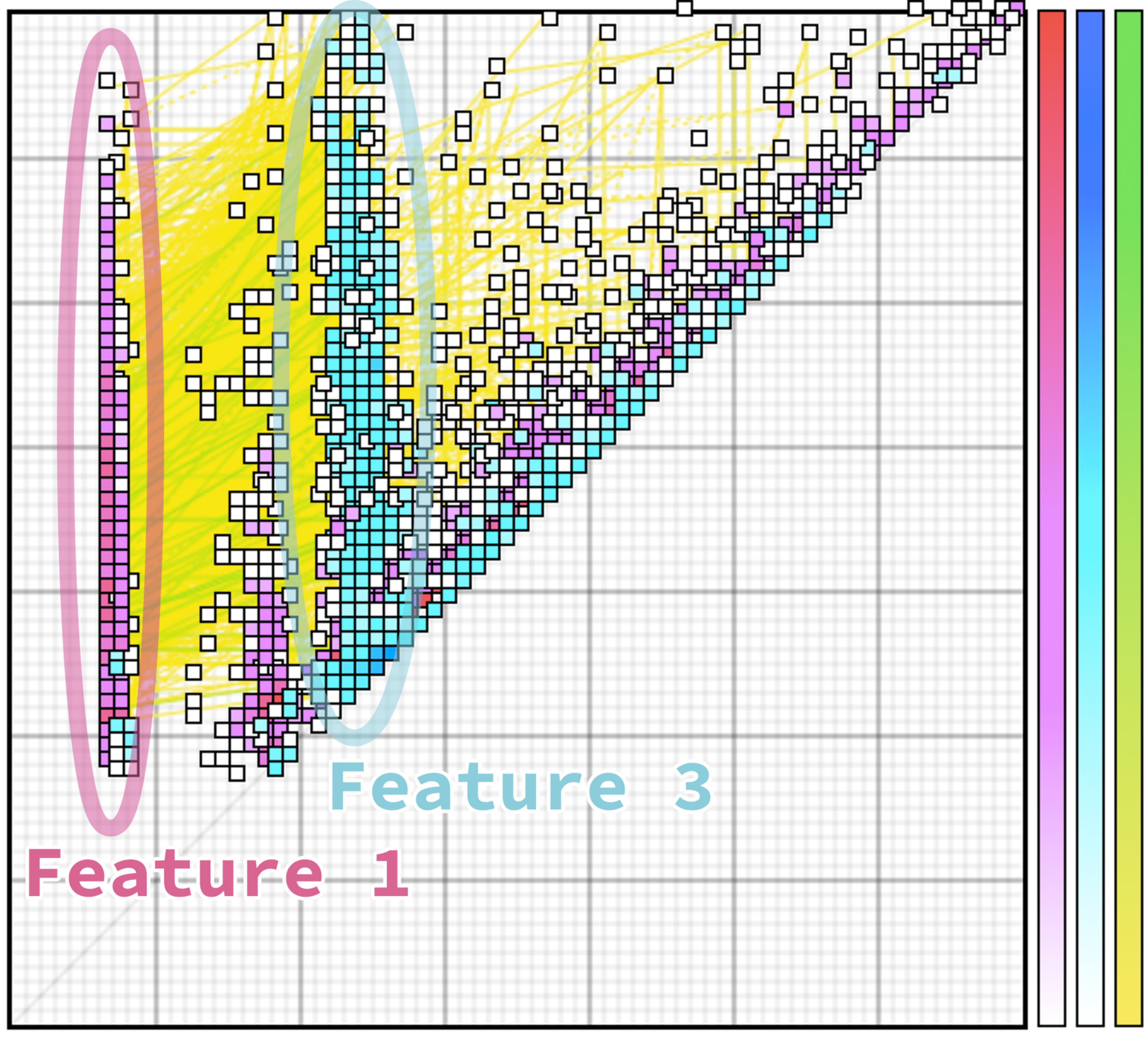}
\caption{} \label{fig:example-silica-case2}
\end{subfigure}
\caption{
(a) Atomic representation of silica with chemical bonds illustrated.
Blue balls represent silicon atoms, and red balls represent oxygen atoms.
(b) The connected persistence diagram of the configuration with all silicon atoms removed.
(c) The connected persistence diagram of the configuration with half of the oxygen atoms removed.
}
\label{fig:silica}
\end{figure}

In this example, we examine the atomic arrangement of amorphous silica
using the layered presentation of connected persistence diagrams.
Our dataset is a point cloud representing silicon and oxygen atoms in silica \cite{le2010isaacs},
as depicted in Figure \ref{fig:silica_input}.
The ring formations in this structure \cite{hiraoka2016hierarchical} motivate us to apply the constant thinning model from Section \ref{subsec:model-point-cloud}.
Specifically, we selectively remove atoms to perturb the structured arrangement
and then employ connected persistence diagrams to evaluate the structural changes. 
We work in homology dimension one throughout this example.

Let $M$ be a persistence module of $\cln$ obtained from a thinning model.
Recall that a connected persistence diagram directly visualizes
$\widetilde{\delta}_M$ (where the superscript for the essential assignment is omitted).
Each layered presentation in Figure \ref{fig:silica} is plotted following the rules below:
\begin{itemize}
    \item $\widetilde{\delta}_M\big|_{\Id}$: corresponds to the standard persistence diagram of the post-thinning point cloud, 
    visualized with the middle color bar where color saturation indicates multiplicity values.
    \item $\widetilde{\delta}_M\big|_{\Iu}$: corresponds to the standard persistence diagram of the original silica molecules, visualized with the left color bar.
    \item $\widetilde{\delta}_M\big|_{\Iud}$: corresponds to vertical persistence between generators, 
    visualized with the right color bar,
    with dashed lines implying negative multiplicity and the color saturation signifying the absolute value of multiplicity.
\end{itemize}
Although generators near the diagonal line overlap, the main features in the diagrams remain clear.

\noindent\textbf{Case 1: Removing All Silicon Atoms}

This scenario assesses topological changes when all silicon atoms are removed while oxygen atoms are left intact. 
The most prominent features before and after the thinning are labeled as 
Feature 1 and 2 in Figure \ref{fig:example-silica-case1}.
These two features appear to persist through the thinning process,
but that cannot be justified using the standard persistence diagrams alone.
This observation is confirmed here by the numerous horizontal green connecting lines 
between the two features. Therefore, we conclude that the main feature's birth is delayed, with its basic structure maintained.

\noindent\textbf{Case 2: Thinning Out Oxygen Atoms at a 50\% Rate}

Half of the oxygen atoms are removed in this case, while silicon atoms are untouched.
This amounts to a similar number of atoms being removed compared with Case 1.
Contrary to the first case, the primary topological feature shows significant deformation,
as can be seen by many more line segments emanating from Feature 1 but not terminating within Feature 3.
For those line segments connecting Feature 1 and Feature 3, 
 we observe that the transition distance is much longer, and the angles are steeper, reflecting delayed birth and death,
 exhibiting a more substantial structural disruption.

This comparative analysis demonstrates that
the connected persistence diagrams effectively reveal 
how different thinning strategies affect the structure.

\subsubsection{Face-Centered Cubic and Hexagonal Close Packing}\label{subsec:FCC-HCP}
Face-centered cubic (FCC) and 
hexagonal close packing (HCP) 
are two packings of equal spheres in three-dimensional space,
seen in various materials.
These two packings share many common properties and cannot be distinguished using standard persistence diagrams in any dimension.
This challenge has sparked numerous research efforts within the community.
Hiraoka et al. \cite{yasuaki_ippei_kazuto_2019} demonstrated that the persistence diagrams of the two structures after a thinning process are topologically distinct.
Meanwhile, Osang et al. \cite{osang2021topological} proposed $k$-fold covers, which modify the growing-radius ball model 
typically used to obtain a filtration of simplicial complexes from a given point cloud, 
to distinguish these two structures. 

\begin{figure}[htbp]
    \centering
    \begin{subfigure}{0.3\textwidth}
    \centering
    \begin{tikzpicture}
    \def\circlesize{1.3cm}  
    \def\circledist{1.3}  
    \pgfmathsetmacro{\sqrtthree}{sqrt(3)}
    \foreach \i in {0,...,2} 
        \node[draw, circle, minimum size=\circlesize, fill=area-color] at ({\i*\circledist+0.5*\circledist},{0*\sqrtthree/2*\circledist}) {};
        
    \foreach \i in {0,...,3} 
        \node[draw, circle, minimum size=\circlesize, fill=area-color] at ({\i*\circledist},{1*\sqrtthree/2*\circledist}) {};
        
    \foreach \i in {0,...,2} 
        \node[draw, circle, minimum size=\circlesize, fill=area-color] at ({\i*\circledist+0.5*\circledist},{2*\sqrtthree/2*\circledist}) {};
    \node at (1.5*\circledist, 0) {$\bullet$};
    \node[left,font=\large] at (1.5*\circledist, 0) {A};
    \node at (1.5*\circledist, 2/3*\sqrtthree/2*\circledist) {$\bullet$};
    \node[above left,font=\large] at (0.95*1.5*\circledist, 2/3*\sqrtthree/2*\circledist) {B};
    \node at (1.5*\circledist, 0.99*4/3*\sqrtthree/2*\circledist) {$\bullet$};
    \node[above,font=\large] at (1.5*\circledist, 1.05*4/3*\sqrtthree/2*\circledist) {C};
    \end{tikzpicture}
    \vspace*{\fill}
    \caption{}
    \label{fig:close-packing-layer}
    \end{subfigure}
    \hspace*{\fill}
    \begin{subfigure}{0.3\textwidth}
    \centering
    \includegraphics[height=4cm,keepaspectratio]{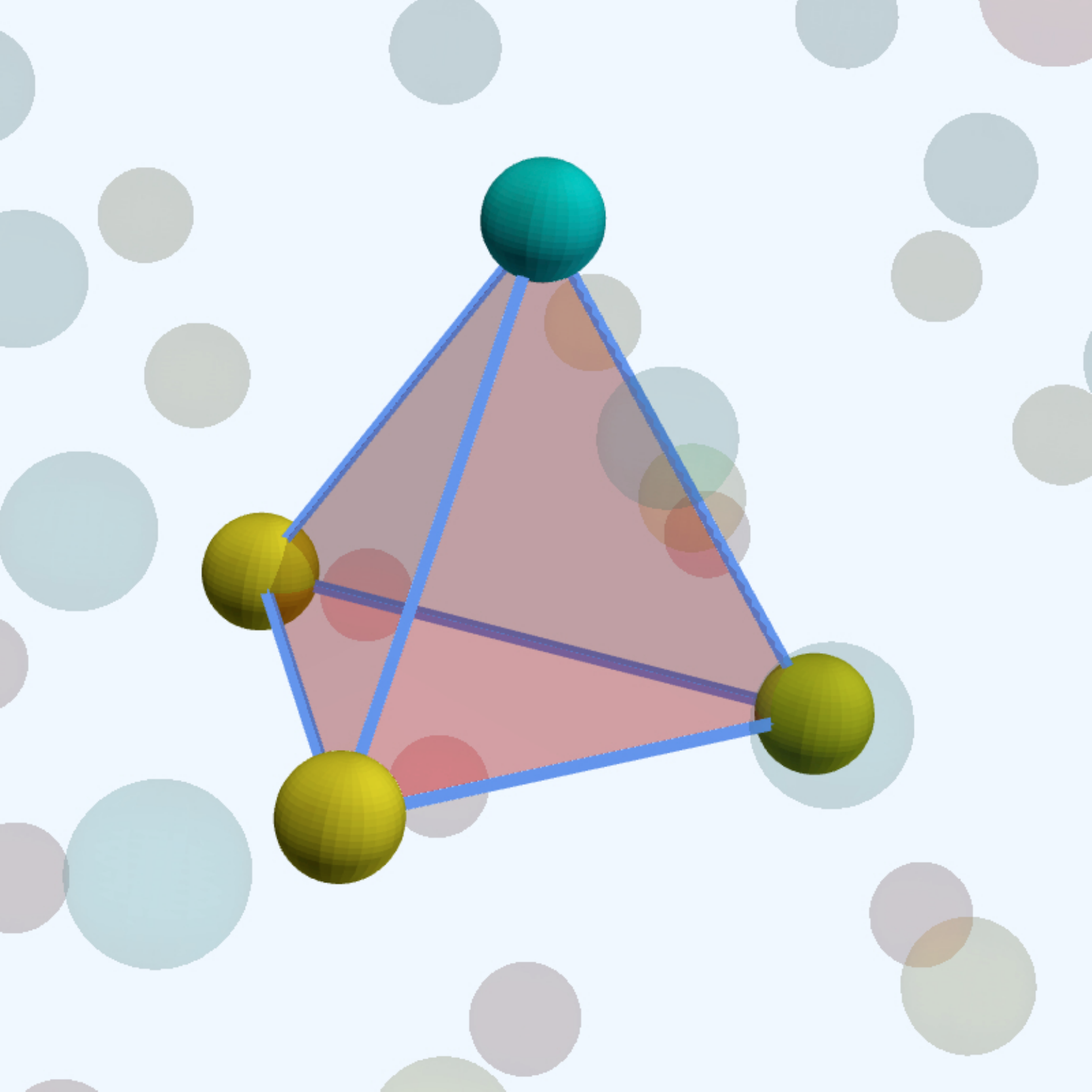}
    \caption{}\label{fig:fcc-hcp-tetrahedron}
    \end{subfigure}
    \hspace*{\fill}
    \begin{subfigure}{0.3\textwidth}
    \centering
    \includegraphics[height=4cm,keepaspectratio]{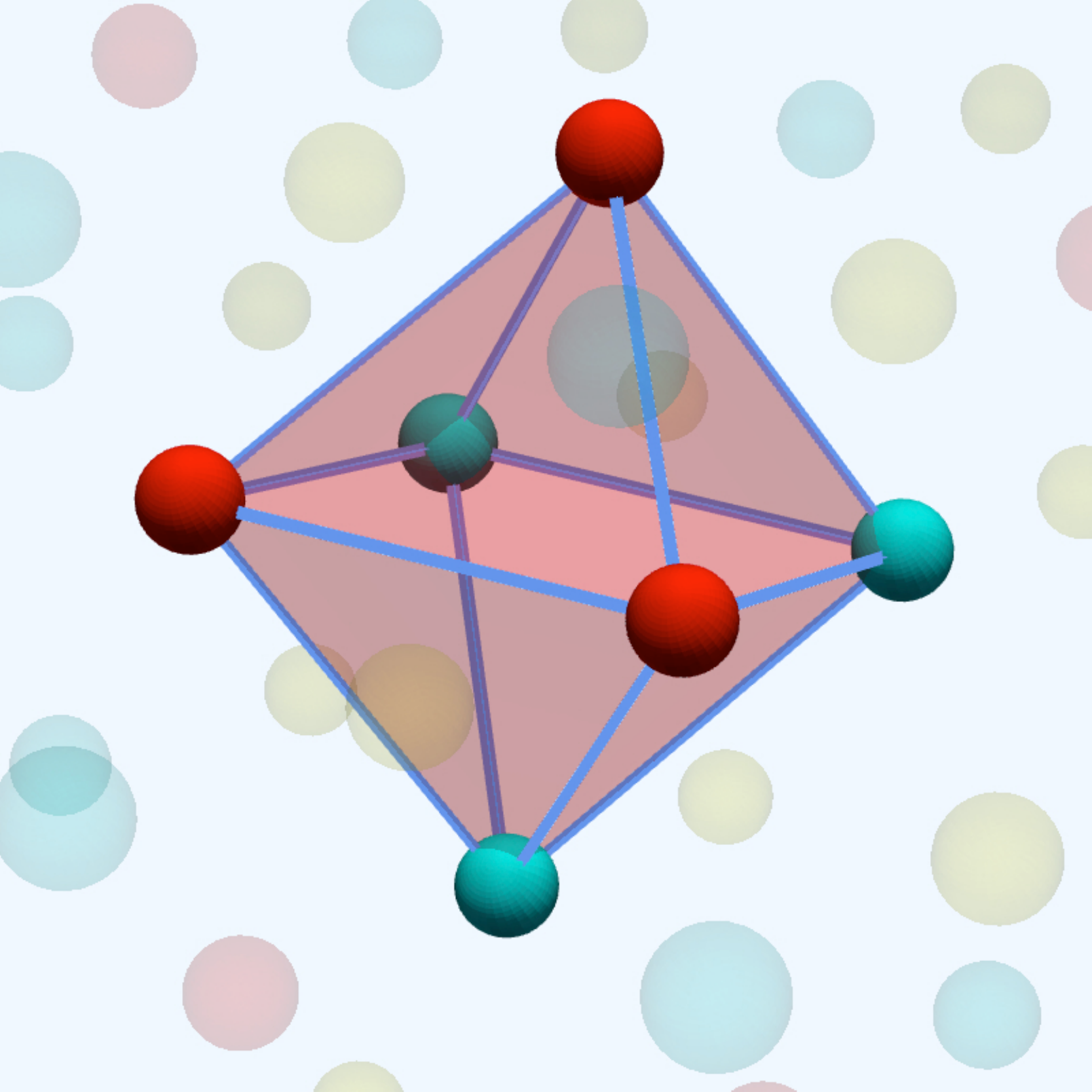}
    \caption{}\label{fig:fcc-hcp-octahedron}
    \end{subfigure}
    \caption{(a) Single layer illustration of a close packing. 
    (b) A two-dimensional homology generator in FCC with spheres colored by layer. Its HCP counterpart shares the same shape. 
    (c) Another type of homology generator in FCC. Notice that only two adjacent layers are involved, so we also have this in HCP.}
    \label{fig:fcc-hcp}
\end{figure}
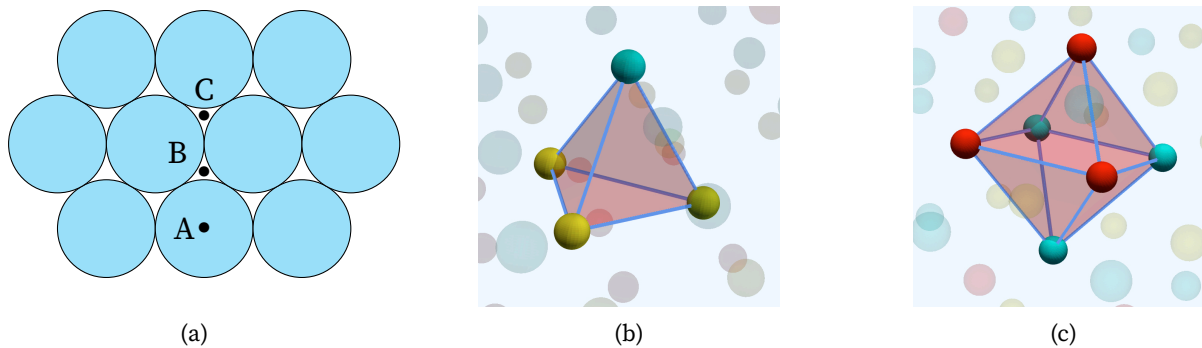

Figure \ref{fig:close-packing-layer} illustrates one layer of a packing.
If this layer serves as the base 2D plane, the relative position of any succeeding layer 
is determined by projecting the center of any sphere from that layer onto this base.
We label potential projection points as A, B, and C.
If the next layer projects at B, then the third layer could be A or C.
Maintaining a periodic packing pattern up until now across all layers, 
a choice of C for the third layer results in a layer sequence ABCABC$\ldots$, 
characteristic of FCC packing.
Conversely, choosing position A for the third layer yields a layer sequence ABABAB$\ldots$, referred to as HCP. 
For the subsequent discussion, we assume spheres with a radius of one and a homology dimension of two.

\begin{figure}[htbp]
    \centering
    \begin{subfigure}{0.24\textwidth}
        \centering
        \includegraphics[width=\linewidth, keepaspectratio]{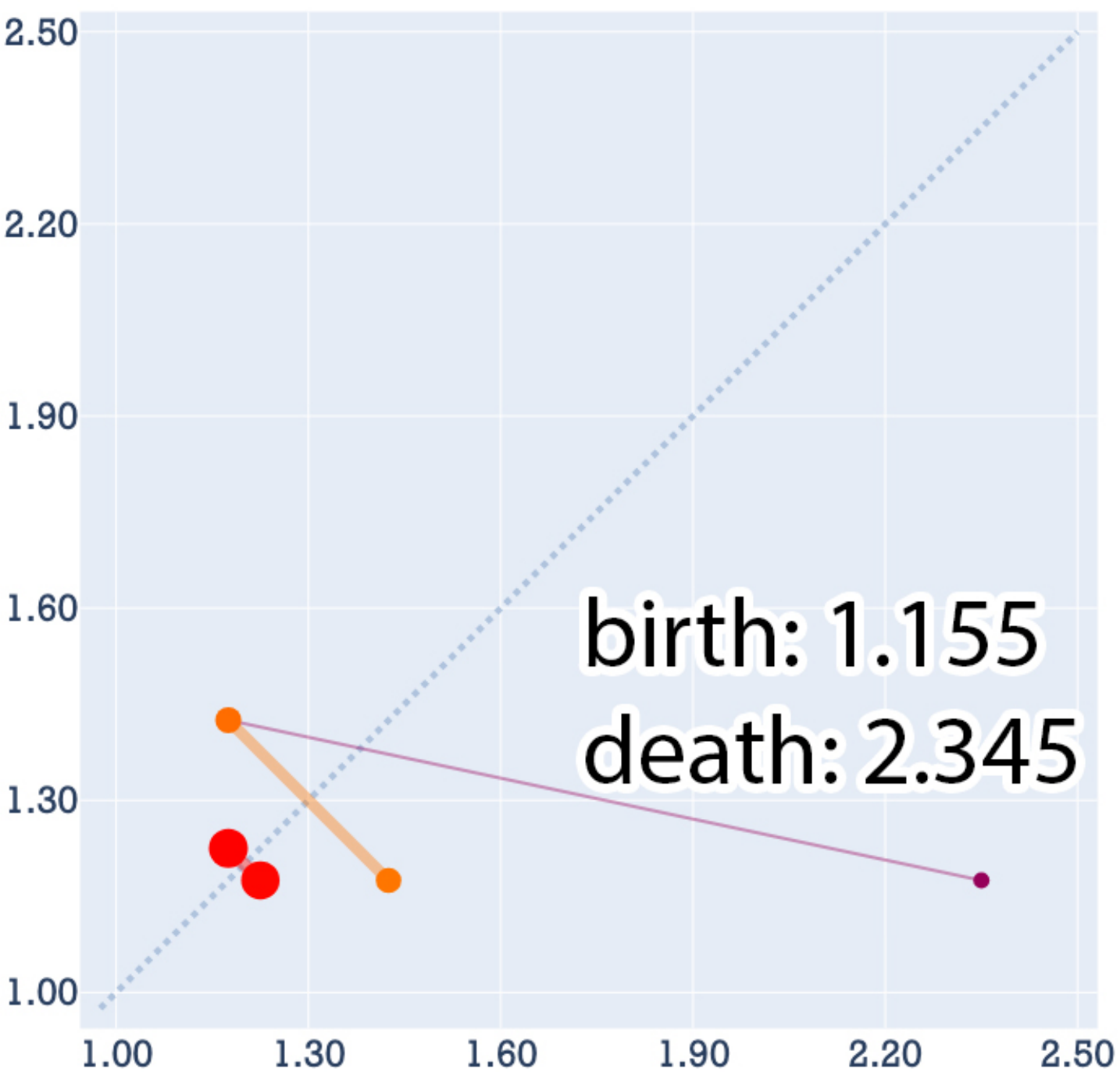}
        \caption{}
        \label{fig:tetrahedron-thinning-fcc-cPD}
    \end{subfigure}
    \hfill
    \begin{subfigure}{0.24\textwidth}
        \centering
        \includegraphics[width=\linewidth, keepaspectratio]{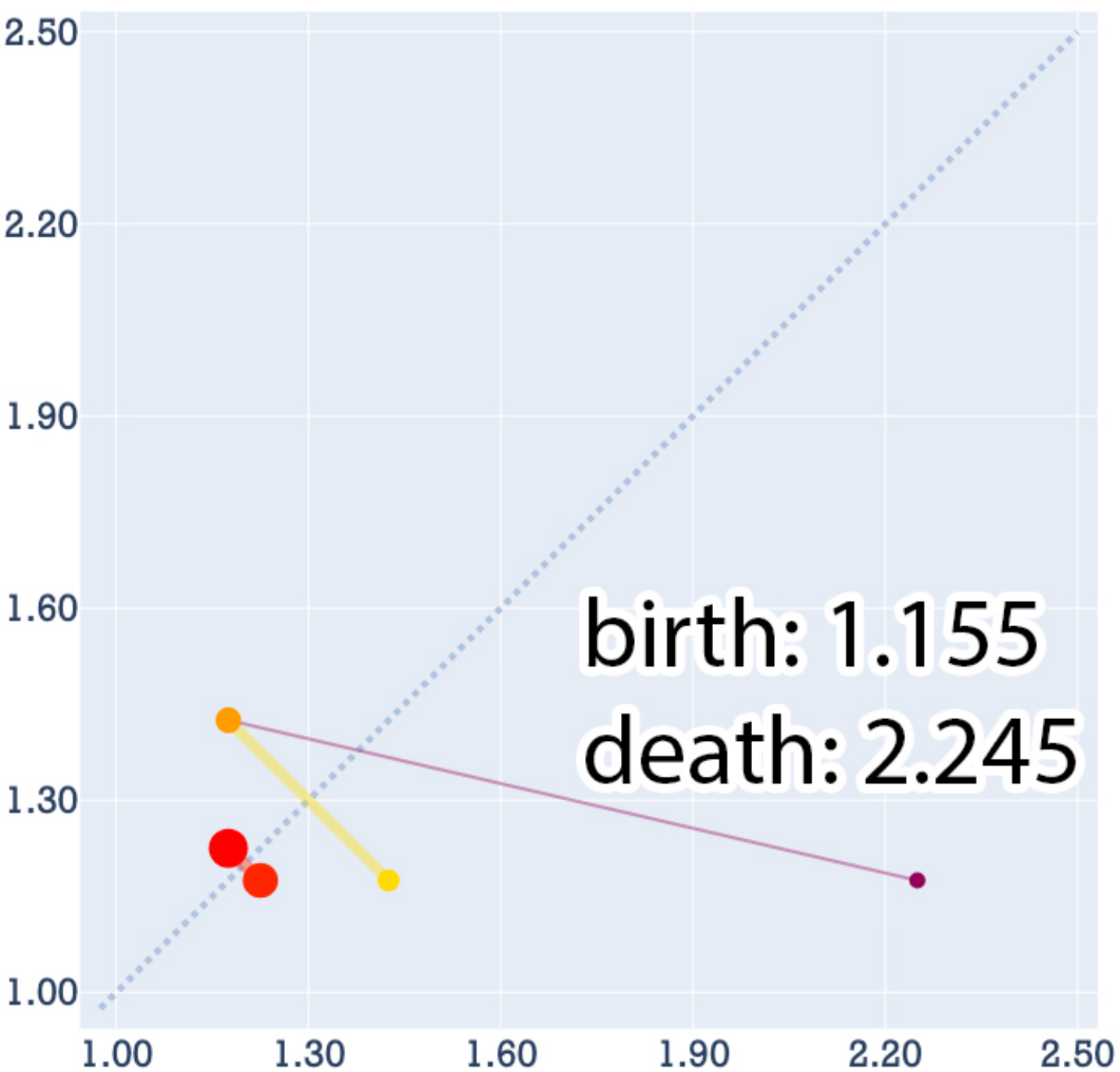}
        \caption{}
        \label{fig:tetrahedron-thinning-hcp-cPD}
    \end{subfigure}
    \hfill
    \begin{subfigure}{0.24\textwidth}
        \centering
        \includegraphics[trim={100pt 90pt 100pt 110pt},clip,width=\linewidth, keepaspectratio]{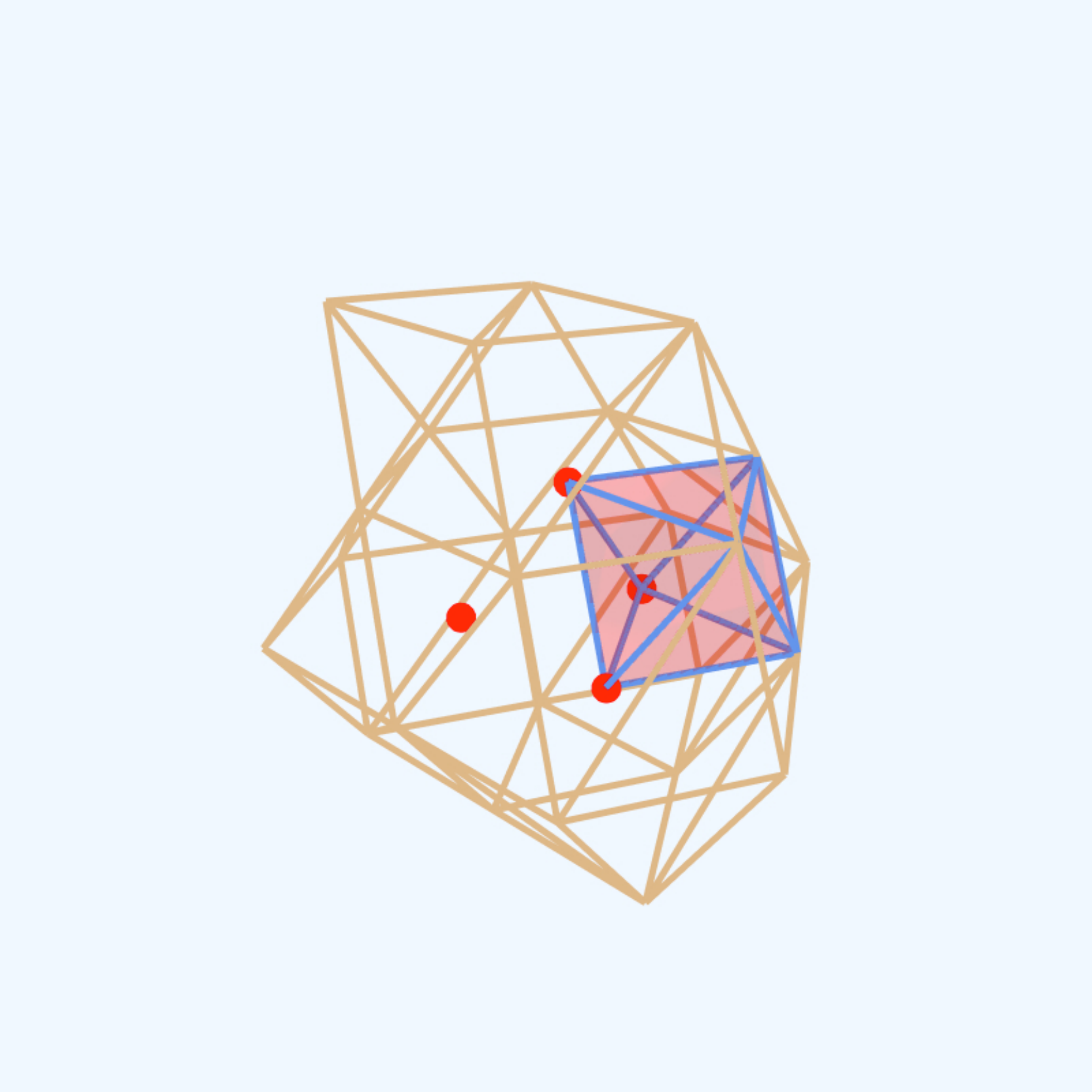}
        \caption{}
        \label{fig:tetrahedron-thinning-fcc-inv}
    \end{subfigure}
    \hfill
    \begin{subfigure}{0.24\textwidth}
        \centering
        \includegraphics[trim={90pt 81pt 90pt 99pt},clip,width=\linewidth, keepaspectratio]{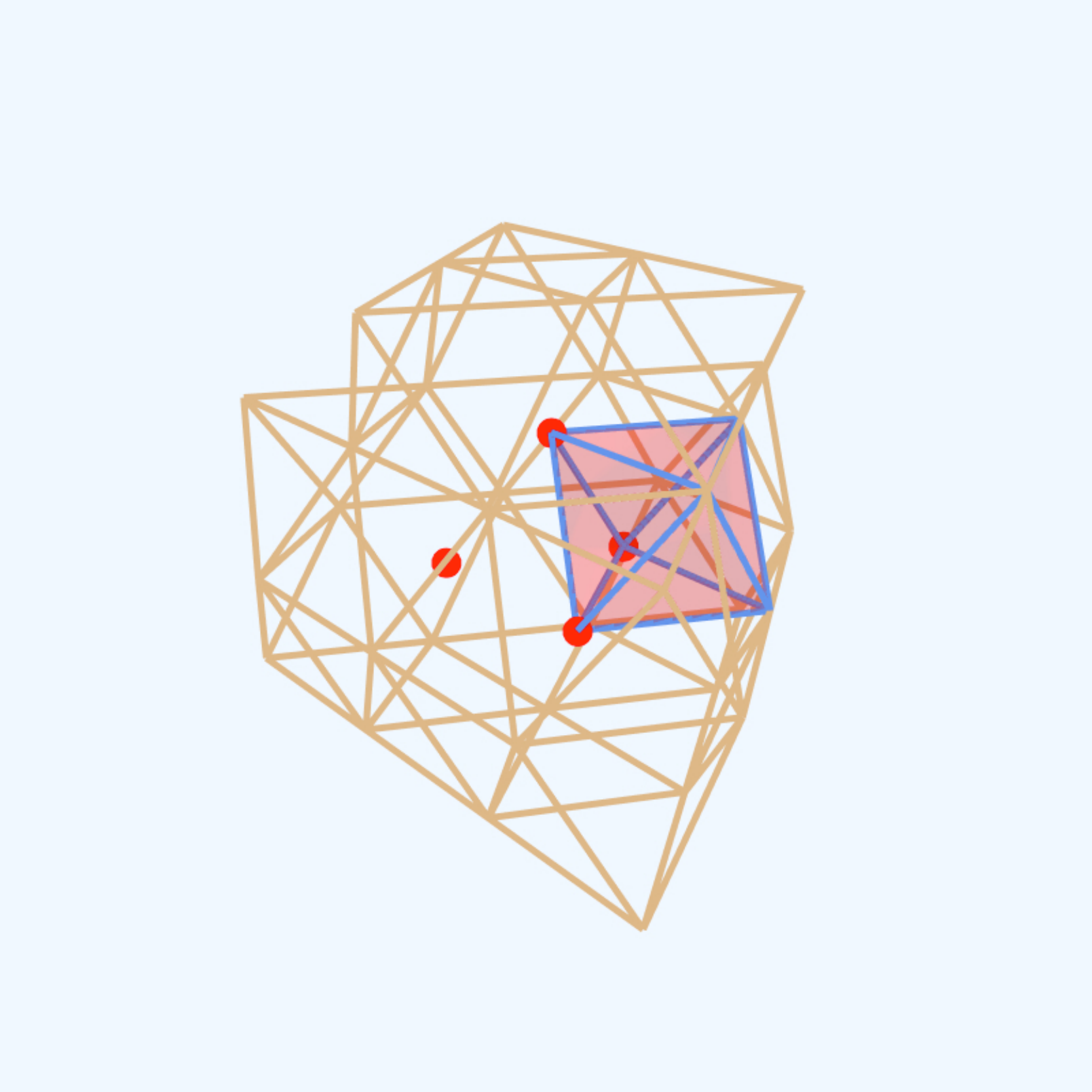}
        \caption{}
        \label{fig:tetrahedron-thinning-hcp-inv}
    \end{subfigure}
    \caption{
    (a) The connected persistence diagram for the FCC configuration.
    (b) The connected persistence diagram for the HCP configuration.
    (c) The new generator post-thinning in FCC. 
    (d) The new generator post-thinning in HCP. 
    }
    \label{fig:fcc-hcp-thinning-results}
\end{figure}

The standard persistence diagrams of FCC and HCP 
are identical to the upper half triangular region of Figure \ref{fig:tetrahedron-thinning-fcc-cPD} and \ref{fig:tetrahedron-thinning-hcp-cPD} respectively.
The generator closer to the diagonal line arises from tetrahedron structures formed by four neighboring atoms, as illustrated in Figure \ref{fig:fcc-hcp-tetrahedron}.
Another generator corresponds to the octahedron shown in Figure \ref{fig:fcc-hcp-octahedron}.
We utilize a patterned thinning model as discussed in Section \ref{subsec:model-point-cloud} to discern the subtle topological differences between these two packings.
Specific to this context, we remove tetrahedral structures randomly
 to create a $\cln$-filtration, and we restrict to removing only one such structure to 
 make the explanation more intuitive.

Figure \ref{fig:fcc-hcp-thinning-results} shows the resulting connected persistence diagrams and the homology generators post-thinning obtained via inverse analysis \cite{obayashi2022persistent}.
Notice that the deaths of the two labeled generators in the diagrams are different
\footnote{We record the numerical values of birth and death obtained in the experiment, but their analytical values are also computable. 
The one for FCC is $(\frac{2\sqrt{3}}{3}, \frac{\sqrt{22}}{2})$, and for HCP it is $(\frac{2\sqrt{3}}{3},\frac{11\sqrt{6}}{12})$.}.
We represent each such generator as a brown cage in 
Figure \ref{fig:tetrahedron-thinning-fcc-inv} and \ref{fig:tetrahedron-thinning-hcp-inv}, 
where faces are omitted for visual clarity.
The red points within each cage represent the tetrahedron structure removed during the thinning,
and we depict how each large cage can contain one original octahedron generator once embedded back into the upper row in the $\cln$-filtration. 
This relation is reflected in the connected persistence diagram via the connecting line between the octahedron generator and the cage generator.


\section{Concluding Remarks}


Refining interval approximations has led to notable advancements and extended developments.
The novel framework has broadened the 
scope of defining invariants on a given interval via accommodating courses of diverse and valid shapes. 
This discovery reveals that numerous ranks can be defined on a given interval,
expanding the previously understood fact. 
Furthermore, the partial interval approximation offers a practical invariant for general commutative grids.

When applied to the commutative ladders, 
our framework soon leads to an efficient solution for the indecomposable decomposition of finite-type commutative ladders. 
The algorithm, explained in Section \ref{subsec:application-cl4}, can potentially be instrumental in computing the decomposition of broader cases beyond its current application.
This indicates a vast horizon of unexplored courses, 
hinting at the potential discovery of more intricate non-interval structures.
For general commutative ladders, we introduced the connected persistence diagrams. 
These diagrams allow for the simultaneous visualization of both horizontal and vertical persistence, 
bridging a crucial gap in the field.
Based on the new toolkit and proposed models for the construction of commutative ladder filtrations,
we provided insights that stand among the early works analyzing the behavior of non-intervals,
and employed connected persistence diagrams to unveil hidden topological information in material structures.



\section*{Acknowledgements}
The authors would like to thank Prof. Asashiba for valuable discussions and Prof. Ochiai for his helpful feedback.
Y.H. was supported by JSPS Grant-in-Aid for Transformative Research Areas (A) (22H05107),
Grant-in-Aid for Scientific Research (A) (JP20H00119), 
and 
JST MIRAI Program (JPMJMI22G1).
K.N. was supported by JSPS Grant-in-Aid for Transformative Research Areas (A) (20H05884)
and 
JSPS KAKENHI JP (19H00834).
I.O. was supported by 
JSPS Grants-in-Aid for Transformative Research Areas (A) (20H05884),
JSPS KAKENHI JP (19H00834),
and JST PRESTO (JPMJPR1923).
C.X. was supported by JST SPRING (JPMJSP2110)
and RIKEN Junior Research Associate Program.

\section*{Statements and Declarations}
The authors declare no competing financial or non-financial interests related to this work.

\appendix
\section{Appendix}

\subsection{Enumeration Algorithm for Finding Alternating Zigzag Courses}\label{appendix:enum-zigzag-courses}

Algorithm \ref{alg:zigzag-courses-enum} 
is an enumeration algorithm 
that illustrates the procedures for discovering alternating zigzag courses. 
This approach involves examining all type $\mathbb{A}_n$ courses with 
an increasing number of paths up to a predetermined threshold value $N$. 

\newcommand\mycommfont[1]{\footnotesize\ttfamily\textcolor{blue}{#1}}
\SetCommentSty{mycommfont}
\begin{algorithm}[htbp]
\DontPrintSemicolon
\caption{Searching for all type $\mathbb{A}_n$ alternating zigzag courses with $n\leq N$ (via enumeration)}\label{alg:zigzag-courses-enum}
\KwData{a fully commutative quiver $G$}
\KwResult{a list of all type $\mathbb{A}_n$ alternating zigzag courses in $G$ with $n\leq N$}
$P \gets \text{the set of all non-trivial paths in } G$ \tcp*[r]{$P$ is finite as $G$ is finite and acyclic}
$J \gets \set{\big(A_1=\underset{1}{\bullet},F\colon 1\mapsto v\big)| v\in G_0}$ \tcp*[r]{initialize with type $\mathbb{A}_1$ courses}
\For{$n$ in $\set{2,\ldots,N}$}{
    \For{$(p_1,\ldots,p_{n-1})$ in $\underbrace{P\times P\times\cdots\times P}_{n-1}$}{
    \tcp{iterate over all $(n-1)$-tuples of paths in $P$}
        \If{$(p_1,\ldots, p_{n-1})$ determines an alternating zigzag course $(A_n,F)$ \tcp*[r]{see Remark \ref{rmk:zigzag-course-from-paths}}}{
            add $(A_n,F)$ to $J$\;
        }
    }
}
\Return {$J$}
\end{algorithm}

\subsection{BFS Algorithm for Finding Alternating Zigzag Courses}\label{appendix:bfs-zigzag-courses}
Algorithm \ref{alg:zigzag-courses-bfs}
provides a more efficient algorithm to find alternating zigzag courses
using a breadth-first search (BFS) approach.


\begin{algorithm}[htbp]
\DontPrintSemicolon
\caption{Searching for all type $\mathbb{A}_n$ alternating zigzag courses with $n\leq N$ (via BFS)}\label{alg:zigzag-courses-bfs}
\KwData{a fully commutative quiver $G$}
\KwResult{a list of all type $\mathbb{A}_n$ alternating zigzag courses on $G$ for $n\leq N$}
$startsFrom \gets \text{a hash map}$\;
$endsAt \gets \text{a hash map}$\;
\For{$v$ in $G_0$}{
$startsFrom(v) \gets \text{the set of non-trivial paths in $G$ that start from $v$}$\; 
$endsAt(v) \gets \text{the set of non-trivial paths in $G$ that end at $v$}$\;
}
$J \gets \set{\big(A_1=\underset{1}{\bullet},F\colon 1\mapsto v\big)| v\in G_0} $ \tcp*[r]{initialize with type $\mathbb{A}_1$ zigzag courses}
$queue \gets \set{\big(A_1=\underset{1}{\bullet},F\colon 1\mapsto v\big)| v\in G_0} $ \tcp*[r]{initialize the queue for BFS}
\While {$queue$ is not empty}{
$ currentCourse \gets queue.popleft()$ \;
$ n \gets \text{$n$ as the value in $currentCourse$ being a type $\mathbb{A}_n$ course}$\;
\If{$n+1>N$}{
\textbf{break} \tcp*[r]{break loop since all required alternating zigzag courses have been processed}
}
$ lastVertex\gets \text{the last vertex in $currentCourse$}$ \tcp*[r]{this vertex can be obtained as $F(n)$}
\eIf{$n\mod 2 == 1$}{
    $pathsToBeAttached=startsFrom(lastVertex)$ \tcp*[r]{attach rightward paths}
}{
    $pathsToBeAttached=endsAt(lastVertex)$ \tcp*[r]{attach leftward paths}
}
\For{$newPath$ in $pathsToBeAttached$}{
    $newCourse\gets \text{ } appendTo(currentCourse,newPath)$ \tcp*[r]{form a longer alternating zigzag course by appending $newPath$ to $currentCourse$}
    add $newCourse$ to $J$\; 
    $queue.append(newCourse)$\; 
}
}
\Return {$J$}
\end{algorithm}

 



\subsection{Algorithm for Extracting Linearly Independent Functions}\label{appendix:linearly-indep-picking}
Algorithm \ref{alg:coefficient-matrix} obtains a linearly independent set from the associated functions
of a set of alternating zigzag courses.
When a function is evaluated on the pre-determined set $\mathcal{L}$, it provides a new row for our coefficient matrix.
Each entry in this row is the multiplicity of the longest interval in a zigzag persistence module, which can be calculated using existing software packages such as \cite{dey2022fast} and \cite{Dionysus2}.
This algorithm then iteratively appends rows to the coefficient matrix until the matrix's rank reaches $|\mathcal{L}|$ or all potential candidates have been considered.

\begin{algorithm}[htbp]
\DontPrintSemicolon
\caption{Extracting linearly independent functions from a set of alternating zigzag courses}\label{alg:coefficient-matrix}
\KwData{a fully commutative quiver $G$; a finite subset $\mathcal{L}$ of the isomorphism classes; a base field $\Bbbk$ for the path algebra; a set $J$ of alternating zigzag courses}
\KwResult{a rank $|\mathcal{L}|$ coefficient matrix $C$ if it exists}
$C \gets \text{an empty matrix}$ \tcp*[r]{initialize the coefficient matrix}
$r \gets 0$ \tcp*[r]{for keeping track of the rank}
\For{$(A_n,F)$ in $J$}{
    $f(-)\gets d_{\tour_{(A_n, F)}(-)}(V_{A_n})$ \tcp*[r]{multiplicity of the longest interval in $\tour_{(A_n, F)}(-)$. }
    $C_\text{new}\gets\left( \begin{array}{c}
                    C \\ 
                    \hline 
                    \\[-1em]
                    \{f(L)\}_{L\in\mathcal{L}} \\
                    \end{array} 
                    \right)$ \tcp*[r]{add a new row to $C$ by incorporating the new function}
            $r_\text{new} \gets \rank(C_\text{new})$\;
            \If{$r_\text{new}>r$ \tcp*[r]{update the coefficient matrix only if the new function increases the rank}}{
                $C\gets C_\text{new}$\;
                $r \gets r_\text{new}$\;
            }
            \If{$r==|\mathcal{L}|$}{
                \Return{$C$}
            }
}
\textbf{raise} Exception(``Insufficient courses in the input to solve a general $\mathcal{L}$-decomposition.'')\;
\end{algorithm}

\FloatBarrier
\addcontentsline{toc}{section}{References}
\printbibliography[heading=bibnonumbered, title=References]

\end{document}